\documentclass{article} % For LaTeX2e
\usepackage{iclr2025_conference,times}

\usepackage{hyperref}
\usepackage{url}

\definecolor{mygreen}{rgb}{0,0.45,0}
\definecolor{myred}{rgb}{0.8,0,0}
\definecolor{maroon}{cmyk}{0,0.87,0.68,0.32}
\definecolor{azure}{rgb}{0.2, 0.5, 1.0}
\definecolor{dkgreen}{rgb}{0.0, 0.40, 0.10}
\definecolor{airforceblue}{rgb}{0.36, 0.54, 0.66}
\newcommand{\FEDAVG}{\dkgreen{\text{FedAvg}}}
\newcommand{\FEDPROX}{\dkgreen{\text{FedProx}}}
\newcommand{\GD}{\dkgreen{\text{GD}}}
\newcommand{\SGD}{\dkgreen{\text{SGD}}}
\newcommand{\PPM}{\dkgreen{\text{PPM}}}
\newcommand{\SPPM}{\dkgreen{\text{SPPM}}}
\newcommand{\FEDEXP}{\dkgreen{\text{FedExP}}}
\newcommand{\FEDEXPROX}{\dkgreen{\text{FedExProx}}}

\newcommand{\CGD}{\dkgreen{\text{CGD}}}
\newcommand{\AGD}{\dkgreen{\text{AGD}}}
\newcommand{\DCGD}{\dkgreen{\text{DCGD}}}
\newcommand{\dkgreen}[1]{{\color{dkgreen} #1}}

\newcommand{\rbrac}[1]{\left(#1\right)}
\newcommand{\sbrac}[1]{\left[#1\right]}
\newcommand{\cbrac}[1]{\left\{#1\right\}}
\newcommand{\inner}[2]{\left\langle #1, #2 \right\rangle}
\newcommand{\ProxSub}[2]{{\rm prox}_{#1}\left(#2\right)}

\newcommand{\MoreauSub}[3]{{M}^{#1}_{#2}\rbrac{#3}}
\newcommand{\Moreau}{{M}}
\newcommand{\M}[1]{{M}^{\gamma}\rbrac{#1}}

\newcommand{\bregman}[3]{{D}_{#1}\rbrac{#2, #3}}
\newcommand{\tx}{\Tilde{x}}
\newcommand{\aloc}[3]{A^{#1}_{#2}\rbrac{#3}}

\usepackage{amsmath,amsfonts,bm}

\def\eqref#1{equation~\ref{#1}}
\def\1{\bm{1}}

\def\mA{{\bm{A}}}

\def\mI{{\bm{I}}}

\def\mO{{\bm{O}}}

\DeclareMathAlphabet{\mathsfit}{\encodingdefault}{\sfdefault}{m}{sl}
\SetMathAlphabet{\mathsfit}{bold}{\encodingdefault}{\sfdefault}{bx}{n}

\newcommand{\R}{\mathbb{R}}

\usepackage{multirow}
\usepackage{calligra}
\usepackage{layout}

\usepackage{amsmath}
\usepackage{amsthm}
\usepackage{amssymb}
\usepackage{amsfonts}
\usepackage{mathtools}

\usepackage{url}
\usepackage{color}
\usepackage{graphicx}
\usepackage{algorithmic}
\usepackage{algorithm}
\usepackage{verbatim}
\usepackage{thmtools} 
\usepackage{thm-restate}
\usepackage{xcolor,colortbl}
\definecolor{midnightblue}{HTML}{0059b3}
\definecolor{darkmidnightblue}{HTML}{154c84}
\definecolor{noonblue}{HTML}{e5eef7}
\definecolor{chromered}{HTML}{f14233}
\definecolor{darkgreen}{HTML}{0e6029}

\usepackage{nicefrac}
\usepackage{adjustbox}

\usepackage{cleveref}

\newcommand{\norm}[1]{\left\| #1 \right\|}

\newcommand{\bB}{\mathbb{B}}

\newcommand{\bbS}{\mathbb{S}}

\newcommand{\cC}{\mathcal{C}}
\newcommand{\cD}{\mathcal{D}}
\newcommand{\cE}{\mathcal{E}}
\newcommand{\cF}{\mathcal{F}}

\newcommand{\cO}{\mathcal{O}}

\newcommand{\cX}{\mathcal{X}}

\newcommand{\del}[1]{}

\newcommand{\eqdef}{:=}

\newcommand{\Exp}[1]{{\mathbb E}\left[#1\right]}

\newcommand{\ExpSub}[2]{{\mathbb E}_{#1}\left[#2\right]}
\newtheorem{assumption}{Assumption}
\newtheorem{lemma}{Lemma}

\newtheorem{theorem}{Theorem}

\newtheorem{definition}{Definition}
\newtheorem{fact}{Fact}

\theoremstyle{definition}

\usepackage[colorinlistoftodos,bordercolor=orange,backgroundcolor=orange!20,linecolor=orange,textsize=scriptsize]{todonotes}

\usepackage{makecell}
\usepackage{booktabs}
\usepackage{enumitem}
\usepackage{subfigure}
\usepackage{threeparttable}

\title{On the Convergence of FedProx with Extrapolation and Inexact Prox}

\author{Hanmin Li \\
Generative AI CoE \\
KAUST \\
Thuwal, Saudi Arabia \\
\texttt{hanmin.li@kaust.edu.sa}
\And
Peter Richt\'{a}rik \\
Generative AI CoE \\
KAUST \\
Thuwal, Saudi Arabia \\
\texttt{peter.richtarik@kaust.edu.sa}
}

\iclrfinalcopy % Uncomment for camera-ready version, but NOT for submission.
\begin{document}

\maketitle

\begin{abstract}
   Enhancing the FedProx federated learning algorithm \citep{li2020federated} with server-side extrapolation, \citet{li2024power} recently introduced the {\FEDEXPROX} method.
   Their theoretical analysis, however, relies on the assumption that each client computes a certain proximal operator exactly, which is impractical since this is virtually never possible to do in real settings.
   In this paper, we investigate the behavior of {\FEDEXPROX} without this exactness assumption in the smooth and globally strongly convex setting.
   We establish a general convergence result, showing that inexactness leads to convergence to a neighborhood of the solution.
   Additionally, we demonstrate that, with careful control, the adverse effects of this inexactness can be mitigated.
   By linking inexactness to biased compression \citep{beznosikov2023biased}, we refine our analysis, highlighting robustness of extrapolation to inexact proximal updates.
   We also examine the local iteration complexity required by each client to achieved the required level of inexactness using various local optimizers.
   Our theoretical insights are validated through comprehensive numerical experiments.
\end{abstract}

\section{Introduction}
Distributed optimization is becoming increasingly essential in modern machine learning, especially as models grow more complex.
Federated learning (FL), a decentralized approach where multiple clients collaboratively train a shared model while keeping their data locally to preserve privacy, is a key example of this trend \citep{Konecny2016federated,mcmahan2017communication}.
Often, a central server coordinates the process by aggregating the locally trained models from each client to update the global model without accessing the raw data.
The federated average algorithm ({\FEDAVG}), introduced by \citet{mcmahan2017communication} and \citet{mangasarian1993backpropagation}, is one of the most popular strategies for tackling federated learning problems.
The algorithm comprises three essential components: client sampling, data sampling, and local training.
During its execution, the server first samples a subset of clients to participate in the training process for a given round.
Each selected client then performs local training using stochastic gradient descent ({\SGD}), with or without random reshuffling, to enhance communication efficiency, as documented by \citet{bubeck2015convex, gower2019sgd, moulines2011non,sadiev2022federated}.
{\FEDAVG} has proven to be highly successful in practice, nevertheless it suffers from client drift when data is heterogeneous \citep{karimireddy2020scaffold}.

Various techniques have been proposed to address the challenges of data heterogeneity, with {\FEDPROX}, introduced by \citet{li2020federated}, being one notable example.
Rather than having each client perform local {\SGD} rounds, {\FEDPROX} requires each client to compute a proximal operator locally.
Computing the proximal operator can be regarded as an optimization problem that each client can solve locally.
Proximal algorithms are advantageous when the proximal operators can be evaluated relatively easily \citep{parikh2014proximal}.
Algorithms based on proximal operators, such as the proximal point method ({\PPM}) \citep{rockafellar1976monotone,parikh2014proximal} and its extension to the stochastic setting ({\SPPM}) \citep{bertsekas2011incremental,asi2019stochastic,khaled2022faster,richtarik2020stochastic,patrascu2018nonasymptotic}, offer greater stability against inaccurately specified step sizes, unlike gradient-based methods.
{\PPM} was introduced by \citet{martinet1972algorithmes} and expanded by \citet{rockafellar1976monotone}.
Its extension into the stochastic setting are often used in federated optimization.
The stability mentioned is particularly useful when problem-specific parameters, such as the smoothness constant of the objective function, are unknown which renders determining the step size for {\SGD} becomes challenging.
Indeed, an excessively large step size in {\SGD} leads to divergence, while a small step size ensures convergence but significantly slows down the training process.

Another approach to mitigating the slowdown caused by heterogeneity is the use of a server step size.
Specifically, in {\FEDAVG}, a local step size is employed by each client to minimize their individual objectives, while a server step size is used to aggregate the `pseudo-gradients' obtained from each client \citep{karimireddy2020scaffold,reddi2021adaptive}.
The local step size is set relatively small to mitigate client drift, while the server step size is set larger to avoid slowdowns.
However, the small step sizes result in a slowdown during the initial phase of training, which cannot be fully compensated by the large server step size \citep{jhunjhunwala2023fedexp}.
Building on the extrapolation technique employed in parallel projection methods to solve the convex feasibility problem \citep{censor2001averaging, combettes1997convex, necoara2019randomized}, \citet{jhunjhunwala2023fedexp} introduced {\FEDEXP} as an extension of {\FEDAVG}, incorporating adaptive extrapolation as the server step size.
Extrapolation involves moving further along the line connecting the most recent iterate, $x_k$, and the average of the projections of $x_k$ onto different convex sets, $\cX_i$, in the parallel projection method, which accelerates the algorithm. 
Extrapolation is also known as over-relaxation \citep{rechardson1911approximate} in fixed point theory. 
It is a common technique to effectively accelerate the convergence of fixed point methods including gradient based algorithms and proximal splitting algorithms \citep{condat2023proximal, iutzeler2019generic}. 
Recently, \citet{li2024power} shows that the combination of extrapolation with {\FEDPROX} also results in better complexity bounds.
The analysis of the resulting algorithm {\FEDEXPROX} reveals the relationship between the extrapolation parameter and the step size of gradient-based methods with respect to the Moreau envelope associated with the original objective function.
However, it relies on the assumption that each proximal operator is solved accurately, which makes it impractical and less advantageous compared to gradient-based algorithms.

\subsection{Contributions}

Our paper makes the following contributions, please refer to \Cref{sec:notations} for notation details.

\begin{itemize}[noitemsep, left=0pt]
   \item We provide a new analysis of {\FEDEXPROX} based on \citet{li2024power}, focusing on the case where the proximal operators are evaluated inexactly in the globally strongly convex setting, removing the need for the assumption of exact proximal operator evaluations.
   By properly defining the notion of approximation, we establish a general convergence guarantee of the algorithm to a neighborhood of the solution utilizing the theory of biased {\SGD} \citep{demidovich2024guide}.
   Specifically, our algorithm achieves a linear convergence rate of $\cO\rbrac{\frac{L_{\gamma}\rbrac{1 + \gamma L_{\max}}}{\mu}}$ to a neighborhood of the solution, matching the rate presented by \citet{li2024power}.

   \item Building on our understanding of how the neighborhood arises, we propose a new method of approximation.
   This alternative characterization of inexactness eliminates the neighborhood from the previous convergence guarantee, provided that the inexactness is properly bounded, and the extrapolation parameter is chosen to be sufficiently small.

   \item By leveraging the similarity between the definitions of inexactness and compression, we enhance our analysis using the theory of biased compression \citep{beznosikov2023biased}.
   The improved analysis offers a faster rate of $\cO\rbrac{\frac{L_{\gamma}\rbrac{1 + \gamma L_{\max}}}{\mu - 4\varepsilon_2 L_{\max}}}$\footnote{The parameter $\varepsilon_2$ is the parameter associated with accuracy of relative approximation as defined in \Cref{def:inexact-2}. We use the notation $\cO\rbrac{\cdot}$ to ignore constant factors and $\tilde{\cO}\rbrac{\cdot}$ when logarithmic factors are also omitted.}, leading to convergence to the exact solution, provided that the inexactness is bounded in a more permissive manner. 
   More importantly, the optimal extrapolation $\nicefrac{1}{\gamma L_{\gamma}}$ matches the exact case.
   This shows that extrapolation aids convergence as long as sufficient accuracy is reached, even with inexact proximal evaluations.

   \item We then analyze how the aforementioned approximations can be obtained by each client.
   As examples, we provide the local iteration complexity when the client employs gradient descent (\GD) or Nesterov's accelerated gradient descent (\AGD), demonstrating that these approximations are readily achievable.
   Specifically, for the $i$-th client, the local iteration complexity is $\tilde{\cO}\rbrac{1 + \gamma L_i}$ when using {\GD}, and $\tilde{\cO}\rbrac{\sqrt{1 + \gamma L_i}}$ when using {\AGD}.
   See \Cref{table-1-paper} and \Cref{table-2-paper} for a detailed comparison of complexities of all relevant quantities.

   \item Finally, we validate our theoretical findings through numerical experiments.
   Our numerical results suggest that the proposed technique of relative approximation effectively eliminates bias. 
   In some cases, the algorithm even outperforms {\FEDPROX} with exact updates, further validating the effectiveness of server extrapolation, even when proximal updates are inexact.
\end{itemize}

\begin{table}[t]

   \caption{
      Comparison of {\FEDEXPROX} \citep{li2024power} and our proposed inexact versions of the algorithms using different approximations. 
      In the convergence column, we present the rate at which each algorithm converges to either the solution or a neighborhood in the globally strongly convex setting. 
      Here, $L_{\gamma}$ represents the smoothness constant of $M^{\gamma}$ as defined before \Cref{thm:1:conv-full-batch-stncvx}. 
      The neighborhood column indicates the size of the neighborhood, while the optimal extrapolation column suggests the best choice of $\alpha$ for each algorithm. 
      The final column outlines the conditions on the inexactness.
      All quantities are presented with constant factors omitted, $K$ is the number of total iterations, $\gamma$ is the local step size for the proximal operator, $S\rbrac{\varepsilon_2}$ defined in \Cref{thm:020-full} is a factor of slowing down due to inexactness in $(0, 1]$.
      For relative approximation, we first present the original theory in the third row and then place the sharper analysis in the following row for comparison.
      }

   \label{table-1-paper}
   \begin{center}
       \begin{threeparttable}
           \footnotesize       
           \begin{tabular}{c!{\vrule width -2.5pt}cccc}
               \toprule 
               \multicolumn{1}{c}{{\bf Algorithm }} & \multicolumn{1}{c}{\bf Convergence} & \multicolumn{1}{c}{\bf Neighborhood} & \multicolumn{1}{c}{\bf \makecell{Optimal \\ Extrapolation}} & \multicolumn{1}{c}{\bf \makecell{Bound on \\ Inexactness}} \\ 
               \midrule \midrule \\
               {\FEDEXPROX} & $\exp\rbrac{-\frac{K\mu}{L_{\gamma}\rbrac{1 + \gamma L_{\max}}}}$ & $0$  & $\frac{1}{\gamma L_{\gamma}}$ & NA \\ \rowcolor{azure!10}
               \makecell{(NEW) {\FEDEXPROX} with \\ $\varepsilon_1$ approximation} & $\exp\rbrac{-\frac{K\mu}{L_{\gamma}\rbrac{1 + \gamma L_{\max}}}}$ & $\varepsilon_1\rbrac{\frac{\frac{1}{\gamma} + L_{\max}}{\mu}}^2$\tnote{(a)} & $\frac{1}{4\gamma L_{\gamma}}$ & NA \tnote{(b)} \\ \rowcolor{azure!15}
               \makecell{(NEW) {\FEDEXPROX} with \\ $\varepsilon_2$ relative  approximation \\ by biased {\SGD}} & $\exp\rbrac{-\frac{K\mu S\rbrac{\varepsilon_2}}{L_{\gamma}\rbrac{1 + \gamma L_{\max}}}}$\tnote{(c)} & $0$ & $< \frac{1}{\gamma L_{\gamma}}$ & $< \frac{\mu^2}{4L_{\max}^2}$ \\ \rowcolor{azure!20}
               \makecell{(NEW) {\FEDEXPROX} with \\ $\varepsilon_2$ relative  approximation \\ by biased compression} & $\exp\rbrac{-\frac{K\rbrac{\mu - 4\varepsilon_2 L_{\max}}}{L_{\gamma}\rbrac{1 + \gamma L_{\max}}}}$ & $0$ & $\frac{1}{\gamma L_{\gamma}}$\tnote{(d)} & $< \frac{\mu}{4L_{\max}}$ \\
               \bottomrule
           \end{tabular}
           \begin{tablenotes}
               \item[(a)]{\footnotesize 
                  Note that when $\varepsilon_1 = 0$, i.e., when the proximal operators are evaluated exactly, the neighborhood diminishes, and we recover the result of {\FEDEXPROX} by \citet{li2024power}, up to a constant factor.
               } 
               \item[(b)]{\footnotesize
                  Unlike relative approximations, the convergence guarantee here is more general, allowing for the analysis of unbounded inexactness. However, as the inexactness increases, the neighborhood grows correspondingly, rendering the result practically useless.
               }
               \item[(c)]{\footnotesize
                  Refer to \Cref{thm:020-full} for the definition of $S\rbrac{\varepsilon_2}$ and the corresponding optimal extrapolation parameter. The theory indicates that inexactness will adversely affect the algorithm's convergence.
               }
               \item[(d)]{\footnotesize
                  Surprisingly, our sharper analysis reveals that the optimal extrapolation parameter in this case remains the same as in the exact setting, highlighting the effectiveness of extrapolation even when the proximal operators are evaluated inexactly.
               }
           \end{tablenotes}
       \end{threeparttable}
   \end{center}
\end{table}

\begin{table}[t]

   \caption{
         Comparison of local iteration complexities of each client in order to obtain an approximation using either {\GD} or {\AGD} \citep{nesterov2013introductory}.
         We use the $i$-th client as an example, where the local objective $f_i: \R^d \mapsto \R$ is $L_i$-smooth and convex, $i \in \cbrac{1, 2, \hdots, n}$.
      }

   \label{table-2-paper}
   \begin{center}
       \begin{threeparttable}
           \footnotesize       
           \begin{tabular}{lcc}
               \toprule 
               \multicolumn{1}{l}{{\bf Algorithm }} & \multicolumn{1}{c}{\bf $\varepsilon_1$ absolute approximation} & \multicolumn{1}{c}{\bf $\varepsilon_2$ relative approximation} \\
               \midrule \midrule \\
               Gradient descent & $\cO\rbrac{\rbrac{1+ \gamma L_i}\log\rbrac{\frac{\norm{x_k - \ProxSub{\gamma f_i}{x_k}}^2}{\varepsilon_1}}}$ \tnote{(a)} & $\cO\rbrac{\rbrac{1+ \gamma L_i}\log\rbrac{\frac{1}{\varepsilon_2}}}$ \\
               Accelerate gradient descent & $\cO\rbrac{\sqrt{1 + \gamma L_i}\log\rbrac{\frac{\norm{x_k - \ProxSub{\gamma f_i}{x_k}}^2}{\varepsilon_1}}}$ & $\cO\rbrac{\sqrt{1 + \gamma L_i}\log\rbrac{\frac{1}{\varepsilon_2}}}$ \\
               \bottomrule
           \end{tabular}
           \begin{tablenotes}
               \item[(a)]{\footnotesize 
                  We can easily provide an upper bound of $\norm{x_k - \ProxSub{\gamma f_i}{x_k}}^2$ for determining the number of local computations needed.
               } 
           \end{tablenotes}
       \end{threeparttable}
   \end{center}
\end{table}

\subsection{Related work}
% \paragraph{{\color{black} Stochastic Gradient Descent.}} 
Arguably, stochastic gradient descent ({\SGD}) \citep{robbins1951stochastic, ghadimi2013stochastic,gower2019sgd,gorbunov2020unified} remains one of the foundational algorithm in the field of machine learning. 
One can simply formulate it as 
\begin{equation*}
   x_{k+1} = x_k - \eta \cdot g(x_k),
\end{equation*}
where $\eta > 0$ is a scalar step size, $g(x_k)$ is a possibly stochastic estimator of the true gradient $\nabla f(x_k)$.
In the case when $g(x_k) = \nabla f(x_k)$, {\SGD} becomes {\GD}.
Various extensions of {\SGD} have been proposed since its introduction, examples include compressed gradient descent ({\CGD}) \citep{alistarh2017qsgd,khirirat2018distributed}, {\SGD} with momentum \citep{loizou2017linearly, liu2020improved}, {\SGD} with matrix step size \citep{li2023det} and variance reduction \citep{gower2020variance, johnson2013accelerating,gorbunov2021marina,tyurin2024dasha,li2023marina}. 
\citet{gower2019sgd} presented a framework for analyzing {\SGD} with unbiased gradient estimator in the convex case based on expected smoothness. 
However, in practice, sometimes the gradient estimator could be biased, examples include {\SGD} with sparsified or delayed update \citep{alistarh2018convergence, recht2011hogwild}.
\citet{beznosikov2023biased} examined biased updates in the context of compressed gradient descent.
\citet{demidovich2024guide} provides a framework for analyzing {\SGD} with biased gradient estimators in the non-convex setting.

% \paragraph{{\color{black} Stochastic Proximal Point Method.}}
Proximal point method ({\PPM}) was originally introduced as a method to solve variational inequalities \citep{martinet1972algorithmes,rockafellar1976monotone}. 
The transition to the stochastic case, driven by the need to efficiently address large-scale optimization problems, leads to the development of {\SPPM}.
Due to its stability and advantage over the gradient based methods, it has been extensively studied, as documented by \citep{patrascu2018nonasymptotic,bianchi2016ergodic,bertsekas2011incremental}.  
For proximal algorithms to be practical, it is commonly assumed that the proximal operator can be solved efficiently, such as in cases where a closed-form solution is available.
However, in large-scale machine learning models, it is rarely possible to find such a solution in closed form.
To address this issue, most proximal algorithms assume that only an approximate solution is obtained, achieving a certain level of accuracy \citep{khaled2022faster, sadiev2022communication, karagulyan2024spam}.
Various notions of inexactness are employed, depending on the assumptions made, the properties of the objective, and the availability of algorithms capable of efficiently finding such approximations.

% \paragraph{{\color{black} Moreau Envelope.}}
Moreau envelope was first introduced to handle non-smooth functions by \citet{moreau1965proximite}. 
It is also known as the Moreau-Yosida regularization.
The use of the Moreau envelope as an analytical tool to analyze proximal algorithms is not novel.
\citet{ryu2014stochastic} noted that running a proximal algorithm on the objective is equivalent to applying gradient methods to its Moreau envelope. 
\citet{davis2019stochastic} analyzed stochastic proximal point method ({\SPPM}) for weakly convex and Lipschitz functions based on this finding.
Recently, \citet{li2024power} provided an analysis of {\FEDPROX} with server-side step size in the convex case, based on the reformulation of the problem using the Moreau envelope.
The role of the Moreau envelope extends beyond analyzing proximal algorithms; it has also been applied in the contexts of personalized federated learning \citep{t2020personalized} and meta-learning \citep{mishchenko2023convergence}.
The mathematical properties of the Moreau envelope are relatively well understood, as documented by 
\citet{jourani2014differential,planiden2019proximal,planiden2016strongly}.

Projection methods initially emerged as an effective tool for solving systems of linear equations or inequalities \citep{kaczmarz1993approximate} and were later generalized to solve the convex feasibility problem \citep{combettes1997convex}. 
The parallel version of this approach involves averaging the projections of the current iterates onto all existing convex sets $\cX_i$ to obtain the next iterate, a process that is empirically known to be accelerated by extrapolation.
Numerous heuristic rules have been proposed to adaptively set the extrapolation parameter, such as those by \citet{bauschke2006extrapolation} and \citet{pierra1984decomposition}.
Only recently, the mechanism behind constant extrapolation was uncovered by \citet{necoara2019randomized}, who developed the corresponding theoretical framework.
Additionally, \citet{li2024power} provides explanations for the effectiveness of adaptive rules, revealing the connection between the extrapolation parameter and the step size of {\SGD} when using the Moreau envelope as the global objective.

\section{Mathematical background}
In this work, we are interested in the distributed optimization problem which is formulated in the following finite-sum form
\begin{equation}
   \label{eq:prob-formulation}
   \min_{x \in \R^d} \cbrac{f(x) \eqdef \frac{1}{n}\sum_{i=1}^n f_i(x)},
\end{equation}
where $x \in \R^d$ is the model, $n$ is the number of devices/clients, $f: \R^d \mapsto \R$ is global objective, each $f_i: \R^d \mapsto \R$ is the empirical risk of model $x$ associated with the $i$-th client. 
Each $f_i\rbrac{x}$ often has the form 
\begin{equation}
   \label{eq:substructure}
   f_i(x) \eqdef \ExpSub{\xi \sim \cD_i}{l\rbrac{x, \xi}},
\end{equation}
where the loss function $l\rbrac{x, \xi}$ represents the loss of model $x$ on data point $\xi$ over the training data $\cD_i$ owned by client $i \in [n] \eqdef \cbrac{1,2,\hdots, n}$.
We first give the definitions for the proximal operator and Moreau envelope, which we will be using in our analysis.
\begin{definition}[Proximal operator]
   \label{def:1:prox}
   The proximal operator of an extended real-valued function $\phi:\R^d \mapsto \R\cup\cbrac{+\infty}$ with step size $\gamma > 0$ and center $x \in \R^d$ is defined as 
   \begin{align*}
      \ProxSub{\gamma \phi}{x} \eqdef \arg\min_{z\in \R^d}\cbrac{\phi\cbrac{z} + \frac{1}{2\gamma}\norm{z - x}^2}.
   \end{align*}
\end{definition}
It is well-known that for any proper, closed, and convex function $\phi$, the proximal operator with any $\gamma > 0$ returns a singleton.
\begin{definition}[Moreau envelope]
   \label{def:2:moreau}
   The Moreau envelope of an extended real-valued function $\phi:\R^d \mapsto \R\cup\cbrac{+\infty}$ with step size $\gamma > 0$ and center $x \in \R^d$ is defined as 
   \begin{align*}
      \MoreauSub{\gamma}{\phi}{x} \eqdef \min_{z \in \R^d}\cbrac{\phi\rbrac{z} + \frac{1}{2\gamma}\norm{z - x}^2}.
   \end{align*}
\end{definition}
By the definition of Moreau envelope, it is easy to see that 
\begin{align}
   \label{eq:ppt-1-func-val}
   \MoreauSub{\gamma}{\phi}{x} = \phi\rbrac{\ProxSub{\gamma \phi}{x}} + \frac{1}{2\gamma}\norm{x - \ProxSub{\gamma \phi}{x}}^2.
\end{align}
Not only are their function values related, but for any proper, closed, and convex function $\phi$, the Moreau envelope is differentiable, specifically, we have:
\begin{align}
   \label{eq:ppt-2-gradient}
   \nabla \MoreauSub{\gamma}{\phi}{x} = \frac{1}{\gamma}\rbrac{x - \ProxSub{\gamma f}{x}}.
\end{align} 
The above identity indicates that $\phi$ and $\Moreau^{\gamma}_{\phi}$ are intrinsically related.
This relationship plays a key role in our analysis.
We also need the following assumptions on $f$ and $f_i$ to carry out our analysis.
\begin{assumption}[Differentiability]
   \label{asp:diff}
   The function $f_i: \R^d \mapsto \R$ in (\ref{eq:prob-formulation}) is differentiable and bounded from below for all $i \in [n]$.
\end{assumption}
\begin{assumption}[Interpolation regime]
   \label{asp:int-pl-rgm}
   There exists $x_{\star} \in \R^d$ such that $\nabla f_i(x_{\star}) = 0$ for all $i \in [n]$.
\end{assumption}
The same as \citet{li2024power}, we assume that we are in the interpolation regime.
This situation arises in modern deep learning scenarios where the number of parameters, $d$, significantly exceeds the number of data points.
For justifications, we refer the readers to \citet{arora2019fine,montanari2022interpolation}. 
The motivation for this assumption stems from the parallel projection methods (\ref{eq:parallel-proj}) used to solve convex feasibility problems, where the intersection of all convex sets $\cX_i$ is assumed to be non-empty, which is precisely the interpolation assumption of each $f_i$ being the indicator function of $\cX_i$.
\begin{align}
   \label{eq:parallel-proj}
   x_{k+1} = \frac{1}{n}\sum_{i=1}^{n} \Pi_{\cX_i}\rbrac{x_k}.
\end{align}
It is known that for (\ref{eq:parallel-proj}), the use of extrapolation would enhance its performance both in theory and practice \citep{necoara2019randomized}. 
Since $\ProxSub{\gamma f_i}{x_k}$ can be viewed as projection to some level set of $f_i$, it is analogous to $\Pi_{\cX_i}\rbrac{x_k}$.
Therefore, it is reasonable to assume that extrapolation would be effective under the same assumption.
\begin{assumption}[Individual convexity]
   \label{asp:cvx}
   The function $f_i: \R^d\mapsto\R$ is convex for all $i \in [n]$.
   This means that for each $f_i$, 
   \begin{equation}
     \label{eq:asp:convexity}
     0 \leq f_i(x) - f_i(y) - \inner{\nabla f_i(y)}{x - y}, \quad \forall x, y\in \R^d. 
   \end{equation}
\end{assumption}

\begin{assumption}[Smoothness]
   \label{asp:smoothness}
   The function $f_i: \R^d\mapsto\R$ is $L_i$-smooth, $L_i > 0$ for all $i \in [n]$. 
   This means that for each $f_i$,
   \begin{equation}
     \label{eq:asp:smoothness}
     f_i(x) - f_i(y) - \inner{\nabla f_i(y)}{x - y} \leq \frac{L_i}{2}\norm{x - y}^2, \quad \forall x, y\in \R^d. 
   \end{equation}
   We will use $L_{\max}$ to denote $\max_{i \in [n]} L_i$.
\end{assumption}

\begin{assumption}[Global strong convexity]
   \label{asp:stn-cvx}
   The function $f$ is $\mu$-strongly convex, $\mu > 0$. That is 
   \begin{equation*}
       f(x) - f(y) - \inner{\nabla f(y)}{x - y} \geq \frac{\mu}{2}\norm{x - y}^2, \quad \forall x, y\in \R^d. 
   \end{equation*}
\end{assumption}
These are all standard assumptions commonly used in convex optimization. 
We first present our algorithm as \Cref{alg:inexact-fedexprox-full}.
In the following sections, we provide the analysis of this algorithm under different definitions of inexactness, respectively in \Cref{sec:inexact-1} and \Cref{sec:inexact-2}.
Details on how these inexactness levels can be achieved by each client are provided in \Cref{sec:oracles}. 
Finally, numerical experiments validating our results are presented in \Cref{sec:app:exp}.

\begin{algorithm}[t]
	\caption{Inexact {\FEDEXPROX}}
	\label{alg:inexact-fedexprox-full}
	\begin{algorithmic}[1]
	\STATE {\bf Parameters:} extrapolation parameter $\alpha_k = \alpha > 0$, step size for the proximal operator $\gamma > 0$, starting point $x_0 \in \R^d$, number of clients $n$, total number of iterations $K$, proximal solution accuracy $\varepsilon \geq 0$.
	\FOR{$k=0,1,2\dotsc K-1$}   
    	\STATE The server broadcasts the current iterate $x_k$ to each client
        \STATE Each client computes an $\varepsilon$ approximation of the solution $\tx_{i, k+1} \simeq \ProxSub{\gamma f_i}{x_k}$, and sends it back to the server
        \STATE The server computes 
        \begin{equation}
            \label{eq:update-rule}
            x_{k+1} = x_k + \alpha_k\rbrac{\frac{1}{n}\sum_{i=1}^n \tx_{i, k+1} - x_k}.
        \end{equation} 
	\ENDFOR
	\end{algorithmic}
\end{algorithm}

\section{Absolute approximation in distance}
\label{sec:inexact-1}
As previously suggested, we assume that each proximal operator is solved inexactly, and we need to quantify this inexactness in some way.
Notice that client $i$ is required to solve the following minimization problem. 
\begin{equation}
   \label{eq:def-local-opt}
   \min_{z \in \R^d} \aloc{\gamma}{k, i}{z} \eqdef f_i\rbrac{z} + \frac{1}{2\gamma}\norm{z - x_k}^2, 
\end{equation}
where $x_k$ is the current iterate and $\gamma > 0$ is a constant. 
Since we have assumed each function $f_i$ is convex, $\aloc{\gamma}{k, i}{z}$ is $\frac{1}{\gamma}$-strongly convex with $\ProxSub{\gamma f_i}{x_k}$ being its unique minimizer.
One of the most straightforward ways to measure inexactness in this case is through the squared distance to the minimizer, leading to the following definition.
\begin{definition}[Absolute approximation]
   \label{def:inexact-1}
   Given a proper, closed and convex function $\phi: \R^d \mapsto \R$, and a step size $\gamma > 0$, we say that a point $y \in \R^d$ is an $\varepsilon_1$-approximation of $\ProxSub{\gamma \phi}{x}$, if for some $\varepsilon_1 \geq 0$, 
   \begin{equation}
      \label{eq:approx-dist}
      \norm{y - \ProxSub{\gamma f}{x}}^2 \leq \varepsilon_1.
   \end{equation}
\end{definition}
In order to analyze \Cref{alg:inexact-fedexprox-full}, we first transform the update rule given in (\ref{eq:update-rule}) in the following way,
\begin{eqnarray}
   \label{eq:trans-update}
   x_{k+1} &=& x_k + \alpha_k\rbrac{\frac{1}{n}\sum_{i=1}^n\rbrac{\tx_{i, k+1} - \ProxSub{\gamma f_i}{x_k}} + \frac{1}{n}\sum_{i=1}^{n}\ProxSub{\gamma f_i}{x_k} - x_k} \notag \\
   &\overset{(\ref{eq:ppt-2-gradient})}{=}& x_k - \alpha_k\cdot g(x_k),
\end{eqnarray}
where 
\begin{equation}
   \label{eq:biased-estimator}
   g(x_k) \eqdef \underbrace{\frac{1}{n}\sum_{i=1}^{n} \gamma\nabla\MoreauSub{\gamma}{f_i}{x_k}}_{\text{Gradient}} - \underbrace{\frac{1}{n}\sum_{i=1}^{n}\rbrac{\tx_{i, k+1} - \ProxSub{\gamma f_i}{x_k}}}_{\text{Bias}}.
\end{equation}
The above reformulation suggests that \Cref{alg:inexact-fedexprox-full} is in fact, {\SGD} with respect to global objective $\gamma M^{\gamma}\rbrac{x} \eqdef \frac{1}{n}\sum_{i=1}^{n}\gamma \MoreauSub{\gamma}{f_i}{x}$ with a biased gradient estimator.
Compared to {\SGD} with an unbiased gradient estimator, its biased counterpart is less well understood.
However, we are still able to obtain the following convergence guarantee using theories for biased {\SGD} from \cite{demidovich2024guide}.
\begin{theorem}
   \label{thm:1:conv-full-batch-stncvx}
   Assume \Cref{asp:diff} (Differentiability), \Cref{asp:int-pl-rgm} (Interpolation Regime), \Cref{asp:cvx} (Individual convexity), \Cref{asp:smoothness} (Smoothness) and \Cref{asp:stn-cvx} (Global strong convexity) hold.
   If each client only computes a $\varepsilon_1$-absolute approximation $\tx_{i, k+1}$ in squared distance of $\ProxSub{\gamma f_i}{x_k}$ at every iteration, such that $\norm{\tx_{i, k+1} - \ProxSub{\gamma f_i}{x_k}}^2 \leq \varepsilon_1$.
   Then we have the following convergence guarantee or \Cref{alg:inexact-fedexprox-full}: For a constant extrapolation parameter satisfying $0 < \alpha \leq \frac{1}{4}\cdot\frac{1}{\gamma L_{\gamma}}$, where $\gamma$ is the step size of the proximity operator, $\alpha_k = \alpha$ is a constant extrapolation parameter, $L_{\gamma}$ is the smoothness constant of $M^{\gamma}$.
   The last iterate $x_K$ satisfy
   \begin{align*}
      \cE_K &\leq \rbrac{1 - \frac{\alpha\gamma\mu}{8\rbrac{1 + \gamma L_{\max}}}}^K \cE_0 + \frac{4\varepsilon_1\rbrac{1 + \gamma L_{\max}}}{\mu}\cdot\rbrac{2\alpha L_{\gamma} + \frac{1}{\gamma}}, 
   \end{align*}
   where $\cE_k = \gamma M^{\gamma}\rbrac{x_k} - \gamma M^{\gamma}_{\inf}$. 
   Specifically, when we choose $\alpha = \frac{1}{4}\cdot \frac{1}{\gamma L_{\gamma}}$, we have 
   \begin{align*}
      \Delta_K \leq \rbrac{1 - \frac{\mu}{32L_{\gamma}\rbrac{1 + \gamma L_{\max}}}}^K\frac{L_{\gamma}\rbrac{1 + \gamma L_{\max}}}{\mu} \cdot \Delta_0 + 12\varepsilon_1\cdot\rbrac{\frac{\nicefrac{1}{\gamma} + L_{\max}}{\mu}}^2,
   \end{align*}
   where $\Delta_K = \norm{x_K - x_\star}^2$, $x_\star$ is a minimizer of $f$.
\end{theorem}
For the sake of brevity in the following discussion, we will use the notation $\cE_k = \gamma M^{\gamma}\rbrac{x_k} - \gamma M^{\gamma}_{\inf}$, where $M^{\gamma}_{\inf}$ denotes the infimum of $M^{\gamma}$, $\Delta_k = \norm{x_k - x_\star}^2$, where $x_\star$ is a minimizer of $M^{\gamma}$.
Notice that since we are in the interpolation regime, according to \Cref{fact:5:global:equivalence}, the minimizer of $M^{\gamma}$ is also a minimizer of $f$.
Note that instead of converging to the exact minimizer $x_\star$, the algorithm converges to a neighborhood whose size depends on both $\varepsilon_1$ and $\gamma$; the smaller $\gamma$ is, the larger the neighborhood becomes.
This can be understood intuitively: A smaller $\gamma$ means less progress is made per iteration, leading to a larger accumulated error as the total number of iterations increases.
The parameter $\varepsilon_1$ can be arbitrarily large, and the convergence guarantee still holds, indicating that the theory presented is quite general.
However, as $\varepsilon_1$ increases, the size of the neighborhood grows proportionally, which limits the practical significance of the result.
When $\varepsilon_1 = 0$, the neighborhood diminishes, and we obtain an iteration complexity of $\tilde{\cO}\rbrac{\frac{L_{\gamma}\rbrac{1 + \gamma L_{\max}}}{\mu}}$\footnote{We leave out the log factor in $\tilde{\cO}\rbrac{\cdot}$ notation.}, which recovers the result of \citet{li2024power} up to a constant factor.
The optimal constant extrapolation parameter is now given by $\alpha_\star = \frac{1}{4}\cdot\frac{1}{\gamma L_{\gamma}}$ which is $4$ times smaller than that of \citet{li2024power}.

\section{Relative approximation in distance}
\label{sec:inexact-2}
\Cref{thm:1:conv-full-batch-stncvx} offers a general theoretical framework for understanding the behavior of \Cref{alg:inexact-fedexprox-full}.
However, a key challenge with \Cref{alg:inexact-fedexprox-full} which utilizes inexact proximal solutions that satisfy \Cref{def:inexact-1}, is that, unless the proximal operators are solved exactly, convergence will always be limited to a neighborhood of the solution.
The underlying reason is that, as the algorithm progresses, the gradient term in the gradient estimator $g(x_k)$ diminishes, whereas the bias term remains unchanged.
Building on this observation, we propose employing a different type of approximation, specifically an approximation in relative distance, as defined below.
\begin{definition}[Relative approximation]
   \label{def:inexact-2}
   Given a convex function $\phi: \R^d \mapsto \R$ and a stepsize $\gamma > 0$, we say that a point $y \in \R^d$ is a $\varepsilon_2$-relative approximation of $\ProxSub{\gamma \phi}{x}$, if for some $\varepsilon_2 \in [0, 1)$, 
   \begin{equation}
      \label{eq:index2-approx}
      \norm{y - \ProxSub{\gamma \phi}{x}}^2 \leq \varepsilon_2 \cdot \norm{x - \ProxSub{\gamma \phi}{x}}^2.
   \end{equation}
\end{definition}
We impose the requirement that the coefficient $\varepsilon_2$ be less than $1$ to ensure that the next iterate is no worse than the current one.
As we can observe, if the approximation of the solution for each proximal operator satisfies \Cref{def:inexact-2}, both the gradient term and the bias term diminish as the algorithm progresses, ensuring convergence to the exact solution.
Using the theory of biased {\SGD}, we can obtain the following theorem.
\begin{theorem}
   \label{thm:020-full}
   Assume all the assumptions mentioned in \Cref{thm:1:conv-full-batch-stncvx} also hold here.
   If each client only computes a $\varepsilon_2$-relative approximation $\tx_{i, k+1}$ in distance with $\varepsilon_2 < \nicefrac{\mu^2}{4L^2_{\max}}$, such that $\norm{\tx_{i, k+1} - \ProxSub{\gamma f_i}{x_k}}^2 \leq \varepsilon_2 \cdot \norm{x_k - \ProxSub{\gamma f_i}{x_k}}^2$.
   If we are running \Cref{alg:inexact-fedexprox-full} with $\alpha_k = \alpha$ satisfying
   \begin{align*}
      0 < \alpha \leq \frac{1}{\gamma L_{\gamma}}\cdot\frac{\mu - 2\sqrt{\varepsilon_2}L_{\max}}{\mu + 4\sqrt{\varepsilon_2}L_{\max} + 4\varepsilon_2 L_{\max}}.
   \end{align*}
   Then the iterates generated by \Cref{alg:inexact-fedexprox-full} satisfies 
   \begin{equation*}
      \cE_K \leq \rbrac{1 - \alpha\cdot\frac{\gamma\rbrac{\mu - 2\sqrt{\varepsilon_2} L_{\max}}}{4\rbrac{1 + \gamma L_{\max}}}}^K \cE_0.
   \end{equation*}
   Specifically, if we choose the largest $\alpha$ possible, we have
   \begin{align*}
      \Delta_K \leq \rbrac{1 - \frac{\mu}{4L_{\gamma}\rbrac{1 + \gamma L_{\max}}} \cdot S\rbrac{\varepsilon_2}}^K\cdot\frac{L{\gamma}\rbrac{1 + \gamma L_{\max}}}{\mu}\Delta_0,
   \end{align*} 
   where $S(\varepsilon_2) \eqdef \frac{\rbrac{\mu - 2\sqrt{\varepsilon_2}L_{\max}}\rbrac{1 - 2\sqrt{\varepsilon_2}\frac{L_{\max}}{\mu}}}{\mu + 4\sqrt{\varepsilon_2}L_{\max} + 4\varepsilon_2 L_{\max}}
   $ satisfies $0 < S(\varepsilon_2) \leq 1$ is the factor of slowing down due to inexact proximal operator evaluation.
\end{theorem}

Observe that when $\varepsilon_2 = 0$, meaning the proximal operators are solved exactly, the optimal extrapolation is $\alpha = \frac{1}{\gamma L_{\gamma}}$ and the iteration complexity is $\tilde{\cO}\rbrac{\frac{L_{\gamma}\rbrac{1 + \gamma L_{\max}}}{\mu}}$.
This recovers the exact result from \citet{li2024power}. 
In the case of an inexact solution, as $\varepsilon_2$ increases, both $\alpha$ and $S(\varepsilon_2)$ decrease, leading to a slower rate of convergence.
Note that arbitrary rough approximations are not permissible in this case, as $\varepsilon_2$ must satisfy $\varepsilon_2 < \frac{\mu^2}{4L_{\max}^2}$.

It is worthwhile noting that \Cref{def:inexact-2} is connected to the concept of compression.
Indeed, in our case we have $x_k - \ProxSub{\gamma f_i}{x_k} = \gamma\nabla\MoreauSub{\gamma}{f_i}{x_k}$, while $\tx_{i, k+1} - \ProxSub{\gamma f_i}{x_k}$ can be interpreted as the gradient after compression, that is, $\cC(\gamma\nabla\MoreauSub{\gamma}{f_i}{x_k})$.
This indicates that \Cref{alg:inexact-fedexprox-full} with approximation satisfying \Cref{def:inexact-2} can be viewed as compressed gradient descent with biased compressor.
We obtain the following convergence guarantee based on theory provided by \citet{beznosikov2023biased}.
\begin{theorem}
   \label{thm:010-stncvx-rela}
   Assume all assumptions of \Cref{thm:1:conv-full-batch-stncvx} hold.
   Let the approximation $\tx_{i, k+1}$ all satisfies \Cref{def:inexact-2} with $\varepsilon_2 < \nicefrac{\mu}{4L_{\max}}$, that is $ \norm{\tx_{i, k+1} - \ProxSub{\gamma f_i}{x_k}}^2 \leq \varepsilon_2 \cdot \norm{x_k - \ProxSub{\gamma f_i}{x_k}}^2
   $.
   If we are running \Cref{alg:inexact-fedexprox-full} with $\alpha_k = \alpha \in (0, \frac{1}{\gamma L_{\gamma}}]$, we have the iterates produced by it satisfying
   \begin{equation*}
      \cE_K \leq \rbrac{1 - \rbrac{1 - \frac{4\varepsilon_2 L_{\max}}{\mu}}\cdot \frac{\gamma\mu}{4\rbrac{1 + \gamma L_{\max}}}\cdot\alpha}^K\cE_0.
   \end{equation*}
   specifically, if we take the largest extrapolation ($\alpha = \frac{1}{\gamma L_{\gamma}} > 1$) possible, we have 
   \begin{align*}
      \Delta_K \leq \rbrac{1 - \rbrac{1 - \frac{4\varepsilon_2 L_{\max}}{\mu}}\cdot \frac{\mu}{4L_{\gamma}\rbrac{1 + \gamma L_{\max}}}}^K \cdot \frac{L_{\gamma}\rbrac{1 + \gamma L_{\max}}}{\mu}\Delta_0.
   \end{align*}
\end{theorem}
The convergence guarantee obtained in this way is sharper, indeed, \Cref{thm:010-stncvx-rela} suggests that as long as $\varepsilon_2 < \nicefrac{\mu}{4L}$, we are able to pick $\alpha = \nicefrac{1}{\gamma L_{\gamma}}$\footnote{It is shown in \citet{li2024power} that $\nicefrac{1}{\gamma L_{\gamma}} > 1$, which justifies why $\alpha$ is called the extrapolation parameter.}which is the optimal extrapolation for exact proximal computation given in \citet{li2024power}.
Notably, this implies that extrapolation is an effective technique for accelerating the algorithm in this setting, regardless of inexact proximal operator evaluations.
Same as \Cref{thm:020-full}, the convergence is slowed down by the approximation, and in the case of $\varepsilon_2 = 0$, we recover the result in \citet{li2024power}

\section{Achieving the level of inexactness}
\label{sec:oracles}
To fully comprehend the overall complexity of \Cref{alg:inexact-fedexprox-full}, it is essential to examine whether the inexactness in evaluating the proximal operators can be effectively achieved.
Since each $\ProxSub{\gamma f_i}{x_k}$ is computed locally by the corresponding client, the client has access to all the necessary data points for the computation.
Thus, the most straightforward approach is to have each client perform {\GD}.
\begin{theorem}[Local computation via {\GD}]
   \label{thm:complexity}
   Assume \Cref{asp:diff} (Differentiability), \Cref{asp:cvx} (Individual convexity) and \Cref{asp:smoothness} (Smoothness) hold.
   The iteration complexity for the $i$-th client to provide an approximation using {\GD} in the $k$-th iteration with local step size $\eta_i = \frac{\gamma}{1 + \gamma L_i}$, satisfying \Cref{def:inexact-1} is 
   $
      \cO\rbrac{\rbrac{1+ \gamma L_i}\log\rbrac{\nicefrac{\norm{x_k - \ProxSub{\gamma f_i}{x_k}}^2}{\varepsilon_1}}},
   $
   and for \Cref{def:inexact-2}, it is 
   $
      \cO\rbrac{\rbrac{1+ \gamma L_i}\log\rbrac{\nicefrac{1}{\varepsilon_2}}}.
   $
   % If we are using {\AGD} with local stepsize $\eta_i = \frac{\gamma}{1 + \gamma L_i}$ and momentum parameter $\alpha_i = \frac{\sqrt{1 + \gamma L_i} - 1}{\sqrt{1 + \gamma L_i} + 1}$, the iteration complexities of \Cref{alg:inexact-fedexprox-full} with approximations satisfying \Cref{def:inexact-1}, \Cref{def:inexact-2} are 
   % {
   %    \small
   %    \begin{align*}
   %       \cO\rbrac{\sqrt{1 + \gamma L_i}\log\rbrac{\frac{\rbrac{1 + \gamma L_i}\cdot\norm{x_k - \ProxSub{\gamma f_i}{x_k}}^2}{\varepsilon_1}}}; \quad \cO\rbrac{\sqrt{1 + \gamma L_i}\log\rbrac{\frac{1 + \gamma L_i}{\varepsilon_2}}},
   %    \end{align*}
   % }
   % respectively.
\end{theorem}
Note that there are no constraints on $\varepsilon_1$, and since $\norm{x_k - \ProxSub{\gamma f_i}{x_k}}^2 \leq \norm{\gamma \nabla f(x_k)}^2$ by (\ref{eq:upperbound-gd-comp}), it is straightforward to adjust {\GD} to optimize the approximation. 
However, for $\varepsilon_2$, we require $\varepsilon_2 < \frac{\mu}{4L_{\max}}$.
In practice, $\varepsilon_2$ can be set to a sufficiently small value to satisfy this condition, though this will increase the number of local iterations performed by each client.
The complexity bounds also indicate that as the local step size $\gamma$ increases, it becomes more challenging to compute the approximation.
Alternatively, other algorithms can be employed to find such an approximation.
For instance, by leveraging the structure in (\ref{eq:substructure}), {\SGD} can be used as a local solver for the proximal operator when computational resources are limited.
We can use the accelerated gradient descent ({\AGD}) of \citet{nesterov2013introductory} to obtain a better iteration complexity for each client.
\begin{theorem}[Local computation via {\AGD}]
   \label{thm:complexity-accelerated}
   Assume all assumptions mentioned in \Cref{thm:complexity} hold.
   The iteration complexities for the $i$-th client to provide an approximation in the $k$-the iteration using {\AGD} with local step size $\eta_i = \frac{\gamma}{1 + \gamma L_i}$ and momentum parameter $\alpha_i = \frac{\sqrt{1 + \gamma L_i} - 1}{\sqrt{1 + \gamma L_i} + 1}$, satisfying \Cref{def:inexact-1}, \Cref{def:inexact-2} are 
   {
      \small
      \begin{align*}
         \cO\rbrac{\sqrt{1 + \gamma L_i}\log\rbrac{\frac{\rbrac{1 + \gamma L_i}\cdot\norm{x_k - \ProxSub{\gamma f_i}{x_k}}^2}{\varepsilon_1}}}; \quad \cO\rbrac{\sqrt{1 + \gamma L_i}\log\rbrac{\frac{1 + \gamma L_i}{\varepsilon_2}}},
      \end{align*}
   }
   respectively.
\end{theorem}
Finally, we provide numerical evidence to support our theoretical findings.
We refer the readers to \Cref{sec:app:exp} for the details of the settings and the corresponding experiments.

\section{Conclusions}
% In this paper, we examine the convergence properties of inexact {\FEDEXPROX}.
% We establish convergence guarantees under a general metric for inexactness, indicating convergence to a neighborhood of the solution.
% Building on this, we propose a method to eliminate the neighborhood and offer a more refined analysis using the theory of biased compression for {\CGD}. 
% Our findings are supported by numerical experiments.

\subsection{Limitations}
\label{sec:limitations}
Despite achieving satisfactory results in the full-batch setting, the client sampling setting did not yield similar outcomes. This may be attributed to the nature of biased compression, which likely requires adjustments to the algorithm itself for resolution.
Nonetheless, we provide the analysis in \Cref{sec:minibatch} for reference.
Unlike \citet{li2024power}, the presence of bias makes it unclear how to incorporate adaptive step-size rules such as gradient diversity in our case.
The only permissible inexactness for gradient diversity arises from client sub-sampling in the interpolation regime.

\subsection{Future work}
\label{sec:future}
There are still open problems to be addressed. 
For example, can \Cref{alg:inexact-fedexprox-full} be modified to incorporate the benefits of error feedback? 
Is it possible to eliminate the interpolation regime assumption while still demonstrating that extrapolation is theoretically beneficial for {\FEDEXPROX}?  
Another direction that may be of independent interest is to develop adaptive rules of determining the step size for {\SGD} with biased update.

\newpage

\bibliography{iclr2025_conference}

\begin{thebibliography}{68}
\providecommand{\natexlab}[1]{#1}
\providecommand{\url}[1]{\texttt{#1}}
\expandafter\ifx\csname urlstyle\endcsname\relax
  \providecommand{\doi}[1]{doi: #1}\else
  \providecommand{\doi}{doi: \begingroup \urlstyle{rm}\Url}\fi

\bibitem[Alistarh et~al.(2017)Alistarh, Grubic, Li, Tomioka, and
  Vojnovic]{alistarh2017qsgd}
Dan Alistarh, Demjan Grubic, Jerry Li, Ryota Tomioka, and Milan Vojnovic.
\newblock {QSGD}: Communication-efficient {SGD} via gradient quantization and
  encoding.
\newblock \emph{Advances in Neural Information Processing Systems}, 30, 2017.

\bibitem[Alistarh et~al.(2018)Alistarh, Hoefler, Johansson, Konstantinov,
  Khirirat, and Renggli]{alistarh2018convergence}
Dan Alistarh, Torsten Hoefler, Mikael Johansson, Nikola Konstantinov, Sarit
  Khirirat, and C{\'e}dric Renggli.
\newblock The convergence of sparsified gradient methods.
\newblock \emph{Advances in Neural Information Processing Systems}, 31, 2018.

\bibitem[Anitescu(2000)]{anitescu2000degenerate}
Mihai Anitescu.
\newblock Degenerate nonlinear programming with a quadratic growth condition.
\newblock \emph{SIAM Journal on Optimization}, 10\penalty0 (4):\penalty0
  1116--1135, 2000.

\bibitem[Arora et~al.(2019)Arora, Du, Hu, Li, and Wang]{arora2019fine}
Sanjeev Arora, Simon Du, Wei Hu, Zhiyuan Li, and Ruosong Wang.
\newblock Fine-grained analysis of optimization and generalization for
  overparameterized two-layer neural networks.
\newblock In \emph{International Conference on Machine Learning}, pp.\
  322--332. PMLR, 2019.

\bibitem[Asi \& Duchi(2019)Asi and Duchi]{asi2019stochastic}
Hilal Asi and John~C Duchi.
\newblock Stochastic (approximate) proximal point methods: Convergence,
  optimality, and adaptivity.
\newblock \emph{SIAM Journal on Optimization}, 29\penalty0 (3):\penalty0
  2257--2290, 2019.

\bibitem[Bauschke et~al.(2006)Bauschke, Combettes, and
  Kruk]{bauschke2006extrapolation}
Heinz~H Bauschke, Patrick~L Combettes, and Serge~G Kruk.
\newblock Extrapolation algorithm for affine-convex feasibility problems.
\newblock \emph{Numerical Algorithms}, 41:\penalty0 239--274, 2006.

\bibitem[Beck(2017)]{beck2017first}
Amir Beck.
\newblock \emph{First-order methods in optimization}.
\newblock SIAM, 2017.

\bibitem[Bertsekas(2011)]{bertsekas2011incremental}
Dimitri~P Bertsekas.
\newblock Incremental proximal methods for large scale convex optimization.
\newblock \emph{Mathematical Programming}, 129\penalty0 (2):\penalty0 163--195,
  2011.

\bibitem[Beznosikov et~al.(2023)Beznosikov, Horv{\'a}th, Richt{\'a}rik, and
  Safaryan]{beznosikov2023biased}
Aleksandr Beznosikov, Samuel Horv{\'a}th, Peter Richt{\'a}rik, and Mher
  Safaryan.
\newblock On biased compression for distributed learning.
\newblock \emph{Journal of Machine Learning Research}, 24\penalty0
  (276):\penalty0 1--50, 2023.

\bibitem[Bianchi(2016)]{bianchi2016ergodic}
Pascal Bianchi.
\newblock Ergodic convergence of a stochastic proximal point algorithm.
\newblock \emph{SIAM Journal on Optimization}, 26\penalty0 (4):\penalty0
  2235--2260, 2016.

\bibitem[Bubeck et~al.(2015)]{bubeck2015convex}
S{\'e}bastien Bubeck et~al.
\newblock Convex optimization: Algorithms and complexity.
\newblock \emph{Foundations and Trends{\textregistered} in Machine Learning},
  8\penalty0 (3-4):\penalty0 231--357, 2015.

\bibitem[Censor et~al.(2001)Censor, Elfving, and Herman]{censor2001averaging}
Y~Censor, T~Elfving, and GT~Herman.
\newblock Averaging strings of sequential iterations for convex feasibility
  problems.
\newblock In \emph{Studies in Computational Mathematics}, volume~8, pp.\
  101--113. Elsevier, 2001.

\bibitem[Combettes(1997)]{combettes1997convex}
Patrick~L Combettes.
\newblock Convex set theoretic image recovery by extrapolated iterations of
  parallel subgradient projections.
\newblock \emph{IEEE Transactions on Image Processing}, 6\penalty0
  (4):\penalty0 493--506, 1997.

\bibitem[Condat et~al.(2023)Condat, Kitahara, Contreras, and
  Hirabayashi]{condat2023proximal}
Laurent Condat, Daichi Kitahara, Andr{\'e}s Contreras, and Akira Hirabayashi.
\newblock Proximal splitting algorithms for convex optimization: A tour of
  recent advances, with new twists.
\newblock \emph{SIAM Review}, 65\penalty0 (2):\penalty0 375--435, 2023.

\bibitem[Davis \& Drusvyatskiy(2019)Davis and
  Drusvyatskiy]{davis2019stochastic}
Damek Davis and Dmitriy Drusvyatskiy.
\newblock Stochastic model-based minimization of weakly convex functions.
\newblock \emph{SIAM Journal on Optimization}, 29\penalty0 (1):\penalty0
  207--239, 2019.

\bibitem[Demidovich et~al.(2024)Demidovich, Malinovsky, Sokolov, and
  Richt{\'a}rik]{demidovich2024guide}
Yury Demidovich, Grigory Malinovsky, Igor Sokolov, and Peter Richt{\'a}rik.
\newblock A guide through the zoo of biased sgd.
\newblock \emph{Advances in Neural Information Processing Systems}, 36, 2024.

\bibitem[Ghadimi \& Lan(2013)Ghadimi and Lan]{ghadimi2013stochastic}
Saeed Ghadimi and Guanghui Lan.
\newblock Stochastic first-and zeroth-order methods for nonconvex stochastic
  programming.
\newblock \emph{SIAM Journal on Optimization}, 23\penalty0 (4):\penalty0
  2341--2368, 2013.

\bibitem[Gong \& Ye(2014)Gong and Ye]{gong2014linear}
Pinghua Gong and Jieping Ye.
\newblock Linear convergence of variance-reduced stochastic gradient without
  strong convexity.
\newblock \emph{arXiv preprint arXiv:1406.1102}, 2014.

\bibitem[Gorbunov et~al.(2020)Gorbunov, Hanzely, and
  Richt{\'a}rik]{gorbunov2020unified}
Eduard Gorbunov, Filip Hanzely, and Peter Richt{\'a}rik.
\newblock A unified theory of {SGD}: Variance reduction, sampling, quantization
  and coordinate descent.
\newblock In \emph{International Conference on Artificial Intelligence and
  Statistics}, pp.\  680--690. PMLR, 2020.

\bibitem[Gorbunov et~al.(2021)Gorbunov, Burlachenko, Li, and
  Richt{\'a}rik]{gorbunov2021marina}
Eduard Gorbunov, Konstantin~P Burlachenko, Zhize Li, and Peter Richt{\'a}rik.
\newblock {MARINA}: Faster non-convex distributed learning with compression.
\newblock In \emph{International Conference on Machine Learning}, pp.\
  3788--3798. PMLR, 2021.

\bibitem[Gower et~al.(2020)Gower, Schmidt, Bach, and
  Richt{\'a}rik]{gower2020variance}
Robert~M Gower, Mark Schmidt, Francis Bach, and Peter Richt{\'a}rik.
\newblock Variance-reduced methods for machine learning.
\newblock \emph{Proceedings of the IEEE}, 108\penalty0 (11):\penalty0
  1968--1983, 2020.

\bibitem[Gower et~al.(2019)Gower, Loizou, Qian, Sailanbayev, Shulgin, and
  Richt{\'a}rik]{gower2019sgd}
Robert~Mansel Gower, Nicolas Loizou, Xun Qian, Alibek Sailanbayev, Egor
  Shulgin, and Peter Richt{\'a}rik.
\newblock {SGD}: General analysis and improved rates.
\newblock In \emph{International Conference on Machine Learning}, pp.\
  5200--5209. PMLR, 2019.

\bibitem[Iutzeler \& Hendrickx(2019)Iutzeler and
  Hendrickx]{iutzeler2019generic}
Franck Iutzeler and Julien~M Hendrickx.
\newblock A generic online acceleration scheme for optimization algorithms via
  relaxation and inertia.
\newblock \emph{Optimization Methods and Software}, 34\penalty0 (2):\penalty0
  383--405, 2019.

\bibitem[Jhunjhunwala et~al.(2023)Jhunjhunwala, Wang, and
  Joshi]{jhunjhunwala2023fedexp}
Divyansh Jhunjhunwala, Shiqiang Wang, and Gauri Joshi.
\newblock Fed{E}x{P}: Speeding up federated averaging via extrapolation.
\newblock In \emph{International Conference on Learning Representations}, 2023.

\bibitem[Johnson \& Zhang(2013)Johnson and Zhang]{johnson2013accelerating}
Rie Johnson and Tong Zhang.
\newblock Accelerating stochastic gradient descent using predictive variance
  reduction.
\newblock \emph{Advances in neural information processing systems}, 26, 2013.

\bibitem[Jourani et~al.(2014)Jourani, Thibault, and
  Zagrodny]{jourani2014differential}
Abderrahim Jourani, Lionel Thibault, and Dariusz Zagrodny.
\newblock Differential properties of the moreau envelope.
\newblock \emph{Journal of Functional Analysis}, 266\penalty0 (3):\penalty0
  1185--1237, 2014.

\bibitem[Kaczmarz(1937)]{kaczmarz1993approximate}
Stefan Kaczmarz.
\newblock Approximate solution of systems of linear equations.
\newblock \emph{International Journal of Control}, 57\penalty0 (6):\penalty0
  1269--1271, 1937.

\bibitem[Karagulyan et~al.(2024)Karagulyan, Shulgin, Sadiev, and
  Richt{\'a}rik]{karagulyan2024spam}
Avetik Karagulyan, Egor Shulgin, Abdurakhmon Sadiev, and Peter Richt{\'a}rik.
\newblock Spam: Stochastic proximal point method with momentum variance
  reduction for non-convex cross-device federated learning.
\newblock \emph{arXiv preprint arXiv:2405.20127}, 2024.

\bibitem[Karimireddy et~al.(2019)Karimireddy, Rebjock, Stich, and
  Jaggi]{karimireddy2019error}
Sai~Praneeth Karimireddy, Quentin Rebjock, Sebastian Stich, and Martin Jaggi.
\newblock Error feedback fixes signsgd and other gradient compression schemes.
\newblock In \emph{International Conference on Machine Learning}, pp.\
  3252--3261. PMLR, 2019.

\bibitem[Karimireddy et~al.(2020)Karimireddy, Kale, Mohri, Reddi, Stich, and
  Suresh]{karimireddy2020scaffold}
Sai~Praneeth Karimireddy, Satyen Kale, Mehryar Mohri, Sashank Reddi, Sebastian
  Stich, and Ananda~Theertha Suresh.
\newblock Scaffold: stochastic controlled averaging for federated learning.
\newblock In \emph{International Conference on Machine Learning}, pp.\
  5132--5143. PMLR, 2020.

\bibitem[Khaled \& Jin(2022)Khaled and Jin]{khaled2022faster}
Ahmed Khaled and Chi Jin.
\newblock Faster federated optimization under second-order similarity.
\newblock In \emph{The Eleventh International Conference on Learning
  Representations}, 2022.

\bibitem[Khirirat et~al.(2018)Khirirat, Feyzmahdavian, and
  Johansson]{khirirat2018distributed}
Sarit Khirirat, Hamid~Reza Feyzmahdavian, and Mikael Johansson.
\newblock Distributed learning with compressed gradients.
\newblock \emph{arXiv preprint arXiv:1806.06573}, 2018.

\bibitem[Kone\v{c}n{\'y} et~al.(2016)Kone\v{c}n{\'y}, McMahan, Yu,
  Richt{\'a}rik, Suresh, and Bacon]{Konecny2016federated}
Jakub Kone\v{c}n{\'y}, H~Brendan McMahan, Felix~X Yu, Peter Richt{\'a}rik,
  Ananda~Theertha Suresh, and Dave Bacon.
\newblock Federated learning: Strategies for improving communication
  efficiency.
\newblock \emph{arXiv preprint arXiv:1610.05492}, 8, 2016.

\bibitem[Li et~al.(2023)Li, Karagulyan, and Richt{\'a}rik]{li2023marina}
Hanmin Li, Avetik Karagulyan, and Peter Richt{\'a}rik.
\newblock Variance reduced distributed non-convex optimization using matrix
  stepsizes.
\newblock \emph{arXiv preprint arXiv:2310.04614}, 2023.

\bibitem[Li et~al.(2024{\natexlab{a}})Li, Acharya, and
  Richt{\'a}rik]{li2024power}
Hanmin Li, Kirill Acharya, and Peter Richt{\'a}rik.
\newblock The power of extrapolation in federated learning.
\newblock \emph{arXiv preprint arXiv:2405.13766}, 2024{\natexlab{a}}.

\bibitem[Li et~al.(2024{\natexlab{b}})Li, Karagulyan, and
  Richt{\'a}rik]{li2023det}
Hanmin Li, Avetik Karagulyan, and Peter Richt{\'a}rik.
\newblock {Det-CGD}: Compressed gradient descent with matrix stepsizes for
  non-convex optimization.
\newblock In \emph{International Conference on Learning Representations},
  2024{\natexlab{b}}.

\bibitem[Li et~al.(2020)Li, Sahu, Zaheer, Sanjabi, Talwalkar, and
  Smith]{li2020federated}
Tian Li, Anit~Kumar Sahu, Manzil Zaheer, Maziar Sanjabi, Ameet Talwalkar, and
  Virginia Smith.
\newblock Federated optimization in heterogeneous networks.
\newblock \emph{Proceedings of Machine Learning and Systems}, 2:\penalty0
  429--450, 2020.

\bibitem[Liu \& Wright(2015)Liu and Wright]{liu2015asynchronous}
Ji~Liu and Stephen~J Wright.
\newblock Asynchronous stochastic coordinate descent: Parallelism and
  convergence properties.
\newblock \emph{SIAM Journal on Optimization}, 25\penalty0 (1):\penalty0
  351--376, 2015.

\bibitem[Liu et~al.(2020)Liu, Gao, and Yin]{liu2020improved}
Yanli Liu, Yuan Gao, and Wotao Yin.
\newblock An improved analysis of stochastic gradient descent with momentum.
\newblock \emph{Advances in Neural Information Processing Systems},
  33:\penalty0 18261--18271, 2020.

\bibitem[Loizou \& Richt{\'a}rik(2017)Loizou and
  Richt{\'a}rik]{loizou2017linearly}
Nicolas Loizou and Peter Richt{\'a}rik.
\newblock Linearly convergent stochastic heavy ball method for minimizing
  generalization error.
\newblock \emph{arXiv preprint arXiv:1710.10737}, 2017.

\bibitem[Mangasarian \& Solodov(1993)Mangasarian and
  Solodov]{mangasarian1993backpropagation}
Olvi~L Mangasarian and Mikhail~V Solodov.
\newblock Backpropagation convergence via deterministic nonmonotone perturbed
  minimization.
\newblock \emph{Advances in Neural Information Processing Systems}, 6, 1993.

\bibitem[Martinet(1972)]{martinet1972algorithmes}
Bernard Martinet.
\newblock \emph{Algorithmes pour la r{\'e}solution de probl{\`e}mes
  d'optimisation et de minimax}.
\newblock PhD thesis, Universit{\'e} Joseph-Fourier-Grenoble I, 1972.

\bibitem[McMahan et~al.(2017)McMahan, Moore, Ramage, Hampson, and
  y~Arcas]{mcmahan2017communication}
Brendan McMahan, Eider Moore, Daniel Ramage, Seth Hampson, and Blaise~Aguera
  y~Arcas.
\newblock Communication-efficient learning of deep networks from decentralized
  data.
\newblock In \emph{Artificial Intelligence and Statistics}, pp.\  1273--1282.
  PMLR, 2017.

\bibitem[Mishchenko et~al.(2023)Mishchenko, Hanzely, and
  Richt{\'a}rik]{mishchenko2023convergence}
Konstantin Mishchenko, Slavomir Hanzely, and Peter Richt{\'a}rik.
\newblock Convergence of first-order algorithms for meta-learning with {M}oreau
  envelopes.
\newblock \emph{arXiv preprint arXiv:2301.06806}, 2023.

\bibitem[Montanari \& Zhong(2022)Montanari and
  Zhong]{montanari2022interpolation}
Andrea Montanari and Yiqiao Zhong.
\newblock The interpolation phase transition in neural networks: Memorization
  and generalization under lazy training.
\newblock \emph{The Annals of Statistics}, 50\penalty0 (5):\penalty0
  2816--2847, 2022.

\bibitem[Moreau(1965)]{moreau1965proximite}
Jean-Jacques Moreau.
\newblock Proximit{\'e} et dualit{\'e} dans un espace {H}ilbertien.
\newblock \emph{Bulletin de la Soci{\'e}t{\'e} Math{\'e}matique de France},
  93:\penalty0 273--299, 1965.

\bibitem[Moulines \& Bach(2011)Moulines and Bach]{moulines2011non}
Eric Moulines and Francis Bach.
\newblock Non-asymptotic analysis of stochastic approximation algorithms for
  machine learning.
\newblock \emph{Advances in Neural Information Processing Systems}, 24, 2011.

\bibitem[Necoara et~al.(2019)Necoara, Richt{\'a}rik, and
  Patrascu]{necoara2019randomized}
Ion Necoara, Peter Richt{\'a}rik, and Andrei Patrascu.
\newblock Randomized projection methods for convex feasibility: Conditioning
  and convergence rates.
\newblock \emph{SIAM Journal on Optimization}, 29\penalty0 (4):\penalty0
  2814--2852, 2019.

\bibitem[Nesterov(2004)]{nesterov2013introductory}
Yurii Nesterov.
\newblock \emph{Introductory lectures on convex optimization: A basic course}.
\newblock Kluwer Academic Publishers, 2004.

\bibitem[Parikh et~al.(2014)Parikh, Boyd, et~al.]{parikh2014proximal}
Neal Parikh, Stephen Boyd, et~al.
\newblock Proximal algorithms.
\newblock \emph{Foundations and {T}rends{\textregistered} in Optimization},
  1\penalty0 (3):\penalty0 127--239, 2014.

\bibitem[Patrascu \& Necoara(2018)Patrascu and
  Necoara]{patrascu2018nonasymptotic}
Andrei Patrascu and Ion Necoara.
\newblock Nonasymptotic convergence of stochastic proximal point methods for
  constrained convex optimization.
\newblock \emph{Journal of Machine Learning Research}, 18\penalty0
  (198):\penalty0 1--42, 2018.

\bibitem[Pierra(1984)]{pierra1984decomposition}
Guy Pierra.
\newblock Decomposition through formalization in a product space.
\newblock \emph{Mathematical Programming}, 28:\penalty0 96--115, 1984.

\bibitem[Planiden \& Wang(2016)Planiden and Wang]{planiden2016strongly}
Chayne Planiden and Xianfu Wang.
\newblock Strongly convex functions, {M}oreau envelopes, and the generic nature
  of convex functions with strong minimizers.
\newblock \emph{SIAM Journal on Optimization}, 26\penalty0 (2):\penalty0
  1341--1364, 2016.

\bibitem[Planiden \& Wang(2019)Planiden and Wang]{planiden2019proximal}
Chayne Planiden and Xianfu Wang.
\newblock Proximal mappings and {M}oreau envelopes of single-variable convex
  piecewise cubic functions and multivariable gauge functions.
\newblock \emph{Nonsmooth Optimization and Its Applications}, pp.\  89--130,
  2019.

\bibitem[Polyak(1964)]{polyak1964gradient}
Boris~T Polyak.
\newblock Gradient methods for solving equations and inequalities.
\newblock \emph{USSR Computational Mathematics and Mathematical Physics},
  4\penalty0 (6):\penalty0 17--32, 1964.

\bibitem[Rechardson(1911)]{rechardson1911approximate}
LF~Rechardson.
\newblock The approximate arithmetical solution by finite difference of
  physical problems involving differential equations, with an application to
  the stresses in a masonary dam.
\newblock \emph{R. Soc. London Phil. Trans. A}, 210:\penalty0 307--357, 1911.

\bibitem[Recht et~al.(2011)Recht, Re, Wright, and Niu]{recht2011hogwild}
Benjamin Recht, Christopher Re, Stephen Wright, and Feng Niu.
\newblock Hogwild!: {A} lock-free approach to parallelizing stochastic gradient
  descent.
\newblock \emph{Advances in Neural Information Processing Systems}, 24, 2011.

\bibitem[Reddi et~al.(2021)Reddi, Charles, Zaheer, Garrett, Rush,
  Kone{\v{c}}n{\`y}, Kumar, and McMahan]{reddi2021adaptive}
Sashank Reddi, Zachary Charles, Manzil Zaheer, Zachary Garrett, Keith Rush,
  Jakub Kone{\v{c}}n{\`y}, Sanjiv Kumar, and H~Brendan McMahan.
\newblock Adaptive federated optimization.
\newblock \emph{International Conference on Learning Representations}, 2021.

\bibitem[Richt{\'a}rik \& Tak{\'a}c(2020)Richt{\'a}rik and
  Tak{\'a}c]{richtarik2020stochastic}
Peter Richt{\'a}rik and Martin Tak{\'a}c.
\newblock Stochastic reformulations of linear systems: algorithms and
  convergence theory.
\newblock \emph{SIAM Journal on Matrix Analysis and Applications}, 41\penalty0
  (2):\penalty0 487--524, 2020.

\bibitem[Richt{\'a}rik et~al.(2021)Richt{\'a}rik, Sokolov, and
  Fatkhullin]{richtarik2021ef21}
Peter Richt{\'a}rik, Igor Sokolov, and Ilyas Fatkhullin.
\newblock Ef21: A new, simpler, theoretically better, and practically faster
  error feedback.
\newblock \emph{Advances in Neural Information Processing Systems},
  34:\penalty0 4384--4396, 2021.

\bibitem[Robbins \& Monro(1951)Robbins and Monro]{robbins1951stochastic}
Herbert Robbins and Sutton Monro.
\newblock A stochastic approximation method.
\newblock \emph{The Annals of Mathematical Statistics}, pp.\  400--407, 1951.

\bibitem[Rockafellar(1976)]{rockafellar1976monotone}
R~Tyrrell Rockafellar.
\newblock Monotone operators and the proximal point algorithm.
\newblock \emph{SIAM Journal on Control and Optimization}, 14\penalty0
  (5):\penalty0 877--898, 1976.

\bibitem[Ryu \& Boyd(2014)Ryu and Boyd]{ryu2014stochastic}
Ernest~K Ryu and Stephen Boyd.
\newblock Stochastic proximal iteration: a non-asymptotic improvement upon
  stochastic gradient descent.
\newblock \emph{Author website, early draft}, 2014.

\bibitem[Sadiev et~al.(2022{\natexlab{a}})Sadiev, Kovalev, and
  Richt{\'a}rik]{sadiev2022communication}
Abdurakhmon Sadiev, Dmitry Kovalev, and Peter Richt{\'a}rik.
\newblock Communication acceleration of local gradient methods via an
  accelerated primal-dual algorithm with an inexact prox.
\newblock \emph{Advances in Neural Information Processing Systems},
  35:\penalty0 21777--21791, 2022{\natexlab{a}}.

\bibitem[Sadiev et~al.(2022{\natexlab{b}})Sadiev, Malinovsky, Gorbunov,
  Sokolov, Khaled, Burlachenko, and Richt{\'a}rik]{sadiev2022federated}
Abdurakhmon Sadiev, Grigory Malinovsky, Eduard Gorbunov, Igor Sokolov, Ahmed
  Khaled, Konstantin Burlachenko, and Peter Richt{\'a}rik.
\newblock Federated optimization algorithms with random reshuffling and
  gradient compression.
\newblock \emph{arXiv preprint arXiv:2206.07021}, 2022{\natexlab{b}}.

\bibitem[Seide et~al.(2014)Seide, Fu, Droppo, Li, and Yu]{seide20141}
Frank Seide, Hao Fu, Jasha Droppo, Gang Li, and Dong Yu.
\newblock 1-bit stochastic gradient descent and its application to
  data-parallel distributed training of speech dnns.
\newblock In \emph{Interspeech}, volume 2014, pp.\  1058--1062. Singapore,
  2014.

\bibitem[T~Dinh et~al.(2020)T~Dinh, Tran, and Nguyen]{t2020personalized}
Canh T~Dinh, Nguyen Tran, and Josh Nguyen.
\newblock Personalized federated learning with {M}oreau envelopes.
\newblock \emph{Advances in Neural Information Processing Systems},
  33:\penalty0 21394--21405, 2020.

\bibitem[Tyurin \& Richt{\'a}rik(2024)Tyurin and
  Richt{\'a}rik]{tyurin2024dasha}
Alexander Tyurin and Peter Richt{\'a}rik.
\newblock {DASHA}: Distributed nonconvex optimization with communication
  compression and optimal oracle complexity.
\newblock In \emph{International Conference on Learning Representations}, 2024.

\end{thebibliography}
\bibliographystyle{iclr2025_conference}

\newpage
\appendix
\tableofcontents
\section{Notations}
\label{sec:notations}
Throughout the paper, we use the notation $\norm{\cdot}$ to denote the standard Euclidean norm defined on $\R^d$ and $\inner{\cdot}{\cdot}$ to denote the standard Euclidean inner product.
Given a differentiable function $f: \R^d \mapsto \R$, its gradient is denoted as $\nabla f(x)$. 
We use the notation $\bregman{f}{x}{y}$ to denote the Bregman divergence associated with a function $f:\R^d \mapsto \R$ between $x$ and $y$. 
The notation $\inf f$ is used to denote the minimum of a function $f: \R^d \mapsto \R$. 
We use $\ProxSub{\gamma \phi}{x}$ to denote the proximity operator of function $\phi:\R^d \mapsto \R$ with $\gamma > 0$ at $x \in \R^d$, and $\MoreauSub{\gamma}{\phi}{x}$ to denote the corresponding Moreau Envelope.
We denote the average of the Moreau envelope of each local objective $f_i$ by the notation ${M^{\gamma}}:\R^d \mapsto \R$. 
Specifically, we define $\M{x} = \frac{1}{n}\sum_{i=1}^n \MoreauSub{\gamma}{f}{x}$.
Note that $\M{x}$ has an implicit dependence on $\gamma$, its smoothness constant is denoted by $L_{\gamma}$.
We say an extended real-valued function $f: \R^d \mapsto \R\cup\cbrac{+\infty}$ is proper if there exists $x \in \R^d$ such that $f(x) < +\infty$.
We say an extended real-valued function $f: \R^d \mapsto \R\cup\cbrac{+\infty}$ is closed if its epigraph is a closed set.
We use the notation $\cE_k = \gamma M^{\gamma}\rbrac{x_k} - \gamma M^{\gamma}_{\inf}$ to denote the function value suboptimality of $\gamma M^{\gamma}$ at $x_k$, and $\Delta_k = \norm{x_k - x_\star}^2$ to denote the squared distance. 
The notation $\cO\rbrac{\cdot}$ is used to describe complexity while omitting constant factors, whereas $\tilde{\cO}\rbrac{\cdot}$ is used when logarithmic factors are also omitted.

\section{Facts and lemmas}
\label{sec:facts}
\begin{fact}[Young's inequality]
   \label{tech-fact:youngs}
   For any two vectors $x, y \in \R^d$, the following inequality holds,
   \begin{equation}
      \label{eq:youngs}
      \norm{x + y}^2 \leq 2\norm{x}^2 + 2\norm{y}^2.
   \end{equation}
\end{fact}

\begin{fact}[Property of convex smooth functions]
   \label{tech-lemma:smooth}
   Let $\phi: \R^d \mapsto \R$ be differentiable.
   The following statements are equivalent:
   \begin{enumerate}
      \item $\phi$ is convex and $L$-smooth.
      \item $0 \leq 2\bregman{\phi}{x}{y} \leq L\norm{x - y}^2$ for all $x, y\in \R^d$.
      \item $\frac{1}{L}\norm{\nabla \phi(x) - \nabla \phi(y)}^2 \leq 2\bregman{\phi}{x}{y}$ for all $x, y \in \R^d$.
   \end{enumerate}
   The notation $\bregman{\phi}{x}{y}$ denotes the Bregman divergence associate with $\phi$ at $x, y\in R^d$, defined as 
   \begin{equation*}
      \bregman{\phi}{x}{y} = \phi\rbrac{x} - \phi\rbrac{y} -\inner{\nabla \phi(y)}{x - y}.
   \end{equation*}
\end{fact}
The following two facts establish that the convexity and smoothness of a function $\phi: \R^d \mapsto \R$ ensure the convexity and smoothness of its Moreau envelope.
\begin{fact}[Convexity of Moreau envelope]
   \label{fact:1:convexity-Moreau}
   \citep[Theorem 6.55]{beck2017first} Let $\phi: \R^d \mapsto \R\cup\cbrac{+\infty}$ be a proper and convex function. 
   Then $\Moreau^{\gamma}_\phi$ is a convex function. 
\end{fact}

\begin{fact}[Smoothness of Moreau envelope]
   \label{fact:2:smoothness-Moreau}
   \citep[Lemma 4]{li2024power} Let $\phi: \R^d \mapsto \R$ be a convex and $L$-smooth function. Then $\Moreau^{\gamma}_{\phi}$ is $\frac{L}{1 + \gamma L}$-smooth.
\end{fact}

The following fact illustrates the relationship between the minimizer of a function $\phi$ and its Moreau envelope $\Moreau^{\gamma}_{\phi}$.
\begin{fact}[Minimizer equivalence]
   \label{fact:3:minimizer-equiv}
   \citep[Lemma 5]{li2024power} Let $\phi: \R^d \mapsto \R\cup\cbrac{+\infty}$ be a proper, closed and convex function. 
   Then for any $\gamma > 0$, $\phi$ and $\Moreau^{\gamma}_\phi$ has the same set of minimizers. 
\end{fact}

In our case, we assume each $f_i$ from (\ref{eq:prob-formulation}) is convex and $L_i$-smooth. 
Therefore by \Cref{fact:1:convexity-Moreau} and \Cref{fact:2:smoothness-Moreau}, we know that each $\Moreau^{\gamma}_{f_i}$ is also convex and $\frac{L_i}{1 + \gamma L_i}$-smooth.
This means that $M_{\gamma} = \frac{1}{n}\sum_{i=1}^{n}\Moreau^{\gamma}_{f_i}$ is also convex and smooth.
We denote its smoothness constant as $L_{\gamma}$, and the following fact provides a range for this constant.
\begin{fact}[Global convexity and smoothness]
   \label{fact:4:global:moreau}
   \citep[Lemma 7]{li2024power} Let each $f_i$ be proper, closed convex and $L_i$-smooth. 
   Then $M^{\gamma}$ is convex and $L_{\gamma}$-smooth with 
   \begin{equation*}
      \frac{1}{n^2}\sum_{i=1}^n \frac{L_i}{1 + \gamma L_i} \leq L_{\gamma} \leq \frac{1}{n}\sum_{i=1}^n \frac{L_i}{1 + \gamma L_i}.
   \end{equation*}
\end{fact}
The following fact establishes that the minimizer of $f$ and $M^{\gamma}$ are the same.
\begin{fact}[Global minimizer equivalence]
   \label{fact:5:global:equivalence}
   \citep[Lemma 8]{li2024power} If we let every $f_i: \R^d \mapsto \R\cup\cbrac{+\infty}$ be proper, closed and convex, then $f(x) = \frac{1}{n}\sum_{i=1}^{n}f_i(x)$ has the same set of minimizers and minimum as 
   \begin{equation*}
     \M{x} = \frac{1}{n}\sum_{i=1}^{n} \MoreauSub{\gamma}{f_i}{x},
   \end{equation*}
   if we are in the interpolation regime and $0 < \gamma < \infty$.
\end{fact}
The above fact demonstrates that running {\SGD} on the objective $\Moreau^{\gamma}$ will lead us to the correct destination, as the minimizers of $\Moreau^{\gamma}$ and $f$ are identical in our setting.
In problem (\ref{eq:prob-formulation}), if we assume that $f$ is strongly convex, then we have $M^{\gamma}$ satisfies the following star strong convexity inequality.
\begin{fact}[Star strong convexity]
   \label{fact:6:quadratic-growth}
   \citep[Lemma 11]{li2024power}
   Assume \Cref{asp:diff} (Differentiability), \Cref{asp:int-pl-rgm} (Interpolation Regime), \Cref{asp:cvx} (Individual convexity), \Cref{asp:smoothness} (Smoothness) and \Cref{asp:stn-cvx} (Global strong convexity) hold, then the convex function $\M{x}$ satisfies the following inequality,
   \begin{equation*}
      \M{x} - M^{\gamma}_{\inf} \geq \frac{\mu}{1 + \gamma L_{\max}} \cdot \frac{1}{2}\norm{x - x_\star}^2,
   \end{equation*}
   for any $x \in \R^d$ and a minimizer $x_\star$ of $\M{x}$.
\end{fact}
The above fact implies that the strong convexity of $f$ translates to the star strong convexity of $M^{\gamma}$.
Star strong convexity is also known as quadratic growth (QG) condition \citep{anitescu2000degenerate}.
In the case of a convex function, it is also known as optimal strong convexity \citep{liu2015asynchronous} and semi-strong convexity \citep{gong2014linear}. 
It is known that for a convex function satisfying quadratic growth condition, it also satisfies the Polyak-Lojasiewicz inequality \citep{polyak1964gradient} which is described by the following lemma.
Notice that since \Cref{alg:inexact-fedexprox-full} can be viewed as running {\SGD} with objective $\gamma M^{\gamma}$ and a fixed step size $\alpha_k = \alpha$, we describe the inequality based on $\gamma M^{\gamma}$ in the following lemma.
\begin{lemma}[PL-inequality]
   \label{lemma:PL-lemma}
   Let \Cref{asp:diff} (Differentiability), \Cref{asp:int-pl-rgm} (Interpolation Regime), \Cref{asp:cvx} (Individual convexity), \Cref{asp:smoothness} (Smoothness) and \Cref{asp:stn-cvx} (Global strong convexity) hold, then $\gamma M^{\gamma}\rbrac{x}$ satisfies the following Polyak-Lojasiewicz inequality,
   \begin{equation}
      \label{eq:PL-lemma}
      \norm{\gamma \nabla M^{\gamma}\rbrac{x}}^2 \geq 2\cdot \frac{\gamma\mu}{4\rbrac{1 + \gamma L_{\max}}}\rbrac{\gamma M^{\gamma}\rbrac{x} - \gamma M^{\gamma}_{\inf}},
   \end{equation}
   where $x \in \R^d$ is an arbitrary vector and $x_\star$ is a minimizer of $M^{\gamma}\rbrac{x}$.
\end{lemma}

\section{Theory of biased SGD}
\label{sec:bSGD}
For completeness, we provide the theory of biased {\SGD} we used to analyze our algorithm in this paper.
It is adapted from \citet{demidovich2024guide}, which offers a comprehensive study of various assumptions employed in the analysis of {\SGD} with biased gradient updates.
In addition, the authors introduced a new set of assumptions, referred to as the Biased ABC assumption, which are less restrictive than all previous assumptions.
The authors provided convergence guarantees for {\SGD} with biased gradient updates in the non-convex and convex setting.
Specifically, they considered the case of minimizing a function $f: \R^d \mapsto \R$,
\begin{equation*}
   \min_{x \in \R^d} f(x),
\end{equation*}
with 
\begin{equation}
   \label{eq:biased-SGD}
   x_{k+1} = x_k - \eta g(x_k), \tag{biased SGD}
\end{equation}
where $\eta > 0$ is the stepsize, $g(x_k)$ is a possibly stochastic and biased gradient estimator.
They introduced the biased ABC assumption, 
\begin{assumption}[Biased-ABC]
   \label{asp:biased-ABC}
   \citep[Assumption 9]{demidovich2024guide} There exists constants $A, B, C, b, c \geq 0$ such that the gradient estimator $g(x)$ for every $x \in \R^d$ satisfies 
   \begin{eqnarray*}
      \inner{\nabla f(x)}{\Exp{g(x)}} &\geq& b\norm{\nabla f(x)}^2 - c \\
      \Exp{\norm{g(x)}^2} &\leq& 2A\rbrac{f(x) - f_{\inf}} + B\norm{\nabla f(x)}^2 + C.
   \end{eqnarray*}
\end{assumption}
A convergence guarantee was provided for \ref{eq:biased-SGD} under \Cref{asp:biased-ABC} given that $f$ is $\widehat{L}$-smooth and $\widehat{\mu}$-PL, that is, there exists $\widehat{\mu} > 0$, such that 
\begin{align*}
   \norm{\nabla f(x)}^2 \geq 2\widehat{\mu} \rbrac{f(x) - f_{\inf}},
\end{align*}
for all $x \in \R^d$.
\begin{theorem}[Theory of biased {\SGD}]
   \label{thm:biased-abc-pl}
   \citep[Theorem 4]{demidovich2024guide}
   Let $f$ be $\widehat{L}$-smooth and $\widehat{\mu}$-PL and \Cref{asp:biased-ABC} hold.
   If we choose a step size $\eta$ satisfying
   \begin{align}
      \label{eq:range-eta-theorem}
      0 < \eta < \min\cbrac{\frac{\widehat{\mu}b}{\widehat{L}\rbrac{A + \widehat{\mu}B}}, \frac{1}{\widehat{\mu}b}}.
   \end{align}
   Then we have 
   \begin{align*}
      \Exp{f(x_k) - f_{\inf}} \leq \rbrac{1 - \eta\widehat{\mu}b}^k\rbrac{f(x_0) - f_{\inf}} + \frac{LC\eta}{2\widehat{\mu}b} + \frac{c}{\widehat{\mu}b}.
   \end{align*}
   Under the special case of 
   \begin{equation*}
      \frac{\widehat{\mu}b}{\widehat{L}\rbrac{A + \widehat{\mu}B}} < \frac{1}{\widehat{\mu}b}, 
   \end{equation*}
   The range of the step size can be simplified to 
   \begin{align*}
      0 < \eta \leq \frac{\widehat{\mu}b}{\widehat{L}\rbrac{A + \widehat{\mu}B}},
   \end{align*}
   and if we take the largest possible step size, we have 
   \begin{align*}
      \Exp{f(x_k) - f_{\inf}} \leq \rbrac{1 - \frac{\widehat{\mu}^2b^2}{\widehat{L}\rbrac{A + \widehat{\mu}B}}}^k\rbrac{f(x_0) - f_{\inf}} + \frac{LC}{2\widehat{L}\rbrac{A + \widehat{\mu}B}} + \frac{c}{\widehat{\mu}b}.
   \end{align*}
\end{theorem}
The constants $C, c$ determine whether the algorithm is converging to the exact solution or just a neighborhood.
For $g(x) = \nabla f(x)$, clearly we have $A = 0, B = 1, b = 1, C =0, c =0$, and there is no neighborhood.
This is expected because the algorithm reduces to standard {\GD}
The iteration complexity is give by $\tilde{\cO}\rbrac{\frac{\widehat{L}}{\widehat{\mu}}}$, which is also expected for {\GD}.

\section{Theory of biased compression}
\label{sec:bcomp}
In this section, we present the theory of {\SGD} with biased compression.
The theory is adapted from \citet{beznosikov2023biased}.
The authors introduced theory for analyzing compressed gradient descent ({\CGD}) with biased compressor, both in the single node case and in the distributed case when the objective function is assumed to be strongly convex.
Here, we are only concerned with the single node case because distributed compressed gradient descent ({\DCGD}) with biased compressor may fail to converge.
To address this issue, error feedback mechanism \citep{seide20141,karimireddy2019error,richtarik2021ef21} is needed. 
In the single node case, the authors considered solving 
\begin{align*}
   \min_{x \in \R^d} f(x),
\end{align*}
where $f: \R^d \mapsto \R$ is $\widehat{L}$-smooth and $\widehat{\mu}$-strongly convex, with the following compressed gradient descent algorithm
\begin{align}
   \label{eq:biased-CGD}
   x_{k+1} = x_k - \eta \cC\rbrac{\nabla f(x_k)}, \tag{CGD}
\end{align}
where $\cC: \R^d \mapsto \R$ are potentially biased compression operators, $\eta > 0$ is a step size. 
The author proved that if certain conditions on $\cC$ is satisfied, a corresponding convergence guarantee can then be established.
Three classes of compressor/mapping were introduced.
\begin{definition}[Class $\bB^1$]
   \label{def:b1}
   We say a mapping $\cC \in \bB^1\rbrac{\alpha, \beta}$ for some $\alpha, \beta > 0$ if 
   \begin{align*}
      \alpha \norm{x}^2 \leq \Exp{\norm{\cC\rbrac{x}}^2} \leq \beta\inner{\Exp{\cC\rbrac{x}}}{x}, \qquad \forall x\in \R^d.
   \end{align*}
\end{definition}
\begin{definition}[Class $\bB^2$]
   \label{def:b2}
   We say a mapping $\cC \in \bB^2\rbrac{\xi, \beta}$ for some $\xi, \beta > 0$ if 
   \begin{align*}
      \max\cbrac{\xi\norm{x}^2, \frac{1}{\beta}\Exp{\norm{\cC\rbrac{x}}^2}} \leq \inner{\Exp{\cC\rbrac{x}}}{x}, \qquad \forall x \in \R^d. 
   \end{align*} 
\end{definition}
\begin{definition}[Class $\bB^3$]
   \label{def:b3}
   We say a mapping $\cC \in \bB^3\rbrac{\delta}$ for some $\delta > 0$, if 
   \begin{align*}
      \Exp{\norm{\cC\rbrac{x} - x}^2} \leq \rbrac{1 - \frac{1}{\delta}}\norm{x}^2.
   \end{align*}
\end{definition}
The authors proved the following theorem about the convergence of the algorithm, the notation $\cF_k$ is used to denote $\Exp{f(x_k)} - f_{\inf}$, with $\cF_0 = f(x_0) - f_{\inf}$,
\begin{theorem}
   \label{thm:biased-compression}
   Let $\cC \in \bB^1\rbrac{\alpha, \beta}$.
   Then we have $\cF_k \leq \rbrac{1 - \nicefrac{\alpha}{\beta} \eta\widehat{\mu}\rbrac{2 - \eta\beta\widehat{L}}}\cF_{k-1}$, as long as $0 \leq \eta \leq \frac{2}{\beta \widehat{L}}$.
   If we choose $\eta = \frac{1}{\beta\widehat{L}}$, we have 
   \begin{equation}
      \label{eq:thmbc-r1}
      \cF_k \leq \rbrac{1 - \frac{\alpha}{\beta^2}\cdot\frac{\widehat{\mu}}{\widehat{L}}}^K\cF_0.
   \end{equation}
   Let $\cC \in \bB^2\rbrac{\xi, \beta}$.
   Then we have $\cF_k \leq \rbrac{1 - \xi\eta\rbrac{2 - \eta\beta}\widehat{L}}\cF_{k-1}$, as long as $0 \leq \eta \leq \frac{2}{\beta\widehat{L}}$.
   If we choose $\eta = \frac{1}{\beta\widehat{L}}$, we have 
   \begin{equation}
      \label{eq:thmbc-r2}
      \cF_k \leq \rbrac{1 - \frac{\xi}{\beta}\cdot\frac{\widehat{\mu}}{\widehat{L}}}^k\cF_0.
   \end{equation}
   Let $\cC \in \bB^3\rbrac{\delta}$.
   Then we have $\cF_k \leq \rbrac{1 - \frac{1}{\delta}\eta\widehat{\mu}} \cF_{k-1}$, as long as $0 \leq \eta \leq \frac{1}{\widehat{L}}$.
   If we choose $\eta = \frac{1}{\widehat{L}}$, we have 
   \begin{equation}
      \label{eq:thmbc-r3}
      \cF_k \leq \rbrac{1 - \frac{1}{\delta} \cdot \frac{\widehat{\mu}}{\widehat{L}}}^k \cF_0.
   \end{equation} 
\end{theorem}
Notice that when $\cC\rbrac{x} = x$, that is, when no compression happens, we have $\alpha = \beta = \xi = \delta = 1$.
In this case, the iteration complexity of \ref{eq:biased-CGD} is given by $\tilde{\cO}\rbrac{\frac{\widehat{L}}{\widehat{\mu}}}$ and we recover the result of {\GD}.
It is worth noting that \Cref{thm:biased-compression} remains valid if the condition of $f$ being $\widehat{\mu}$-strongly convex is replaced with $f$ being $\widehat{\mu}$-PL.

\section{\texorpdfstring{Analysis of inexact {\FEDEXPROX} in the client sampling setting}{Analysis of inexact FedExProx in the client sampling setting}}
\label{sec:minibatch}
In this section, we will discuss the case where we do client sampling in \cref{alg:inexact-fedexprox-full}, we first formulate the algorithm as below.
For the sake of simplicity, we use $\tau$-nice sampling as an example.
\begin{algorithm}[h]
	\caption{Inexact {\FEDEXPROX} with $\tau$-nice sampling}
	\label{alg:inexact-fedexprox-mini}
	\begin{algorithmic}[1]
	\STATE {\bf Parameters:} extrapolation parameter $\alpha_k = \alpha > 0$, step size for the proximal operator $\gamma > 0$, starting point $x_0 \in \R^d$, number of clients $n$, size of minibatch $\tau$, total number of iterations $K$, proximal solution accuracy $\varepsilon_2 \geq 0$.
	\FOR{$k=0,1,2\dotsc K-1$}   
    	\STATE The server broadcasts the current iterate $x_k$ to a selected set of client $S_k$ of size $\tau$
        \STATE Each selected client computes a $\varepsilon$ approximation of the solution $\tx_{i, k+1} \simeq \ProxSub{\gamma f_i}{x_k}$, and sends it back to the server
        \STATE The server computes 
        \begin{equation}
            \label{eq:update-inexact-2}
            x_{k+1} = x_k + \alpha_k\rbrac{\frac{1}{\tau}\sum_{i\in S_k} \tx_{i, k+1} - x_k}.
        \end{equation} 
	\ENDFOR
	\end{algorithmic}
\end{algorithm}

\subsection{Relative approximation in distance}

\paragraph{The failure of biased compression theory:}
Similar to \Cref{thm:biased-compression}, we initially apply the theory from \citet{beznosikov2023biased}, as it provides improved results in the full-batch scenario.
We first define the compressing mapping $\cC_{\tau}$ in this case, 
\begin{equation}
   \label{eq:comp-mini-single}
   \cC_{\tau}\rbrac{\gamma\nabla M^{\gamma}\rbrac{x_k}} = \frac{1}{\tau}\sum_{i \in S_k}\rbrac{\gamma\nabla \MoreauSub{\gamma}{f_i}{x_k} - \rbrac{\tx_{i, k+1} - \ProxSub{\gamma f_i}{x_k}}}.
\end{equation}
One can verify for every $x_k$ and $\varepsilon_2$-approximation ${\tx_{i, k+1}}$ of $\ProxSub{\gamma f_i}{x_k}$, we have 
\begin{equation*}
   \cC_{\tau} \in \bB^3\rbrac{\delta = \frac{\mu}{\mu - 4\varepsilon_2L_{\max} - \frac{n - \tau}{\tau\rbrac{n - 1}}\sbrac{4\rbrac{2 + \varepsilon_2}L_{\max} -2\mu}}}
\end{equation*}
In the case of $\tau = n$, we have $\cC_{n} \in \bB^3\rbrac{\frac{\mu}{\mu - 4\varepsilon_2L_{\max}}}$, which recovers the result of (\ref{eq:check-verify}).
When $\tau = 1, \varepsilon_2 = 0$, however, this is problematic, as $\cC_1\in \bB^3\rbrac{\delta = \frac{\mu}{3\mu - 8L_{\max}}}$. 
Notice that we require $\delta > 0$, so we require $3\mu > 8L_{\max}$ which only holds in a very restrictive setting.
This is due to the stochasticity contained in (\ref{eq:comp-mini-single}), which arises from client sampling.

\paragraph{Theory of biased {\SGD}:}
The algorithm does converge, however, and one can use the theory of \citet{demidovich2024guide} to obtain a convergence guarantee.
\begin{theorem}
   \label{thm:020-mini}
   Assume \Cref{asp:diff} (Differentiability), \Cref{asp:int-pl-rgm} (Interpolation regime), \Cref{asp:cvx} (Individual convexity), \Cref{asp:smoothness} (Smoothness) and \Cref{asp:stn-cvx} (Global strong convexity) hold.
   Let the approximation $\tx_{i, k+1}$ all satisfies \Cref{def:inexact-2} with $\varepsilon_2 < \frac{\mu^2}{4L^2_{\max}}$, that is 
   \begin{equation*}
      \norm{\tx_{i, k+1} - \ProxSub{\gamma f_i}{x_k}}^2 \leq \varepsilon_2 \cdot \norm{x_k - \ProxSub{\gamma f_i}{x_k}}^2,
   \end{equation*}
   holds for all client $i$ at iteration $k$.
   If we are running \Cref{alg:inexact-fedexprox-mini} with minibatch size $\tau$ and extrapolation parameter $\alpha_k = \alpha > 0$ satisfying
   \begin{align*}
      \alpha \leq \frac{1}{\gamma L_{\gamma}} \cdot \frac{\mu - 2\sqrt{\varepsilon_2}L_{\max}}{\mu + 4\varepsilon_2 L_{\max} + 4\sqrt{\varepsilon_2}L_{\max} + \frac{n-\tau}{\tau\rbrac{n -1}}\cdot\rbrac{4L_{\max} + 4\sqrt{\varepsilon_2}L_{\max} - \mu}}
   \end{align*}
   Then the iterates generated by \Cref{alg:inexact-fedexprox-mini} satisfies 
   \begin{equation}
      \label{eq:thm020-guarantee-MINI}
      \Exp{\cE_K} \leq \rbrac{1 - \alpha\cdot\frac{\gamma\rbrac{\mu - 2\sqrt{\varepsilon_2} L_{\max}}}{4\rbrac{1 + \gamma L_{\max}}}}^K \cE_0.
   \end{equation}
   Specifically, if we choose the largest $\alpha$ possible, we have 
   \begin{align*}
      \Exp{\Delta_K} \leq \rbrac{1 - \frac{\mu}{4L_{\gamma}\rbrac{1 + \gamma L_{\max}}}\cdot S\rbrac{\varepsilon_2, \tau}}^K \cdot\frac{L{\gamma}\rbrac{1 + \gamma L_{\max}}}{\mu} \Delta_0,
   \end{align*}
   where $S\rbrac{\varepsilon_2, \tau}$ is defined as 
   \begin{align*}
      S\rbrac{\varepsilon_2, \tau} \eqdef \frac{\rbrac{\mu - 2\sqrt{\varepsilon_2}L_{\max}}\rbrac{1 - 2\sqrt{\varepsilon_2}\frac{L_{\max}}{\mu}}}{\mu + 4\varepsilon_2 L_{\max} + 4\sqrt{\varepsilon_2}L_{\max} + \frac{n-\tau}{\tau\rbrac{n -1}}\cdot\rbrac{4L_{\max} + 4\sqrt{\varepsilon_2}L_{\max} - \mu}},
   \end{align*}
   satisfying 
   \begin{align*}
      0 < S\rbrac{\varepsilon_2, \tau} \leq 1.
   \end{align*}
\end{theorem}
Notice that we have $S\rbrac{\varepsilon_2, \tau=n} = S\rbrac{\varepsilon_2}$, which appears in \Cref{thm:020-full}.
For the special case when $\varepsilon_2 = 0$, every proximal operator is solved exactly.
The range of $\alpha$ becomes,
\begin{equation*}
   0 < \alpha \leq \frac{1}{\gamma L_{\gamma}} \cdot \frac{\mu}{\frac{n - \tau}{\tau\rbrac{n - 1}} \cdot 4L_{\max} + \frac{n\rbrac{\tau - 1}}{\tau\rbrac{n - 1}}\mu}.
\end{equation*} 
According to \citet{li2024power}, 
\begin{equation*}
   0 < \alpha \leq \frac{1}{\gamma L_{\gamma}}\cdot\frac{L_{\gamma}\rbrac{1 + \gamma L_{\max}}}{\frac{n - \tau}{\tau\rbrac{n - 1}}L_{\max} + \frac{n\rbrac{\tau - 1}}{\tau\rbrac{n - 1}}\cdot L_{\gamma}\rbrac{1 + \gamma L_{\max}}}.
\end{equation*}
Clearly the bound we obtain here is suboptimal, since we have $\mu \leq L_{\gamma}\rbrac{1 + \gamma L_{\max}}$ according to (\ref{eq:mu-est}). 
This is due to the previously mentioned issue: the nature of biased compression.
When client sampling is used together with biased compressors, it does not necessarily guarantee any benefits.
To solve this, the modification of the algorithm itself may be needed, which we consider as a future work direction.

\subsection{Absolute approximation in distance}
Similarly to \Cref{thm:020-mini}, by applying the theory of biased {\SGD} \citep{demidovich2024guide}, we can derive a convergence guarantee for the minibatch case, though with a suboptimal convergence rate. For brevity and clarity, we do not include the details here.

\section{Proof of theorems and lemmas}

\subsection{\texorpdfstring{Proof of \Cref{lemma:PL-lemma}}{Proof of Lemma~\ref{lemma:PL-lemma}}}
Using \Cref{fact:6:quadratic-growth}, we have 
\begin{equation}
   \label{eq:ineq-pr1-1}
   M^{\gamma}\rbrac{x} - M^{\gamma}_{\inf} \geq \frac{\mu}{1 + \gamma L_{\max}}\cdot \frac{1}{2} \norm{x - x_\star}^2,
\end{equation}
where $x \in \R^d$ is any vector, $x_\star$ is a minimizer of $M^{\gamma}$, by \Cref{fact:3:minimizer-equiv}, it is also a minimizer of $f$.
Since we assume each function $f_i$ is convex, by \Cref{fact:1:convexity-Moreau}, we know that $\Moreau^{\gamma}_{f_i}$ is also convex.
As a result, the average of $\Moreau^{\gamma}_{f_i}$, $M^{\gamma}$ is also a convex function.
Utilizing the convexity of $M^{\gamma}$, we have, 
\begin{equation*}
   M^{\gamma}_{\inf} \geq M^{\gamma}\rbrac{x} + \inner{\nabla M^{\gamma}\rbrac{x}}{x_\star - x}.
\end{equation*}
Rearranging terms we get,
\begin{equation}
   \label{eq:ineq-pr1-2}
   \inner{\nabla M^{\gamma}\rbrac{x}}{x - x_\star} \geq M^{\gamma}\rbrac{x} - M^{\gamma}_{\inf}.
\end{equation}
As a result, we have 
\begin{equation*}
   \inner{\nabla M^{\gamma}\rbrac{x}}{x - x_\star} \overset{(\ref{eq:ineq-pr1-1}) + (\ref{eq:ineq-pr1-2})}{\geq} \frac{\mu}{1 + \gamma L_{\max}} \cdot \frac{1}{2}\norm{x - x_\star}^2.
\end{equation*}
Using Cauchy-Schwarz inequality, we have 
\begin{equation*}
   \norm{\nabla M^{\gamma}\rbrac{x}}\norm{x - x_\star} \geq \inner{\nabla M^{\gamma}\rbrac{x}}{x - x_\star} \geq \frac{\mu}{1 + \gamma L_{\max}} \cdot \frac{1}{2}\norm{x - x_\star}^2.
\end{equation*}
When $\norm{x - x_\star} > 0$, the above inequality leads to 
\begin{equation}
   \label{eq:ineq-pr1-4}
   \norm{\nabla M^{\gamma}\rbrac{x}} \geq \frac{\mu}{2\rbrac{1 + \gamma L_{\max}}} \cdot \norm{x - x_\star},
\end{equation}
which also holds when $\norm{x - x_\star} = 0$.
Now using (\ref{eq:ineq-pr1-2}) and (\ref{eq:ineq-pr1-4}), we obtain 
\begin{eqnarray*}
   M^{\gamma}\rbrac{x} - M^{\gamma}_{\inf} &\overset{(\ref{eq:ineq-pr1-2})}{\leq}& \inner{\nabla M^{\gamma}\rbrac{x}}{x - x_\star} \\
   &\leq& \norm{\nabla M^{\gamma}\rbrac{x}}\norm{x - x_\star} \\
   &\overset{(\ref{eq:ineq-pr1-4})}{\leq}& \frac{2\rbrac{1 + \gamma L_{\max}}}{\mu}\norm{\nabla M^{\gamma}\rbrac{x}}^2. 
\end{eqnarray*}
A simple rearranging of terms result in 
\begin{equation*}
   \norm{\gamma \nabla M^{\gamma}\rbrac{x}}^2 \geq 2\cdot \frac{\gamma\mu}{4\rbrac{1 + \gamma L_{\max}}}\rbrac{\gamma M^{\gamma}\rbrac{x} - \gamma M^{\gamma}_{\inf}}.
\end{equation*}
Up till here we have already proved the statement in the lemma, but we want to look at the strongly constant $\mu$ of $f$ a little bit.
In order to provide an upper bound of $\mu$, we notice that due to \Cref{fact:2:smoothness-Moreau}, each $\Moreau^{\gamma}_{f_i}$ is $\frac{L_i}{1 + \gamma L_i}$-smooth and therefore $M^{\gamma}$ is smooth. 
We use the notation $L_{\gamma}$ to denote its smoothness constant.
Applying the smoothness of $M^{\gamma}\rbrac{x}$, we have 
\begin{align*}
   M^{\gamma}\rbrac{x} \leq M^{\gamma}\rbrac{x_\star} + \inner{\nabla M^{\gamma}\rbrac{x_\star}}{x - x_\star} + \frac{L_{\gamma}}{2}\norm{x - x^\star}^2.
\end{align*}
Utilizing the fact that $\nabla M^{\gamma}\rbrac{x_\star} = 0$, we have 
\begin{align}
   \label{eq:comment-2}
   M^{\gamma}\rbrac{x} - M^{\gamma}_{\inf} &\leq \frac{L_{\gamma}}{2}\norm{x - x_\star}^2 
\end{align}
Combining (\ref{eq:comment-2}) and (\ref{eq:ineq-pr1-1}), we can deduce that 
\begin{equation*}
   \frac{\mu}{1 + \gamma L_{\max}}\cdot \frac{1}{2}\norm{x - x_\star}^2 \leq M^{\gamma}\rbrac{x} - M^{\gamma}_{\inf} \leq \frac{L_{\gamma}}{2}\norm{x - x_\star}^2.
\end{equation*}
which results in the estimate that 
\begin{equation}
   \label{eq:mu-est}
   \mu \leq L_{\gamma}\rbrac{1 + \gamma L_{\max}}.
\end{equation}

\subsection{\texorpdfstring{Proof of \Cref{thm:1:conv-full-batch-stncvx}}{Proof of Theorem~\ref{thm:1:conv-full-batch-stncvx}}}
Let us first recall that after reformulation, \Cref{alg:inexact-fedexprox-full} can be written as 
\begin{align*}
   x_{k+1} = x_k - \alpha \cdot g(x_k),
\end{align*}
where $g(x_k)$ is defined as 
\begin{align*}
   g(x_k) \eqdef \frac{1}{n}\sum_{i=1}^{n} \gamma\nabla\MoreauSub{\gamma}{f_i}{x_k} - \frac{1}{n}\sum_{i=1}^{n}\rbrac{\tx_{i, k+1} - \ProxSub{\gamma f_i}{x_k}}.
\end{align*}
We view this as running full batch biased {\SGD} with stepsize $\alpha$ and global objective $\gamma M^{\gamma}\rbrac{x}$.
We first examine if \Cref{asp:biased-ABC} (Biased-ABC) holds for arbitrary $x_k$.
Since we are in the full batch case, it is easy to see that 
\begin{align*}
   \Exp{g(x_k)} = g(x_k).
\end{align*}
Since our objective now is $\gamma M^{\gamma}\rbrac{x}$, we have that 
\begin{align*}
   \inner{ \gamma \nabla M^{\gamma}\rbrac{x_k}}{g(x_k)} &= \inner{\gamma\nabla M^{\gamma}\rbrac{x_k}}{\gamma\nabla M^{\gamma}\rbrac{x_k} - \frac{1}{n}\sum_{i=1}^n\rbrac{\tx_{i, k+1} - \ProxSub{\gamma f_i}{x_k}}} \\
   &= \norm{\gamma \nabla M^{\gamma}\rbrac{x_k}}^2 - \underbrace{\inner{\gamma \nabla M^{\gamma}\rbrac{x_k}}{\frac{1}{n}\sum_{i=1}^n\rbrac{\tx_{i, k+1} - \ProxSub{\gamma f_i}{x_k}}}}_{\eqdef P_1}.
\end{align*}
Now let us focus on $P_1$, we have the following upper bound,
\begin{eqnarray*}
   P_1 &\leq& \frac{1}{2}\norm{\gamma \nabla M^{\gamma}\rbrac{x_k}}^2 + \frac{1}{2}\norm{\frac{1}{n}\sum_{i=1}^{n}\rbrac{\tx_{i, k+1} - \ProxSub{\gamma f_i}{x_k}}}^2 \\
   &\overset{(\ref{eq:approx-dist})}{\leq}& \frac{1}{2}\norm{\gamma \nabla M^{\gamma}\rbrac{x_k}}^2 + \frac{\varepsilon_1}{2}.
\end{eqnarray*}
As a result, we have 
\begin{align*}
   \inner{\gamma \nabla M^{\gamma}\rbrac{x_k}}{g(x_k)} \geq \frac{1}{2}\norm{\gamma \nabla M^{\gamma}\rbrac{x_k}} - \frac{\varepsilon_1}{2}, 
\end{align*}
which holds for arbitrary $x_k$.
This suggests that $b = \frac{1}{2}, c = \frac{\varepsilon_1}{2}$.
On the other hand,
\begin{eqnarray*}
   \Exp{\norm{g(x_k)}^2} &=& \norm{\gamma\nabla M^{\gamma}\rbrac{x_k} + \frac{1}{n}\sum_{i=1}^n \rbrac{\tx_{i, k+1} - \ProxSub{\gamma f_i}{x_k}}}^2 \\
   &\overset{(\ref{eq:youngs})}{\leq}& 2\norm{\gamma \nabla M^{\gamma}\rbrac{x_k}}^2 + 2\norm{\frac{1}{n}\sum_{i=1}^n \rbrac{\tx_{i, k+1} - \ProxSub{\gamma f_i}{x_k}}}^2 \\
   &\overset{(\ref{eq:approx-dist})}{\leq}& 2\norm{\gamma \nabla M^{\gamma}\rbrac{x_k}}^2 + 2 \varepsilon_1.
\end{eqnarray*}
Thus, we can choose $A = 0, B=2, C = 2\varepsilon_1$.
Since we have assumed \Cref{asp:cvx} (Individual convexity) and \Cref{asp:smoothness} (Smoothness), it is easy to see that $M^{\gamma}$ is smooth, and we denote its smoothness constant as $L_{\gamma}$.
It is therefore straightforward to see that our global objective $\gamma M^{\gamma}$ is $\gamma L_{\gamma}$-smooth.
We also assume $f$ is $\mu$-strongly convex, which by \Cref{fact:6:quadratic-growth} indicates that $M^{\gamma}$ is $\frac{\mu}{1 + \gamma L_{\max}}$ star strongly convex.
We immediately obtain using \Cref{lemma:PL-lemma} that $\gamma M^{\gamma}$ is $\frac{\gamma\mu}{4\rbrac{1 + \gamma L_{\max}}}$-PL.
Now, we have validated all the assumptions for using \Cref{thm:biased-abc-pl}.
Applying \Cref{thm:biased-abc-pl}, we obtain that when the extrapolation parameter satisfies 
\begin{equation*}
   0 < \alpha < \frac{1}{4} \cdot \min \cbrac{\frac{1}{\gamma L_{\gamma}}, \frac{2\rbrac{1 + \gamma L_{\max}}}{\gamma\mu}}, 
\end{equation*}
the last iterate $x_K$ of \Cref{alg:inexact-fedexprox-full} with each proximal operator solved inexactly according to \Cref{def:1:prox} satisfies 
\begin{align*}
   \cE_K &\leq \rbrac{1 - \frac{\alpha\gamma\mu}{8\rbrac{1 + \gamma L_{\max}}}}^K\cE_0 + \frac{8\varepsilon_1 \alpha L_{\gamma}\rbrac{1 + \gamma L_{\max}}}{\mu} + \frac{4\varepsilon_1\rbrac{1 + \gamma L_{\max}}}{\gamma\mu},
\end{align*}
where $\cE_k = \gamma M^{\gamma}\rbrac{x_k} - M^{\gamma}_{\inf}$.
Let us now prove that  
\begin{align*}
   \frac{1}{\gamma L_{\gamma}} < \frac{2\rbrac{1 + \gamma L_{\max}}}{\gamma\mu}.
\end{align*}
This is equivalent to prove 
\begin{align*}
   \mu < 2L_{\gamma}\rbrac{1 + \gamma L_{\max}},
\end{align*}
which is always true since (\ref{eq:mu-est}) holds.
As a result, we can simplify the range of the extrapolation parameter to 
\begin{align*}
   0 < \alpha \leq \frac{1}{4\gamma L_{\gamma}}.
\end{align*}
If we pick the largest possible $\alpha$, we have 
\begin{align*}
   \cE_K &\leq \rbrac{1 - \frac{\mu}{32L_{\gamma}\rbrac{1 + \gamma L_{\max}}}}^K \cE_0 + \frac{6\varepsilon_1\rbrac{1 + \gamma L_{\max}}}{\gamma\mu}. 
\end{align*}
This result is not directly comparable to that of \citet{li2024power}.
However, using smoothness of $\gamma L_{\gamma}$, if we denote $\Delta_k = \norm{x_k - x_\star}^2$ where $x_\star$ is a minimizer of both $M^{\gamma} $ and $f$ since we assume we are in the interpolation regime (\Cref{asp:int-pl-rgm}), we have 
\begin{align*}
   \cE_0 \leq \frac{\gamma L_{\gamma}}{2}\Delta_0.
\end{align*} 
Using star strong convexity, we have 
\begin{align*}
   \cE_K \geq \frac{\gamma\mu}{2\rbrac{1 + \gamma L_{\max}}} \Delta_K.
\end{align*}
As a result, we can transform the above convergence guarantee into 
\begin{align*}
   \Delta_K \leq \rbrac{1 - \frac{\mu}{32L_{\gamma}\rbrac{1 + \gamma L_{\max}}}}^K\frac{L_{\gamma}\rbrac{1 + \gamma L_{\max}}}{\mu} \cdot \Delta_0 + 12\varepsilon_1\cdot\rbrac{\frac{\nicefrac{1}{\gamma} + L_{\max}}{\mu}}^2.
\end{align*}
This completes the proof.

\subsection{\texorpdfstring{Proof of \Cref{thm:020-full}}{Proof of Theorem~\ref{thm:020-full}}}
Since we based our analysis on the theory of biased {\SGD}, we first verify the validity of \Cref{asp:biased-ABC}.

\paragraph{Finding $b$ and $c$:}
Let us start with finding a lower bound on $\inner{\gamma \nabla M^{\gamma}\rbrac{x_k}}{\Exp{g(x_k)}}$.
We have 
\begin{align*}
   \inner{\gamma M^{\gamma}\rbrac{x_k}}{\Exp{g(x_k)}} &= \inner{\gamma M^{\gamma}\rbrac{x_k}}{\gamma M^{\gamma}\rbrac{x_k} - \frac{1}{n}\sum_{i=1}^{n}\rbrac{\tx_{i, k+1} - \ProxSub{\gamma f_i}{x_k}}} \\
   &= \norm{\gamma M^{\gamma}\rbrac{x_k}}^2 - \inner{\gamma M^{\gamma}\rbrac{x_k}}{\frac{1}{n}\sum_{i=1}^n\rbrac{\tx_{i, k+1} - \ProxSub{\gamma f_i}{x_k}}} \\
   &\geq \norm{\gamma M^{\gamma}\rbrac{x_k}}^2 - \norm{\gamma M^{\gamma}\rbrac{x_k}} \cdot \norm{\frac{1}{n}\sum_{i=1}^n\rbrac{\tx_{i, k+1} - \ProxSub{\gamma f_i}{x_k}}},
\end{align*}
where the last inequality is obtained using Cauchy-Schwarz inequality.
We then utilize the convexity of $\norm{\cdot}$ and obtain, 
\begin{eqnarray*}
   \inner{\gamma M^{\gamma}\rbrac{x_k}}{\Exp{g(x_k)}} &\geq& \norm{\gamma M^{\gamma}\rbrac{x_k}}^2 - \norm{\gamma M^{\gamma}\rbrac{x_k}} \cdot \frac{1}{n}\sum_{i=1}^n\norm{\rbrac{\tx_{i, k+1} - \ProxSub{\gamma f_i}{x_k}}} \\
   & \overset{(\ref{eq:index2-approx})}{\geq}& \norm{\gamma M^{\gamma}\rbrac{x_k}}^2 - \sqrt{\varepsilon_2}\norm{\gamma M^{\gamma}\rbrac{x_k}}\cdot\frac{1}{n}\sum_{i=1}^{n} \norm{x_k - \ProxSub{\gamma f_i}{x_k}} \\
   &=& \norm{\gamma M^{\gamma}\rbrac{x_k}}^2 - \sqrt{\varepsilon_2}\norm{\gamma M^{\gamma}\rbrac{x_k}}\cdot\frac{1}{n}\sum_{i=1}^{n} \norm{\gamma\nabla \MoreauSub{\gamma}{f_i}{x_k}}.
\end{eqnarray*}
Notice that 
\begin{align*}
   \norm{\gamma\nabla \MoreauSub{\gamma}{f_i}{x_k}} = \norm{\gamma\nabla \MoreauSub{\gamma}{f_i}{x_k} - \gamma\nabla \MoreauSub{\gamma}{f_i}{x_\star}},
\end{align*}
holds for any $x_\star$ that is a minimizer of $M^{\gamma}\rbrac{x}$ due to interpolation regime assumption.
As a result, we can provide an upper bound based on smoothness of each individual $\gamma \MoreauSub{\gamma}{f_i}{x}$ using \Cref{tech-lemma:smooth},
\begin{align}
   \label{eq:lemma-5t0}
   \norm{\gamma\nabla \MoreauSub{\gamma}{f_i}{x_k} - \gamma\nabla \MoreauSub{\gamma}{f_i}{x_\star}} \leq \frac{\gamma L_i}{1 + \gamma L_i}\norm{x_k - x_\star}. 
\end{align}
Thus, 
\begin{align*}
   \frac{1}{n}\sum_{i=1}^n \norm{\gamma\nabla \MoreauSub{\gamma}{f_i}{x_k}} &\leq \frac{1}{n}\sum_{i=1}^{n} \frac{\gamma L_i}{1 + \gamma L_i}\norm{x_k - x_\star} \leq \frac{\gamma L_{\max}}{1 + \gamma L_{\max}}\cdot \norm{x_k - x_\star}.
\end{align*}
In addition, we have due to Cauchy-Schwarz inequality and the convexity of $M^{\gamma}\rbrac{x}$ 
\begin{equation}
   \label{eq:lemma-5t1}
   \norm{\nabla M^{\gamma}\rbrac{x_k}}\cdot \norm{x_k - x_\star} \geq \inner{\nabla M^{\gamma}\rbrac{x_k}}{x_k - x_\star} \geq M^{\gamma}\rbrac{x_k} - M^{\gamma}_{\inf},
\end{equation}
and due to quadratic growth condition that 
\begin{equation}
   \label{eq:lemma-5t2}
   M^{\gamma}\rbrac{x_k} - M^{\gamma}_{\inf} \geq \frac{\mu}{1 + \gamma L_{\max}} \cdot \frac{1}{2}\norm{x_k - x_\star}^2.
\end{equation}
Combining (\ref{eq:lemma-5t1}) and (\ref{eq:lemma-5t2}), we have 
\begin{align*}
   \frac{\mu}{2\rbrac{1 + \gamma L_{\max}}} \cdot \norm{x_k - x_\star}^2 \overset{(\ref{eq:lemma-5t1}) + (\ref{eq:lemma-5t2})}{\leq} \norm{\nabla M^{\gamma}\rbrac{x_k}}\cdot\norm{x_k - x_\star}.
\end{align*}
This indicates that 
\begin{align}
   \label{eq:lemma-5t3}
   \norm{x_k - x_\star} \leq \frac{2\rbrac{1 + \gamma L_{\max}}}{\mu}\norm{\nabla M^{\gamma}\rbrac{x_k}}.
\end{align}
Combining (\ref{eq:lemma-5t0}) and (\ref{eq:lemma-5t3}), we generate the following lower bound 
\begin{align*}
   \inner{\gamma M^{\gamma}\rbrac{x_k}}{\Exp{g(x_k)}} &\overset{(\ref{eq:lemma-5t0})}{\geq} \norm{\gamma M^{\gamma}\rbrac{x_k}}^2 - \sqrt{\varepsilon_2}\norm{\gamma M^{\gamma}\rbrac{x_k}} \cdot \frac{\gamma L_{\max}}{1 + \gamma L_{\max}}\norm{x_k - x_\star} \\
   &\overset{(\ref{eq:lemma-5t3})}{\geq} \norm{\gamma M^{\gamma}\rbrac{x_k}}^2 - \sqrt{\varepsilon_2} \cdot \frac{ L_{\max}}{1 + \gamma L_{\max}} \cdot \frac{2\rbrac{1 + \gamma L_{\max}}}{\mu}\norm{\gamma M^{\gamma}\rbrac{x_k}}^2 \\
   &= \rbrac{1 - \sqrt{\varepsilon_2}\cdot \frac{2L_{\max}}{\mu}} \cdot \norm{\gamma M^{\gamma}\rbrac{x_k}}^2.
\end{align*}
Thus, as long as $\varepsilon_2 < \frac{\mu^2}{4L_{\max}^2}$, we have $b = 1 - \sqrt{\varepsilon_2}\cdot \frac{2L_{\max}}{\mu}$, and $c = 0$.

\paragraph{Finding $A, B$ and $C$:}
We start with expanding $\norm{g(x_k)}^2$, 
\begin{align}
   \label{eq:lemma-6-t1}
   \Exp{\norm{g(x_k)}^2} &= \norm{\gamma M^{\gamma}\rbrac{x_k} - \frac{1}{n}\sum_{i=1}^{n}\rbrac{\tx_{i, k+1} - \ProxSub{\gamma f_i}{x_k}}}^2 \notag \\
   &= \norm{\gamma M^{\gamma}\rbrac{x_k}}^2 + \underbrace{\norm{\frac{1}{n}\sum_{i=1}^{n}\rbrac{\tx_{i, k+1} - \ProxSub{\gamma f_i}{x_k}}}^2}_{\eqdef T_2} \notag \\
   &\quad \underbrace{-2\inner{\gamma M^{\gamma}\rbrac{x_k}}{\frac{1}{n}\sum_{i=1}^{n}\rbrac{\tx_{i, k+1} - \ProxSub{\gamma f_i}{x_k}}}}_{\eqdef T_3}.
\end{align}
It is easy to bound $T_2$ utilizing the convexity of $\norm{\cdot}^2$,
\begin{align*}
   T_2 &\leq \frac{1}{n}\sum_{i=1}^{n}\norm{\tx_{i, k+1} - \ProxSub{\gamma f_i}{x_k}}^2 \\
   &\overset{(\ref{eq:index2-approx})}{\leq} \frac{\varepsilon_2}{n}\sum_{i=1}^n \norm{x_k - \ProxSub{\gamma f_i}{x_k}}^2 = \frac{\varepsilon_2}{n}\sum_{i=1}^n \norm{\gamma \MoreauSub{\gamma}{f_i}{x_k}}^2.
\end{align*}
Let $x_\star$ be a minimizer of $M^{\gamma}$, since we assume \Cref{asp:int-pl-rgm} holds, it is also a minimizer of each $\Moreau^{\gamma}_{f_i}$.
As a result, 
\begin{align}
   \label{eq:lemma-6-t2}
   T_2 &\leq \frac{\varepsilon_2}{n}\sum_{i=1}^{n} \norm{\gamma \MoreauSub{\gamma}{f_i}{x_k} - \gamma \MoreauSub{\gamma}{f_i}{x_\star}}^2 \notag \\
   &\leq \frac{\varepsilon_2}{n}\sum_{i=1}^n \frac{2\gamma L_i}{1 + \gamma L_i}\rbrac{\gamma M^{\gamma}_{f_i}\rbrac{x_k} - \gamma M^{\gamma}_{f_i}\rbrac{x_\star}} \leq \frac{2\varepsilon_2 \gamma L_{\max}}{1 + \gamma L_{\max}}\cdot\rbrac{\gamma M^{\gamma}\rbrac{x_k} - \gamma M^{\gamma}_{\inf}}.
\end{align}
We then consider $T_3$, and start with applying Cauchy-Schwarz inequality
\begin{align}
    \label{eq:lemma-6-t3-1}
    T_3 \leq 2\norm{\gamma \nabla M^{\gamma}\rbrac{x_k}}\norm{\frac1n\sum_{i=1}^{n}\rbrac{\tx_{i, k+1} - \ProxSub{\gamma f_i}{x_k}}}.
\end{align}
Using the convexity of $\norm{\cdot}$, we have 
\begin{eqnarray*}
   \norm{\frac1n\sum_{i=1}^{n}\rbrac{\tx_{i, k+1} - \ProxSub{\gamma f_i}{x_k}}} &\leq& \frac{1}{n}\sum_{i=1}^n \norm{\tx_{i, k+1} - \ProxSub{\gamma f_i}{x_k}} \\
   &\overset{(\ref{eq:index2-approx})}{\leq}& \frac{\sqrt{\varepsilon_2}}{n}\sum_{i=1}^{n} \norm{x_k - \ProxSub{\gamma f_i}{x_k}} \\
   &\overset{(\ref{eq:ppt-2-gradient})}{=}& \frac{\sqrt{\varepsilon_2}}{n}\sum_{i=1}^{n}\norm{\gamma\nabla\MoreauSub{\gamma}{f_i}{x_k} - \gamma\nabla\MoreauSub{\gamma}{f_i}{x_\star}} \\
   &\overset{\text{\Cref{tech-lemma:smooth}}}{\leq}& \frac{\sqrt{\varepsilon_2}}{n}\sum_{i=1}^n \frac{\gamma L_i}{1 + \gamma L_i}\norm{x_k - x_\star} \\
   &\leq& \frac{\sqrt{\varepsilon_2}\gamma L_{\max}}{1 + \gamma L_{\max}}\cdot\norm{x_k - x_\star}.
\end{eqnarray*}
Utilizing (\ref{eq:lemma-5t3}), we have 
\begin{align}
   \label{eq:lemma-6-t3}
   \norm{\frac1n\sum_{i=1}^{n}\rbrac{\tx_{i, k+1} - \ProxSub{\gamma f_i}{x_k}}} &\leq \frac{\sqrt{\varepsilon_2}\gamma L_{\max}}{1 + \gamma L_{\max}}\cdot \frac{2\rbrac{1 + \gamma L_{\max}}}{\mu}\norm{\nabla M^{\gamma}\rbrac{x_k}} \notag \\
   &= \frac{2\sqrt{\varepsilon_2}L_{\max}}{\mu}\cdot\norm{\gamma \nabla M^{\gamma}\rbrac{x_k}}
\end{align}
Plug the above inequality into (\ref{eq:lemma-6-t3-1}), we have 
\begin{align}
   \label{eq:lemma-6-t4}
   T_3 \leq \frac{4\sqrt{\varepsilon_2} L_{\max}}{\mu}\cdot \norm{\gamma \nabla M^{\gamma}\rbrac{x_k}}^2.
\end{align}
Combining (\ref{eq:lemma-6-t4}) and (\ref{eq:lemma-6-t2}), plug them into (\ref{eq:lemma-6-t1}), we have 
\begin{align*}
   \Exp{\norm{g\rbrac{x_k}}^2} \leq \frac{2\varepsilon_2\gamma L_{\max}}{1 + \gamma L_{\max}}\cdot \rbrac{\gamma M^{\gamma}\rbrac{x_k} - \gamma M^{\gamma}_{\inf}} + \rbrac{1 + \frac{4\sqrt{\varepsilon_2}L_{\max}}{\mu}}\cdot\norm{\gamma\nabla M^{\gamma}\rbrac{x_k}}^2.
\end{align*}
Thus, we have 
\begin{align*}
   A = \frac{\varepsilon_2\gamma L_{\max}}{1 + \gamma L_{\max}}, \quad B = \frac{\mu + 4\sqrt{\varepsilon_2} L_{\max}}{\mu}, \quad C=0.
\end{align*}

\paragraph{Applying \Cref{thm:biased-abc-pl}:}
First, we list our the values appeared respectively,
\begin{align*}
    & A = \frac{\varepsilon_2\gamma L_{\max}}{1 + \gamma L_{\max}}, \quad B = \frac{\mu + 4\sqrt{\varepsilon_2} L_{\max}}{\mu}, \quad b = \frac{\mu - 2\sqrt{\varepsilon_2}L_{\max}}{\mu}, \\
    & C = c = 0. 
\end{align*}
We know that the PL constant of $\gamma M^{\gamma}$ is given by $\frac{\gamma\mu}{4\rbrac{1 + \gamma L_{\max}}}$ and the corresponding smoothness constant is $\gamma L_{\gamma}$.
Applying \Cref{thm:biased-abc-pl}, the range of $\alpha$ is given by 
\begin{align}
   \label{eq:defB-1}
   0 < \alpha < \min \cbrac{\underbrace{\frac{1}{\gamma L_{\gamma}}\cdot\frac{\mu - 2\sqrt{\varepsilon_2}L_{\max}}{\mu + 4\sqrt{\varepsilon_2}L_{\max} + 4\varepsilon_2 L_{\max}}}_{\eqdef B_1}, \underbrace{\frac{4\rbrac{1 + \gamma L_{\max}}}{\gamma\rbrac{\mu - 2\sqrt{\varepsilon_2} L_{\max}}}}_{\eqdef B_2}}.
\end{align}
Now notice that actually we can prove that for $\varepsilon_2 < \frac{\mu^2}{4L_{\max}^2}$, we have $B_2 > B_1$, and we can simplify the range of $\alpha$ to 
\begin{align*}
   0 < \alpha \leq \frac{1}{\gamma L_{\gamma}}\cdot\frac{\mu - 2\sqrt{\varepsilon_2}L_{\max}}{\mu + 4\sqrt{\varepsilon_2}L_{\max} + 4\varepsilon_2 L_{\max}}.
\end{align*}

\paragraph{Proof of $B_2 > B_1$}:
It is easy to verify that the above inequality ($B_2 > B_1$) can be equivalently written as 
\begin{align*}
   4L_{\gamma}\rbrac{1 + \gamma L_{\max}}\rbrac{\mu + 4\sqrt{\varepsilon_2}L_{\max} + 4\varepsilon_2 L_{\max}} > \rbrac{\mu - 2\sqrt{\varepsilon_2} L_{\max}}^2,
\end{align*}
since when $\sqrt{\varepsilon_2} < \frac{\mu}{2L_{\max}}$, we have $\mu - 2\sqrt{\varepsilon_2}L_{\max} > 0$.
We expand the right-hand side and obtain:
\begin{align*}
   \rbrac{\mu - 2\sqrt{\varepsilon_2} L_{\max}}^2 &= \mu^2 - 4\sqrt{\varepsilon_2} L_{\max} + 4\varepsilon_2 L_{\max}^2 < 2\mu^2 - 4\sqrt{\varepsilon_2}L_{\max} < 2\mu^2.
\end{align*}
For the left-hand side, as we have already shown in \ref{eq:mu-est}, we have 
\begin{align*}
   4L_{\gamma}\rbrac{1 + \gamma L_{\max}}\rbrac{\mu + 4\sqrt{\varepsilon_2}L_{\max} + 4\varepsilon_2 L_{\max}} &\geq 4\mu\rbrac{\mu + 4\sqrt{\varepsilon_2}L_{\max} + 2\varepsilon_2 L_{\max}} > 4\mu^2.
\end{align*}
Combining the above inequality we arrive at $B_2 > B_1$.

\paragraph{The convergence guarantee}:
Given that we select $\alpha$ properly, we have 
\begin{align*}
   \cE_K \leq \rbrac{1 - \alpha\cdot\frac{\gamma\rbrac{\mu - 2\sqrt{\varepsilon_2} L_{\max}}}{4\rbrac{1 + \gamma L_{\max}}}}^K \cE_0,
\end{align*}
where $\cE_k = \gamma M^{\gamma}\rbrac{x_k} - \gamma M^{\gamma}_{\inf}$.
We do not have expectation here since we are in the full batch case.
Specifically, if we choose the largest $\alpha$ possible, we have 
\begin{equation*}
   \cE_K \leq \rbrac{1 - \frac{\mu}{4L_{\gamma}\rbrac{1 + \gamma L_{\max}}} \cdot S\rbrac{\varepsilon_2}}^k\cE_0,
\end{equation*}
where 
\begin{equation*}
   S(\varepsilon_2) = \frac{\rbrac{\mu - 2\sqrt{\varepsilon_2}L_{\max}}\rbrac{1 - 2\sqrt{\varepsilon_2}\frac{L_{\max}}{\mu}}}{\mu + 4\sqrt{\varepsilon_2}L_{\max} + 4\varepsilon_2 L_{\max}},
\end{equation*}
satisfies $0 < S(\varepsilon_2) \leq 1$ is the factor of slowing down due to inexact proximity operator evaluation.
Using smoothness of $\gamma L_{\gamma}$, if we denote $\Delta_k = \norm{x_k - x_\star}^2$ where $x_\star$ is a minimizer of both $M^{\gamma} $ and $f$ since we assume we are in the interpolation regime (\Cref{asp:int-pl-rgm}), we have 
\begin{align*}
   \cE_0 \leq \frac{\gamma L_{\gamma}}{2}\Delta_0.
\end{align*} 
Using star strong convexity (quadratic growth property), we have 
\begin{align*}
   \cE_K \geq \frac{\gamma\mu}{2\rbrac{1 + \gamma L_{\max}}} \Delta_K.
\end{align*}
As a result, we can transform the above convergence guarantee into 
\begin{align*}
   \Delta_K \leq \rbrac{1 - \frac{\mu}{4L_{\gamma}\rbrac{1 + \gamma L_{\max}}} \cdot S\rbrac{\varepsilon_2}}^K\cdot\frac{L{\gamma}\rbrac{1 + \gamma L_{\max}}}{\mu}\Delta_0.
\end{align*}
This completes the proof.

\subsection{\texorpdfstring{Proof of \Cref{thm:010-stncvx-rela}}{Proof of Theorem~\ref{thm:010-stncvx-rela}}}

We start with formalizing the problem.
Using (\ref{eq:trans-update}) and (\ref{eq:biased-estimator}), we can write the update rule of \Cref{alg:inexact-fedexprox-full} as
\begin{equation}
   \label{eq:thm333}
   x_{k+1} = x_k - \alpha \cdot \rbrac{\frac{1}{n}\sum_{i=1}^{n} \gamma\nabla\MoreauSub{\gamma}{f_i}{x_k} - \frac{1}{n}\sum_{i=1}^{n}\rbrac{\tx_{i, k+1} - \ProxSub{\gamma f_i}{x_k}}}.
\end{equation}
Since by \Cref{def:inexact-2}, we have $\norm{\tx_{i, k+1} - \ProxSub{\gamma f_i}{x_k}}^2 \leq \varepsilon_2\norm{\gamma \nabla \MoreauSub{\gamma}{f_i}{x_k}}^2$, we can view the left hand side as a compressed version of the true gradient.
Specifically, there are two possible perspectives:
\begin{enumerate}
   \item[(I).] Let $\cC_i\rbrac{\cdot}$ be the compressing mapping with the $i$-th client, $i \in \cbrac{1, 2, \hdots, n}$, defined as 
   \begin{align*}
      \cC_i\rbrac{\gamma \nabla \MoreauSub{\gamma}{f_i}{x_k}} \eqdef \gamma \nabla \MoreauSub{\gamma}{f_i}{x_k} - \rbrac{\tx_{i, k+1} - \ProxSub{\gamma f_i}{x_k}}.
   \end{align*}
   In this way, we reformulate (\ref{eq:thm333}) as 
   \begin{align}
      \label{eq:thm334}
      x_{k+1} = x_k - \alpha\cdot\frac{1}{n}\sum_{i=1}^{n}\cC_i\rbrac{\gamma \nabla \MoreauSub{\gamma}{f_i}{x_k}}.
   \end{align}
   (\ref{eq:thm334}) is exactly {\DCGD} with biased compression.
   We can easily prove that
   \begin{align*}
      \cC_i &\in \bB^1\rbrac{\alpha = 1-2\sqrt{\varepsilon_2}, \beta = \frac{1 - \sqrt{\varepsilon_2}}{1 + \varepsilon_2}} \\
      \cC_i &\in \bB^2\rbrac{\xi = 1 - \sqrt{\varepsilon_2}, \beta = \frac{1 - \sqrt{\varepsilon_2}}{1 + \varepsilon_2}} \\
      \cC_i &\in \bB^3\rbrac{\delta = \frac{1}{1 - \varepsilon_2}}.
   \end{align*}
   However, {\DCGD} with biased compression may fail to converge even if the above formulation of compression mapping seems quite nice.
   For an example of such failure, we refer the readers to \citet[Example 1]{beznosikov2023biased}.
   This limitation can be circumvented by employing an error feedback mechanism; however, this approach requires modifications to the original algorithm.
   We therefore leave it as a future research direction.

   \item[(II).] We can also view it as if we are in the single node case.
   Let $\cC\rbrac{\cdot}$ be the compressing mapping defined as 
   \begin{align}
      \label{eq:thm335}
      \cC\rbrac{\nabla M^{\gamma}\rbrac{x_k}} &\eqdef \frac{1}{n}\sum_{i=1}^{n} \gamma\nabla\MoreauSub{\gamma}{f_i}{x_k} - \frac{1}{n}\sum_{i=1}^{n}\rbrac{\tx_{i, k+1} - \ProxSub{\gamma f_i}{x_k}} \notag \\
      &= \gamma \nabla M^{\gamma}\rbrac{x_k} - \frac{1}{n}\sum_{i=1}^{n}\rbrac{\tx_{i, k+1} - \ProxSub{\gamma f_i}{x_k}} .
   \end{align}
   This formulation leads us to the convergence guarantee appeared in \Cref{thm:010-stncvx-rela}, as we illustrate below.
\end{enumerate}
Let us first analyze $\cC$ defined in (\ref{eq:thm335}).
We will verify it belongs to $\bB^3\rbrac{\delta}$.
The inequality we want to prove can be written equivalently as 
\begin{align}
   \label{eq:tiewwtp}
   \norm{\gamma\nabla M^{\gamma}\rbrac{x_k} - \frac{1}{n}\sum_{i=1}^n\rbrac{\tx_{i, k+1} - \ProxSub{\gamma f_i}{x_k}} - \gamma\nabla M^{\gamma}\rbrac{x_k}}^2 \leq \rbrac{1 - \frac{1}{\delta}}\norm{\gamma\nabla M^{\gamma}\rbrac{x_k}}^2,
\end{align}
which is exactly 
\begin{align*}
   \norm{\frac{1}{n}\sum_{i=1}^n\rbrac{\tx_{i, k+1} - \ProxSub{\gamma f_i}{x_k}}}^2 \leq \norm{\gamma\nabla M^{\gamma}\rbrac{x_k}}^2
\end{align*}
For the left-hand side, using the convexity of $\norm{\cdot}^2$ in combination with \Cref{def:inexact-2}, we obtain 
\begin{align*}
   \norm{\frac{1}{n}\sum_{i=1}^n\rbrac{\tx_{i, k+1} - \ProxSub{\gamma f_i}{x_k}}}^2 &\leq \frac{1}{n}\sum_{i=1}^n\norm{\tx_{i, k+1} - \ProxSub{\gamma f_i}{x_k}}^2 \\
   &\leq \frac{\varepsilon_2}{n}\sum_{i=1}^n \norm{x_k - \ProxSub{\gamma f_i}{x_k}}^2.
\end{align*} 
Let $x_\star$ be a minimizer of $f$, since we assume \Cref{asp:int-pl-rgm} holds, by \Cref{fact:5:global:equivalence}, it is also a minimizer of $\gamma M^{\gamma}$,
\begin{eqnarray*}
   \frac{\varepsilon_2}{n}\sum_{i=1}^n \norm{x_k - \ProxSub{\gamma f_i}{x_k}}^2 &\overset{(\ref{eq:ppt-2-gradient})}{=}& \frac{\varepsilon_2}{n}\sum_{i=1}^n \norm{\gamma \nabla \MoreauSub{\gamma}{f_i}{x_k}}^2 \\
   &=& \frac{\varepsilon_2}{n}\sum_{i=1}^n \norm{\gamma \nabla \MoreauSub{\gamma}{f_i}{x_k} - \gamma \nabla \MoreauSub{\gamma}{f_i}{x_\star}}^2 \\
   &\overset{\text{\Cref{tech-lemma:smooth}}}{\leq}& \frac{2\varepsilon_2}{n}\sum_{i=1}^n \frac{\gamma L_i}{1 + \gamma L_i}\rbrac{\gamma\MoreauSub{\gamma}{f_i}{x_k} - \gamma\MoreauSub{\gamma}{f_i}{x_\star}} \\
   &\leq& \frac{2\varepsilon_2\gamma L_{\max}}{1 + \gamma L_{\max}}\rbrac{\gamma M^{\gamma}\rbrac{x_k} - \gamma M^{\gamma}\rbrac{x_\star}}.
\end{eqnarray*}
We then notice that as it is illustrated by \Cref{lemma:PL-lemma}, we have 
\begin{align*}
   \rbrac{1 - \frac{1}{\delta}}\norm{\gamma \nabla M^{\gamma}\rbrac{x_k}}^2 \geq \rbrac{1 - \frac{1}{\delta}}\frac{\gamma\mu}{2\rbrac{1 + \gamma L_{\max}}}\rbrac{\gamma M^{\gamma}\rbrac{x_k} - \gamma M^{\gamma}\rbrac{x_\star}}.
\end{align*}
Combining the above two inequalities, we know that the following inequality is a sufficient condition for (\ref{eq:tiewwtp}),
\begin{align*}
   \frac{2\varepsilon_2\gamma L_{\max}}{1 + \gamma L_{\max}}\rbrac{\gamma M^{\gamma}\rbrac{x_k} - \gamma M^{\gamma}\rbrac{x_\star}} \leq \rbrac{1 - \frac{1}{\delta}}\frac{\gamma\mu}{2\rbrac{1 + \gamma L_{\max}}}\rbrac{\gamma M^{\gamma}\rbrac{x_k} - \gamma M^{\gamma}\rbrac{x_\star}}.
\end{align*}
It is easy to check that if we pick 
\begin{equation}
   \label{eq:check-verify}
   \delta = \frac{\mu}{\mu - 4\varepsilon_2L_{\max}} > 0,
\end{equation}
the condition is met.
However, for this to hold, we must ensure that $\varepsilon_2 < \frac{\mu}{4L_{\max}}$.

As we mentioned in \Cref{sec:bcomp}, \citet{beznosikov2023biased} provided the theory of {\CGD} with biased compressor belongs to $\bB^3\rbrac{\delta}$.
We have already shown that $\cC \in \bB^3\rbrac{\delta = \frac{\mu}{\mu - 4\varepsilon_2 L_{\max}}}$, when $\varepsilon_2 < \frac{4L_{\max}}{\mu}$.
Notice that our objective $\gamma M^{\gamma}$ is $\gamma L_{\gamma}$-smooth and $\frac{\gamma\mu}{1 + \gamma L_{\max}}$-PL.\footnote{\Cref{thm:biased-compression} remains valid if we replace $f$ being strongly convex with PL.}
Therefore, as long as $0 < \alpha \leq \frac{1}{\gamma L_{\gamma}}$ and $\varepsilon_2 < \frac{\mu}{4L_{\max}}$, we have 
\begin{align*}
   \cE_K \leq \rbrac{1 - \frac{\mu - 4\varepsilon_2L_{\max}}{\mu}\cdot \frac{\gamma\mu}{4\rbrac{1 + \gamma L_{\max}}}\cdot \alpha}^K\cE_0,
\end{align*}
Taking $\alpha = \frac{1}{\gamma L_{\gamma}}$, which is the largest step size possible, we can further simplify the above convergence into 
\begin{align*}
   M^{\gamma}\rbrac{x_k} - M^{\gamma}_\star \leq \rbrac{1 - \rbrac{1 - \frac{4\varepsilon_2 L_{\max}}{\mu}}\cdot \frac{\mu}{4L_{\gamma}\rbrac{1 + \gamma L_{\max}}}}^K\rbrac{M^{\gamma}\rbrac{x_0} - M^{\gamma_\star}}.
\end{align*}
Using smoothness of $\gamma L_{\gamma}$, if we denote $\Delta_k = \norm{x_k - x_\star}^2$ where $x_\star$ is a minimizer of both $M^{\gamma} $ and $f$ since we assume we are in the interpolation regime (\Cref{asp:int-pl-rgm}), we have 
\begin{align*}
   \cE_0 \leq \frac{\gamma L_{\gamma}}{2}\Delta_0.
\end{align*} 
Using star strong convexity (quadratic growth property), we have 
\begin{align*}
   \cE_K \geq \frac{\gamma\mu}{2\rbrac{1 + \gamma L_{\max}}} \Delta_K.
\end{align*}
As a result, we can transform the above convergence guarantee into 
\begin{align*}
   \Delta_K \leq \rbrac{1 - \rbrac{1 - \frac{4\varepsilon_2 L_{\max}}{\mu}}\cdot \frac{\mu}{4L_{\gamma}\rbrac{1 + \gamma L_{\max}}}}^K \cdot \frac{L_{\gamma}\rbrac{1 + \gamma L_{\max}}}{\mu}\Delta_0.
\end{align*}
This completes the proof.

\subsection{\texorpdfstring{Proof of \Cref{thm:complexity}}{Proof of Theorem~\ref{thm:complexity}}}

Notice that we assume each $f_i$ is $L_i$-smooth and convex.
The local optimization of each client can be written as 
\begin{equation*}
   \min_{z \in \R^d} \cbrac{\aloc{\gamma}{k, i}{z} = f_i\rbrac{z} + \frac{1}{2\gamma}\norm{z - x_k}^2}, 
\end{equation*}
It is easy to see that $\aloc{\gamma}{k, i}{z}$ is $L_i + \frac{1}{\gamma}$-smooth and $\frac{1}{\gamma}$-strongly convex.
We first provide the convergence theory of {\GD} for reference.

\paragraph{Theory of {\GD}:}
For a $\widehat{\mu}$-strongly convex, $\widehat{L}$-smooth function $\phi$, the algorithm can be formulated as
\begin{equation}
   \label{eq:GD}
   z_{t+1} = z_t - \eta \nabla \phi(z_t), \tag{GD}
\end{equation}
where $z_t$ is the iterate in the $t$-th iteration, and $\eta > 0$ is the step size.
{\GD} with step size $\eta \in (0, \frac{1}{\widehat{L}}]$ generates iterates that satisfy 
\begin{equation*}
   \norm{z_t - z_\star}^2 \leq \rbrac{1 - \eta\widehat{\mu}}^t\norm{z_0 - z_\star}^2,
\end{equation*}
where $z_\star$ is a minimizer of $\phi$, $t$ is the number of iterations (number of gradient evaluations).

\paragraph{Approximation satisfying \Cref{def:inexact-1}:}
Notice that $\ProxSub{\gamma f_i}{x_k}$ is the minimizer of $\aloc{\gamma}{k, i}{z}$ and $z_0 = x_k$.
As a result, if we run {\GD} with the largest step size $\frac{\gamma}{1 + \gamma L_i}$, 
\begin{align}
   \label{eq:conv-gd-001}
   \norm{z_t - \ProxSub{\gamma f_i}{x_k}}^2 \leq \rbrac{1 - \frac{1}{1 + \gamma L_i}}^t\norm{x_k - \ProxSub{\gamma f_i}{x_k}}^2
\end{align}
We have 
\begin{equation*}
   t = \cO\rbrac{\rbrac{1 + \gamma L_i}\log\rbrac{\frac{\norm{x_k - \ProxSub{\gamma f_i}{x_k}}^2}{\varepsilon_1}}}.
\end{equation*}
The unknown term $\norm{x_k - \ProxSub{\gamma f_i}{x_k}}^2$ within the log can be bounded by 
\begin{align}
   \label{eq:upperbound-gd-comp}
   \norm{x_k - \ProxSub{\gamma f_i}{x_k}}^2 &= \norm{z_0 - z_\star}^2 \notag \\
   &\leq \gamma^2 \norm{\nabla \aloc{\gamma}{k, i}{z_0} - \nabla \aloc{\gamma}{k, i}{z_\star}}^2 = \norm{\gamma\nabla f_i \rbrac{x_k}}^2,
\end{align}
which can be easily calculated.

\paragraph{Approximation satisfying \Cref{def:inexact-2}:}
According to (\ref{eq:conv-gd-001}), we have 
\begin{equation*}
   t = \cO\rbrac{\rbrac{1 + \gamma L_i}\log\rbrac{\frac{1}{\varepsilon_2}}}.
\end{equation*}
This completes the proof.

\subsection{\texorpdfstring{Proof of \Cref{thm:complexity-accelerated}}{Proof of Theorem~\ref{thm:complexity-accelerated}}}

We first provide the theory of {\AGD} \citep{nesterov2013introductory}.

\paragraph{Theory of {\AGD}:}
For a $\widehat{\mu}$-strongly convex, $\widehat{L}$-smooth function $\phi$, the algorithm can be formulated as 
\begin{align}
   \label{eq:fomu-AGD}
   y_{t+1} &= z_t + \alpha\rbrac{z_t - z_{t-1}} \notag \\
   z_{t+1} &= y_{t+1} - \eta\nabla \phi\rbrac{y_{t+1}} \tag{AGD},
\end{align}
where $z_{t}, y_t$ are iterates, $\eta > 0$ is the step size, $\alpha > 0$ is the momentum parameter. 
{\AGD} with step size $\eta = \frac{1}{\widehat{L}}$, momentum $\alpha = \frac{\sqrt{\widehat{L}} - \sqrt{\widehat{\mu}}}{\sqrt{\widehat{L}} + \sqrt{\widehat{\mu}}}$ generates iterates that satisfy 
\begin{align*}
   \norm{z_t - z_\star}^2 \leq \frac{2\widehat{L}}{\widehat{\mu}}\cdot\rbrac{1 - \sqrt{\frac{\widehat{\mu}}{\widehat{L}}}}^t \norm{z_0 - z_\star}^2,
\end{align*}
where $z_\star$ is a minimizer of $\phi$, $t$ is the number of iterations (number of gradient evaluations).

\paragraph{Approximation satisfying \Cref{def:inexact-1}:}
Notice that $\ProxSub{\gamma f_i}{x_k}$ is the minimizer of $\aloc{\gamma}{k, i}{z}$ and $z_0 = x_k$.
As a result, if we run {\AGD} with the step size $\frac{\gamma}{1 + \gamma L_i}$ and momentum $\alpha = \frac{\sqrt{1 + \gamma L_i} - 1}{\sqrt{1 + \gamma L_i} + 1}$, 
\begin{align}
   \label{eq:conv-agd-001}
   \norm{z_t - \ProxSub{\gamma f_i}{x_k}}^2 \leq 2\cdot\rbrac{1 + \gamma L_i}\rbrac{1 - \frac{1}{\sqrt{1 + \gamma L_i}}}^t\norm{x_k - \ProxSub{\gamma f_i}{x_k}}^2.
\end{align}
We have 
\begin{equation*}
   t = \cO\rbrac{\sqrt{1 + \gamma L_i}\log\rbrac{\frac{\rbrac{1 + \gamma L_i}\cdot\norm{x_k - \ProxSub{\gamma f_i}{x_k}}^2}{\varepsilon_1}}}
\end{equation*}
Similar to the proof of \Cref{thm:complexity}, since we have according to (\ref{eq:upperbound-gd-comp}),
\begin{equation*}
   \norm{x_k - \ProxSub{\gamma f_i}{x_k}}^2 \leq \norm{\gamma\nabla f_i\rbrac{x_k}}^2,
\end{equation*}
it is straightforward to determine the number of local iterations needed.

\paragraph{Approximation satisfying \Cref{def:inexact-2}:}
Using (\ref{eq:conv-agd-001}), we have 
\begin{align*}
   t = \cO\rbrac{\sqrt{1 + \gamma L_i}\log\rbrac{\frac{1 + \gamma L_i}{\varepsilon_2}}}.
\end{align*}

\subsection{\texorpdfstring{Proof of \Cref{thm:020-mini}}{Proof of Theorem~\ref{thm:020-mini}}}

In this case, the gradient estimator is defined as 
\begin{equation}
   \label{eq:gxk-mini}
   g(x_k) = \frac{1}{\tau}\sum_{i \in S_k}\rbrac{\gamma\nabla\MoreauSub{\gamma}{f_i}{x_k} - \rbrac{\tx_{i, k+1} - \ProxSub{\gamma f_i}{x_k}}}.
\end{equation}
Notice that we have 
\begin{align*}
   & \inner{\gamma \nabla M^{\gamma}\rbrac{x_k}}{\Exp{g(x_k)}} \\
   &= \inner{\gamma \nabla M^{\gamma}\rbrac{x_k}}{\Exp{\frac{1}{\tau}\sum_{i \in S_k}\gamma\nabla\MoreauSub{\gamma}{f_i}{x_k} - \frac{1}{\tau}\sum_{i\in S_k}\rbrac{\tx_{i, k+1} - \ProxSub{\gamma f_i}{x_k}}}} \\
   &= \inner{\gamma \nabla M^{\gamma}\rbrac{x_k}}{\gamma \nabla M^{\gamma}\rbrac{x_k} - \frac{1}{n}\sum_{i=1}^{n}\rbrac{\tx_{i, k+1} - \ProxSub{\gamma f_i}{x_k}}}.
\end{align*} 
Using the same technique in the proof of \Cref{thm:020-full}, we are able to obtain that 
\begin{align*}
   \inner{\gamma \nabla M^{\gamma}\rbrac{x_k}}{\Exp{g(x_k)}} \geq \rbrac{1 - \frac{2\sqrt{\varepsilon_2}L_{\max}}{\mu}}\cdot\norm{\gamma \nabla M^{\gamma}\rbrac{x_k}}^2.
\end{align*}
Thus, as long as we pick $\varepsilon_2 < \frac{\mu^2}{4L_{\max}^2}$, we can pick $b = 1 - \sqrt{\varepsilon_2}\cdot\frac{2L_{\max}}{\mu}$ and $c = 0$.
We then compute $\Exp{\norm{g(x_k)}^2}$,
\begin{align*}
   \Exp{\norm{g(x_k)}^2} &= \Exp{\norm{\frac{1}{\tau}\sum_{i \in S_k} \gamma \nabla \MoreauSub{\gamma}{f_i}{x_k} - \frac{1}{\tau}\sum_{i \in S_k}\rbrac{\tx_{i, k+1} - \ProxSub{\gamma f_i}{x_k}}}^2} \\
   &= \underbrace{\Exp{\norm{\frac{1}{\tau}\sum_{i \in S_k} \gamma \nabla \MoreauSub{\gamma}{f_i}{x_k}}^2}}_{\eqdef T_1} + \underbrace{\Exp{\norm{\frac{1}{\tau}\sum_{i \in S_k}\rbrac{\tx_{i, k+1} - \ProxSub{\gamma f_i}{x_k}}}^2}}_{\eqdef T_2} \\
   &\quad  \underbrace{- 2\Exp{\inner{\frac{1}{\tau}\sum_{i \in S_k} \gamma \nabla \MoreauSub{\gamma}{f_i}{x_k}}{\frac{1}{\tau}\sum_{i \in S_k}\rbrac{\tx_{i, k+1} - \ProxSub{\gamma f_i}{x_k}}}}}_{\eqdef T_3}.
\end{align*}
We try to provide upper bounds for those terms separately.

\paragraph{Term $T_1$:}
We have 
\begin{align*}
   T_1 &= \frac{n - \tau}{\tau\rbrac{n - 1}}\cdot \frac{1}{n}\sum_{i=1}^{n}\norm{\gamma \nabla \MoreauSub{\gamma}{f_i}{x_k}}^2 + \frac{n\rbrac{\tau - 1}}{\tau\rbrac{n - 1}}\cdot\norm{\gamma\nabla M^{\gamma}\rbrac{x_k}}^2.
\end{align*}
Using smoothness of $\gamma \Moreau^{\gamma}_{f_i}$ and the fact that we are in the interpolation regime, we have 
\begin{align}
   \label{eq:sss-t1}
   T_1 &= \frac{n - \tau}{\tau\rbrac{n - 1}}\cdot \frac{1}{n}\sum_{i=1}^n\norm{\gamma \nabla \MoreauSub{\gamma}{f_i}{x_k} - \gamma \nabla \MoreauSub{\gamma}{f_i}{x_\star}}^2 + \frac{n\rbrac{\tau - 1}}{\tau\rbrac{n - 1}}\cdot\norm{\gamma\nabla M^{\gamma}\rbrac{x_k}}^2 \notag \\
   &\leq \frac{n - \tau}{\tau\rbrac{n - 1}}\cdot \frac{1}{n}\sum_{i=1}^n\frac{2\gamma L_{i}}{1 + \gamma L_{i}} \cdot \rbrac{\gamma M^{\gamma}_{f_i}\rbrac{x_k} - \gamma \rbrac{M^{\gamma}_{f_i}}_{\inf}} + \frac{n\rbrac{\tau - 1}}{\tau\rbrac{n - 1}}\cdot\norm{\gamma\nabla M^{\gamma}\rbrac{x_k}}^2 \notag \\
   &\leq \frac{n - \tau}{\tau\rbrac{n - 1}}\cdot \frac{2\gamma L_{\max}}{1 + \gamma L_{\max}} \cdot \rbrac{\gamma M^{\gamma}\rbrac{x_k} - \gamma M^{\gamma}_{\inf}} + \frac{n\rbrac{\tau - 1}}{\tau\rbrac{n - 1}}\cdot\norm{\gamma\nabla M^{\gamma}\rbrac{x_k}}^2.
\end{align}

\paragraph{Term $T_2$:}
It is easy to see that using convexity of the squared Euclidean norm, we have 
\begin{align*}
   T_2 &\leq \Exp{\frac{1}{\tau}\sum_{i\in S_k}\norm{\tx_{i, k+1} - \ProxSub{\gamma f_i}{x_k}}^2} \\
   &= \frac{1}{n}\sum_{i=1}^{n} \norm{\tx_{i, k+1} - \ProxSub{\gamma f_i}{x_k}}^2 \overset{(\ref{eq:index2-approx})}{\leq} \frac{\varepsilon_2}{n}\sum_{i=1}^n \norm{\gamma \nabla \MoreauSub{\gamma}{f_i}{x_k}}^2.
\end{align*}
Using smoothness of each individual $\gamma \MoreauSub{\gamma}{f_i}{x_k}$ and the fact we are in the interpolation regime, we have 
\begin{align}
   \label{eq:sss-t2}
   T_2 & \leq \frac{2\varepsilon_2\gamma L_{\max}}{1 + \gamma L_{\max}}\rbrac{\gamma M^{\gamma}\rbrac{x_k} - \gamma M^{\gamma}_{\inf}}.
\end{align}
\paragraph{Term $T_3$:}
We have 
\begin{align*}
   T_3 &= -2 \cdot \frac{n-\tau}{\tau\rbrac{n - 1}}\cdot\frac{1}{n}\sum_{i=1}^{n}\inner{\gamma\nabla\MoreauSub{\gamma}{f_i}{x_k}}{\tx_{i, k+1} - \ProxSub{\gamma f_i}{x_k}} \\
    &\quad -2 \cdot \frac{n\rbrac{\tau - 1}}{\tau\rbrac{n - 1}}\cdot\inner{\gamma \nabla M^{\gamma}\rbrac{x_k}}{\frac{1}{n}\sum_{i=1}^{n}\rbrac{\tx_{i, k+1} - \ProxSub{\gamma f_i}{x_k}}}. \\
\end{align*}
Using Cauchy-Schwarz inequality and convexity, we further obtain 
\begin{align*}
   T_3 &\leq 2 \cdot \frac{n-\tau}{\tau\rbrac{n - 1}} \cdot\frac{1}{n}\sum_{i=1}^{n}\norm{\gamma\nabla\MoreauSub{\gamma}{f_i}{x_k}}\norm{\tx_{i, k+1} - \ProxSub{\gamma f_i}{x_k}} \\
   &\quad + 2\cdot\frac{n\rbrac{\tau - 1}}{\tau\rbrac{n - 1}} \norm{\gamma \nabla M^{\gamma}\rbrac{x_k}} \cdot \frac{1}{n}\sum_{i=1}^{n}\norm{\tx_{i, k+1} - \ProxSub{\gamma f_i}{x_k}}.
\end{align*}
Using similar approaches in the previous paragraphs, we have 
\begin{align}
   \label{eq:sss-t3}
   &T_3 \notag \\
   &\overset{(\ref{eq:index2-approx})}{\leq} \frac{2\rbrac{n-\tau}}{\tau\rbrac{n - 1}} \cdot \frac{\sqrt{\varepsilon_2}}{n}\sum_{i=1}^{n}\norm{\gamma \nabla \MoreauSub{\gamma}{f_i}{x_k}}^2 + \frac{2n\rbrac{\tau - 1}}{\tau\rbrac{n - 1}}\norm{\gamma M^{\gamma}\rbrac{x_k}} \frac{\sqrt{\varepsilon_2}}{n}\cdot\sum_{i=1}^{n}\norm{\gamma \nabla \MoreauSub{\gamma}{f_i}{x_k}} \notag \\
   &\leq \frac{2\rbrac{n-\tau}}{\tau\rbrac{n - 1}} \cdot \frac{\sqrt{\varepsilon_2}}{n}\sum_{i=1}^{n}\norm{\gamma \nabla \MoreauSub{\gamma}{f_i}{x_k} - \gamma \nabla \MoreauSub{\gamma}{f_i}{x_\star}}^2 \notag \\
   &\quad + \frac{2n\rbrac{\tau - 1}}{\tau\rbrac{n - 1}}\norm{\gamma M^{\gamma}\rbrac{x_k}} \frac{\sqrt{\varepsilon_2}}{n}\cdot\sum_{i=1}^{n}\norm{\gamma \nabla \MoreauSub{\gamma}{f_i}{x_k} - \gamma \nabla \MoreauSub{\gamma}{f_i}{x_k}} \notag \\
   &\leq \frac{4\sqrt{\varepsilon_2}\rbrac{n-\tau}}{\tau\rbrac{n-1}}\cdot \frac{\gamma L_{\max}}{1 + \gamma L_{\max}}\rbrac{\gamma M^{\gamma}\rbrac{x_k} - \gamma M^{\gamma}_{\inf}} \notag \\
   &\quad + \frac{4\sqrt{\varepsilon_2}n\rbrac{\tau - 1}}{\tau\rbrac{n - 1}}\cdot \frac{\gamma L_{\max}}{1 + \gamma L_{\max}} \norm{x_k - x_\star}\norm{\gamma \nabla M^{\gamma}\rbrac{x_k}} \notag \\
   &\overset{(\ref{eq:ineq-pr1-4})}{\leq} \frac{4\sqrt{\varepsilon_2}\rbrac{n-\tau}}{\tau\rbrac{n-1}}\cdot \frac{\gamma L_{\max}}{1 + \gamma L_{\max}}\rbrac{\gamma M^{\gamma}\rbrac{x_k} - \gamma M^{\gamma}_{\inf}} \notag \\
   &\quad + \frac{4\sqrt{\varepsilon_2}n\rbrac{\tau - 1}}{\tau\rbrac{n - 1}}\cdot \frac{L_{\max}}{\mu} \norm{\gamma \nabla M^{\gamma}\rbrac{x_k}}^2.
\end{align}
Combining (\ref{eq:sss-t1}), (\ref{eq:sss-t2}) and (\ref{eq:sss-t3}), we have 
\begin{align}
   \label{eq:sss-result}
   \sum_{i=1}^3 T_i &\leq 2\rbrac{\varepsilon_2 + \frac{2\sqrt{\varepsilon_2}\rbrac{n - \tau}}{\tau\rbrac{n - 1}} + \frac{\rbrac{n - \tau}}{\tau\rbrac{n - 1}}}\cdot\frac{\gamma L_{\max}}{1 + \gamma L_{\max}}\cdot\rbrac{\gamma M^{\gamma}\rbrac{x_k} - \gamma M^{\gamma}_{\inf}} \notag \\
   &\quad + \rbrac{\frac{n\rbrac{\tau - 1}}{\tau\rbrac{n - 1}} + \frac{4\sqrt{\varepsilon_2}n\rbrac{\tau - 1}}{\tau\rbrac{n - 1}}}\cdot\frac{L_{\max}}{\mu}\cdot \norm{\gamma M^{\gamma}\rbrac{x_k}}^2.
\end{align}
Therefore, it is easy to see that we can pick 
\begin{align*}
   A &= \rbrac{\varepsilon_2 + \frac{2\sqrt{\varepsilon_2}\rbrac{n - \tau}}{\tau\rbrac{n - 1}} + \frac{\rbrac{n - \tau}}{\tau\rbrac{n - 1}}}\cdot\frac{\gamma L_{\max}}{1 + \gamma L_{\max}} \\
   B &= \rbrac{\frac{n\rbrac{\tau - 1}}{\tau\rbrac{n - 1}} + \frac{4\sqrt{\varepsilon_2}n\rbrac{\tau - 1}}{\tau\rbrac{n - 1}}}\cdot\frac{L_{\max}}{\mu}, \qquad C = 0.
\end{align*}Applying Theorem 4 of \cite{demidovich2024guide}, we list the corresponding values of $A, B, C, b, c \geq 0$ below,
\begin{align*}
   A &= \frac{\gamma L_{\max}}{1 + \gamma L_{\max}}\rbrac{\varepsilon_2 + \frac{2\sqrt{\varepsilon_2}\rbrac{n-\tau}}{\tau\rbrac{n - 1}} + \frac{\rbrac{n - \tau}}{\tau\rbrac{n - 1}}} \\
   B &= \frac{n\rbrac{\tau - 1}}{\tau\rbrac{n - 1}}\rbrac{1 + \frac{4\sqrt{\varepsilon_2}L_{\max}}{\mu}}, \quad C = 0 \\
   b &= \frac{\mu - 2\sqrt{\varepsilon_2}L_{\max}}{\mu}, \quad c = 0.
\end{align*}
We know that the PL constant of $\gamma M^{\gamma}$ is given by $\frac{\gamma\mu}{4\rbrac{1 + \gamma L_{\max}}}$ and the corresponding smoothness constant is $\gamma L_{\gamma}$.
As a result, when $\alpha > 0$ satisfies 
\begin{align*}
   \alpha < \underbrace{\frac{1}{\gamma L_{\gamma}} \cdot \frac{\mu - 2\sqrt{\varepsilon_2}L_{\max}}{\mu + 4\varepsilon_2 L_{\max} + 4\sqrt{\varepsilon_2}L_{\max} + \frac{n-\tau}{\tau\rbrac{n -1}}\cdot\rbrac{4L_{\max} + 4\sqrt{\varepsilon_2}L_{\max} - \mu}}}_{\eqdef B_1^\prime},
\end{align*}
and 
\begin{align*}
   \alpha < \underbrace{\frac{4\rbrac{1 + \gamma L_{\max}}}{\gamma\rbrac{\mu - 2\sqrt{\varepsilon_2} L_{\max}}}}_{= B_2},
\end{align*}
we can obtain a convergence guarantee for the algorithm.
Notice that $B_1^{\prime} \leq B_1 < B_2$\footnote{The definition of $B_1$ is given in (\ref{eq:defB-1})}, thus we can further simplify the range of $\alpha$ to 
\begin{align*}
   \alpha \leq \underbrace{\frac{1}{\gamma L_{\gamma}} \cdot \frac{\mu - 2\sqrt{\varepsilon_2}L_{\max}}{\mu + 4\varepsilon_2 L_{\max} + 4\sqrt{\varepsilon_2}L_{\max} + \frac{n-\tau}{\tau\rbrac{n -1}}\cdot\rbrac{4L_{\max} + 4\sqrt{\varepsilon_2}L_{\max} - \mu}}}_{\eqdef B_1^\prime}.
\end{align*}
Given that we select $\alpha$ properly, we have 
\begin{align*}
   \Exp{\cE_K} \leq \rbrac{1 - \alpha\cdot\frac{\gamma\rbrac{\mu - 2\sqrt{\varepsilon_2} L_{\max}}}{4\rbrac{1 + \gamma L_{\max}}}}^K \cE_0.
\end{align*}
Specifically, if we choose the largest $\alpha$ possible, we have 
\begin{align*}
   \Exp{\cE_K} \leq \rbrac{1 - \frac{\mu}{4L_{\gamma}\rbrac{1 + \gamma L_{\max}}}\cdot S\rbrac{\varepsilon_2, \tau}}^K \cE_0,
\end{align*}
where $S\rbrac{\varepsilon_2, \tau}$ is defined as 
\begin{align*}
   S\rbrac{\varepsilon_2, \tau} = \frac{\rbrac{\mu - 2\sqrt{\varepsilon_2}L_{\max}}\rbrac{1 - 2\sqrt{\varepsilon_2}\frac{L_{\max}}{\mu}}}{\mu + 4\varepsilon_2 L_{\max} + 4\sqrt{\varepsilon_2}L_{\max} + \frac{n-\tau}{\tau\rbrac{n -1}}\cdot\rbrac{4L_{\max} + 4\sqrt{\varepsilon_2}L_{\max} - \mu}},
\end{align*}
satisfying 
\begin{align*}
   0 < S\rbrac{\varepsilon_2, \tau} \leq 1.
\end{align*}
Using smoothness of $\gamma L_{\gamma}$, if we denote $\Delta_k = \norm{x_k - x_\star}^2$ where $x_\star$ is a minimizer of both $M^{\gamma} $ and $f$ since we assume we are in the interpolation regime (\Cref{asp:int-pl-rgm}), we have 
\begin{align*}
   \cE_0 \leq \frac{\gamma L_{\gamma}}{2}\Delta_0.
\end{align*} 
Using star strong convexity (quadratic growth property), we have 
\begin{align*}
   \cE_K \geq \frac{\gamma\mu}{2\rbrac{1 + \gamma L_{\max}}} \Delta_K.
\end{align*}
As a result, we can transform the above convergence guarantee into 
\begin{align*}
   \Exp{\Delta_K} \leq \rbrac{1 - \frac{\mu}{4L_{\gamma}\rbrac{1 + \gamma L_{\max}}}\cdot S\rbrac{\varepsilon_2, \tau}}^K \cdot\frac{L{\gamma}\rbrac{1 + \gamma L_{\max}}}{\mu} \Delta_0.
\end{align*}
This completes the proof.

\section{Experiments}
\label{sec:app:exp}
We describe the settings for the numerical experiments and the corresponding results to validate our theoretical findings.
We are interested in the following optimization problem in the distributed setting,
\begin{align*}
   \min_{x \in \R^d} \cbrac{f(x) = \frac{1}{n}\sum_{i=1}^{n} f_i\rbrac{x}}.
\end{align*}
Here $n$ denotes the number of clients, $d$ is the dimension, each function $f_i: \R^d \mapsto \R$ has the following form 
\begin{equation*}
   f_i(x) = \frac{1}{2} x^{\top}\mA_ix + b_i^{\top}x + c_i,
\end{equation*}
where $\mA_i \in \bbS_{+}^d, b_i \in \R^{d}, c_i \in \R$.
Specifically, we pick $n = 20$ and $d = 300$ for the experiments.
Notice that we have 
\begin{align*}
   \nabla f_i(x) = \mA_i x - b_i; \qquad \nabla^2 f_i(x) = \mA_i \succeq \mO_d,
\end{align*}
which suggests that each $f_i$ is convex and smooth.
We can easily compute that in this case, we have 
\begin{align*}
   \ProxSub{\gamma f_i}{x} = \rbrac{\mA_i + \frac{1}{\gamma}\mI_d}^{-1}\rbrac{\frac{1}{\gamma}x - b_i}.
\end{align*}
All experiment codes were implemented in Python 3.11 using the NumPy and SciPy libraries. 
The computations were performed on a system powered by an AMD Ryzen 9 5900HX processor with Radeon Graphics, featuring 8 cores and 16 threads, running at 3.3 GHz.
Code availability: \url{https://anonymous.4open.science/r/Inexact-FedExProx-code-E783/}

\subsection{\texorpdfstring{Comparison of {\FEDPROX}, {\FEDEXPROX}, {\FEDEXPROX} with absolute approximation and relative approximation}{Comparison of FedProx, FedExProx, FedExProx with absolute approximation and relative approximation}}

In this section, we compare the convergence of {\FEDPROX}, {\FEDEXPROX} and {\FEDEXPROX} with absolute approximation and relative approximation.
For {\FEDPROX}, we simply set the server extrapolation to be $1$ while for {\FEDEXPROX}, we set its extrapolation parameter to be $\frac{1}{\gamma L_{\gamma}}$.
We assume exact proximal evaluation for the above two algorithms.
For {\FEDEXPROX} with approximations, we fix $\varepsilon_1$ and $\varepsilon_2$ to be reasonable values, respectively.
We then set their extrapolation parameter to be the optimal value under the specific setting.
Throughout the experiment, we vary the value of the local step size $\gamma$ to see its effect on all the algorithms.
Specifically, we select $\gamma$ from the set $\cbrac{1000, 100, 10, 1, 0.1. 0.01}$, and we fix $\varepsilon_1 = 0.001$, $\varepsilon_2 = 0.01$ first, then we set them to $\varepsilon_1 = 1e-6$, $\varepsilon_2 = 0.001$.

\begin{figure}
	\centering
   \subfigure{
	   \begin{minipage}[t]{0.98\textwidth}
		   \includegraphics[width=0.33\textwidth]{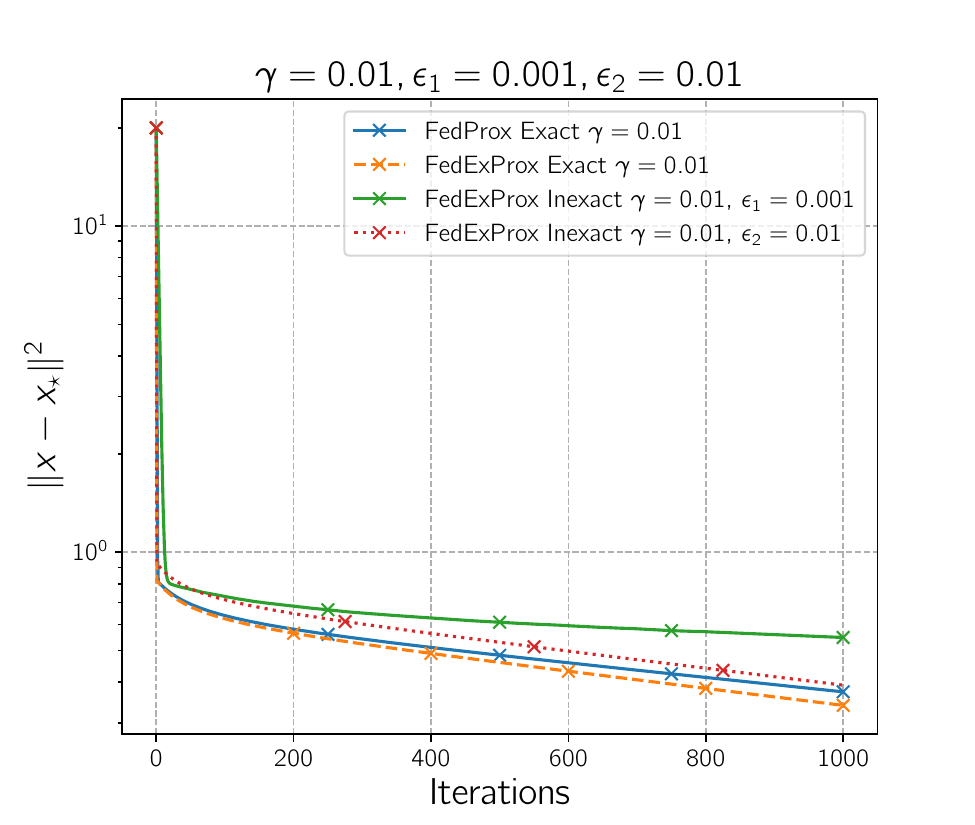} 
		   \includegraphics[width=0.33\textwidth]{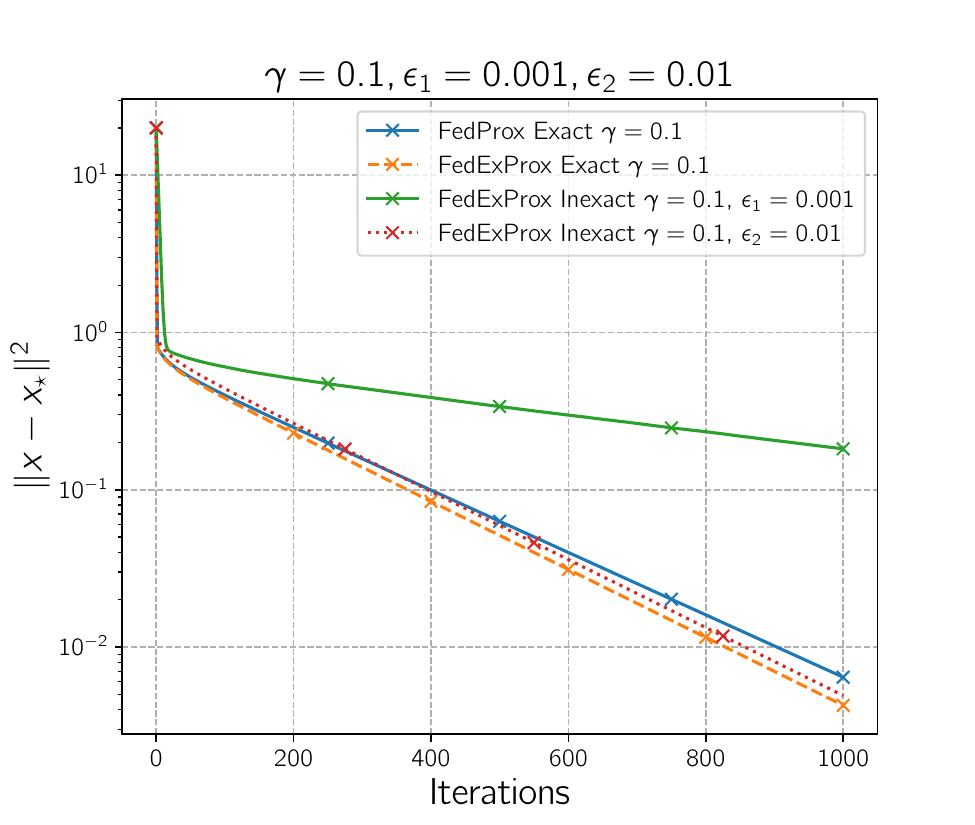}
         \includegraphics[width=0.33\textwidth]{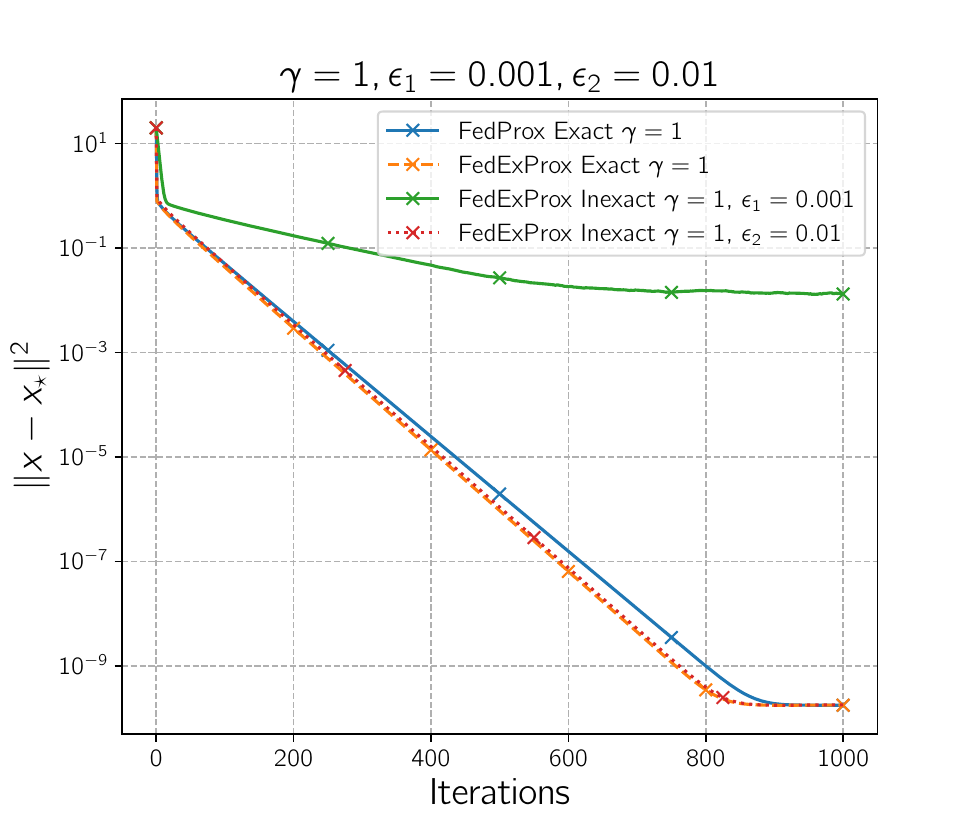}
	   \end{minipage}
   }

   \subfigure{
	   \begin{minipage}[t]{0.98\textwidth}
		   \includegraphics[width=0.33\textwidth]{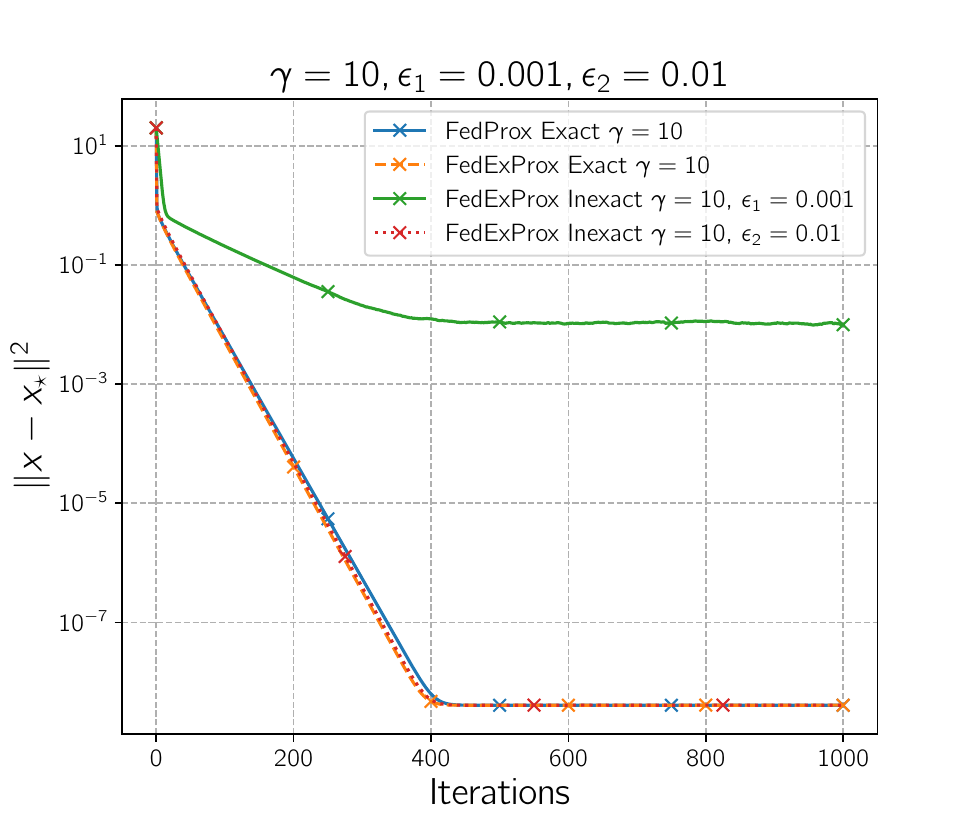} 
		   \includegraphics[width=0.33\textwidth]{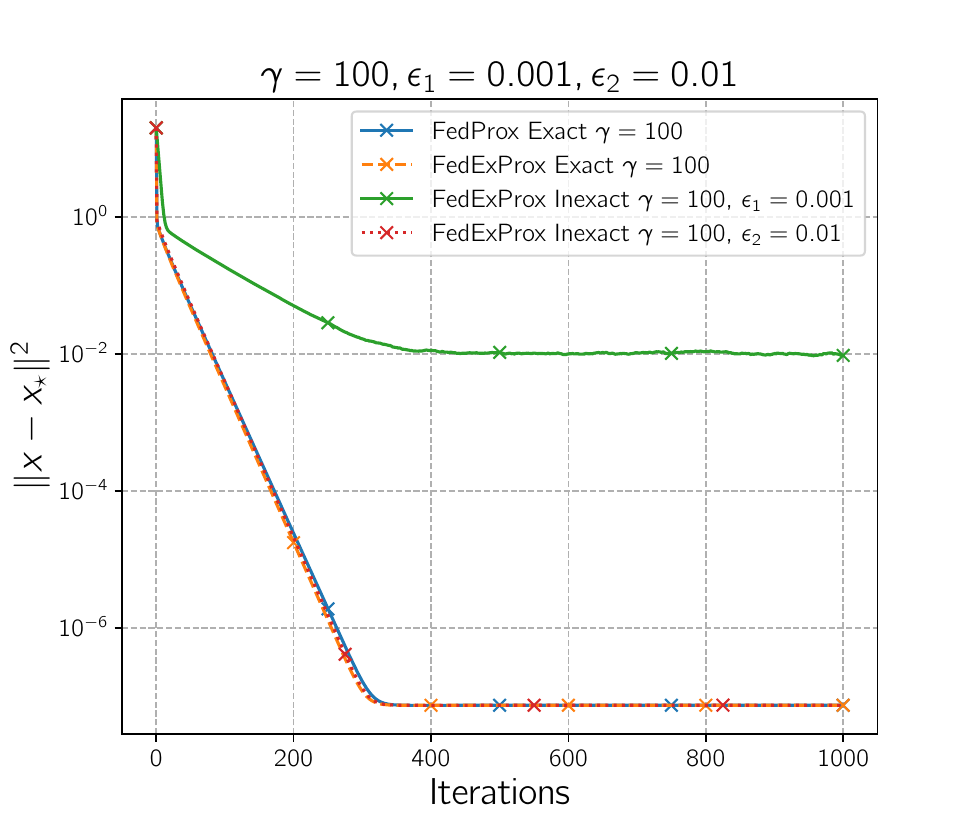}
         \includegraphics[width=0.33\textwidth]{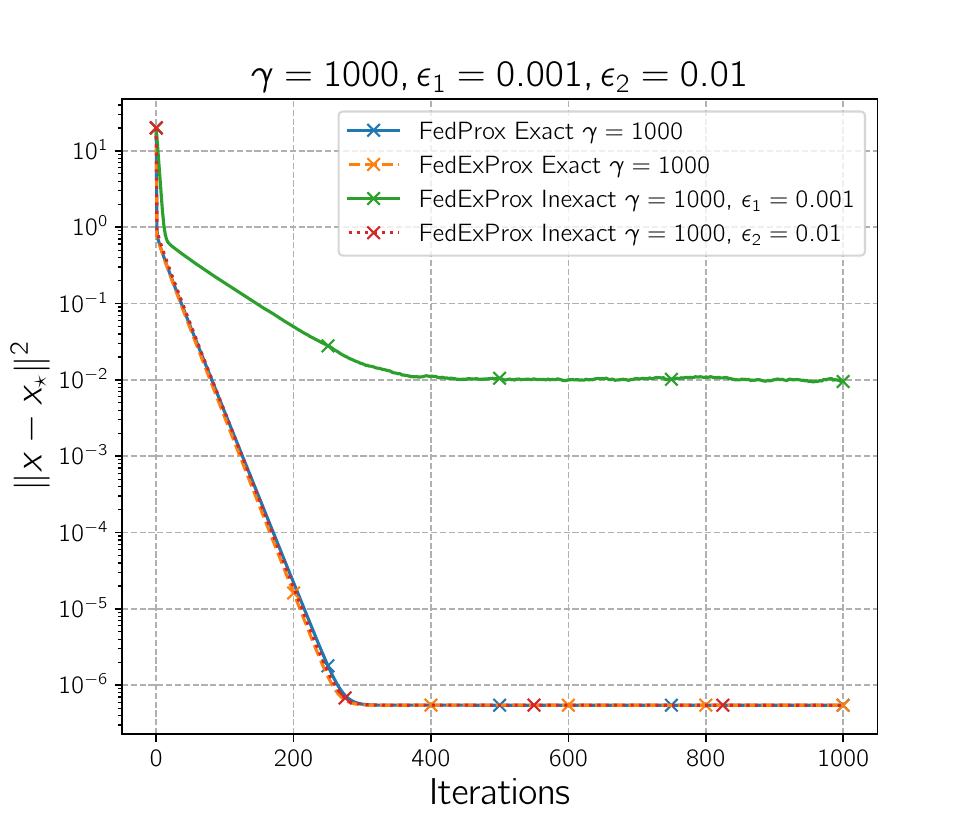}
	   \end{minipage}
   }
    
   \caption{Comparison of {\FEDPROX}, {\FEDEXPROX} with exact proximal evaluations, {\FEDEXPROX} with $\varepsilon_1$-absolute approximation and {\FEDEXPROX} with $\varepsilon_2$-relative approximation. 
   In this case, we fix $\varepsilon_1 = 0.001$, $\varepsilon_2 = 0.01$ and pick the local step size $\gamma \in \cbrac{1000, 100, 10, 1, 0.1. 0.01}$.
   The $y$-axis is the squared distance to the minimizer of $f$, and the $x$-axis denotes the iterations. 
   }

   \label{fig:1-1}
\end{figure}

\begin{figure}
	\centering
   \subfigure{
	   \begin{minipage}[t]{0.98\textwidth}
		   \includegraphics[width=0.33\textwidth]{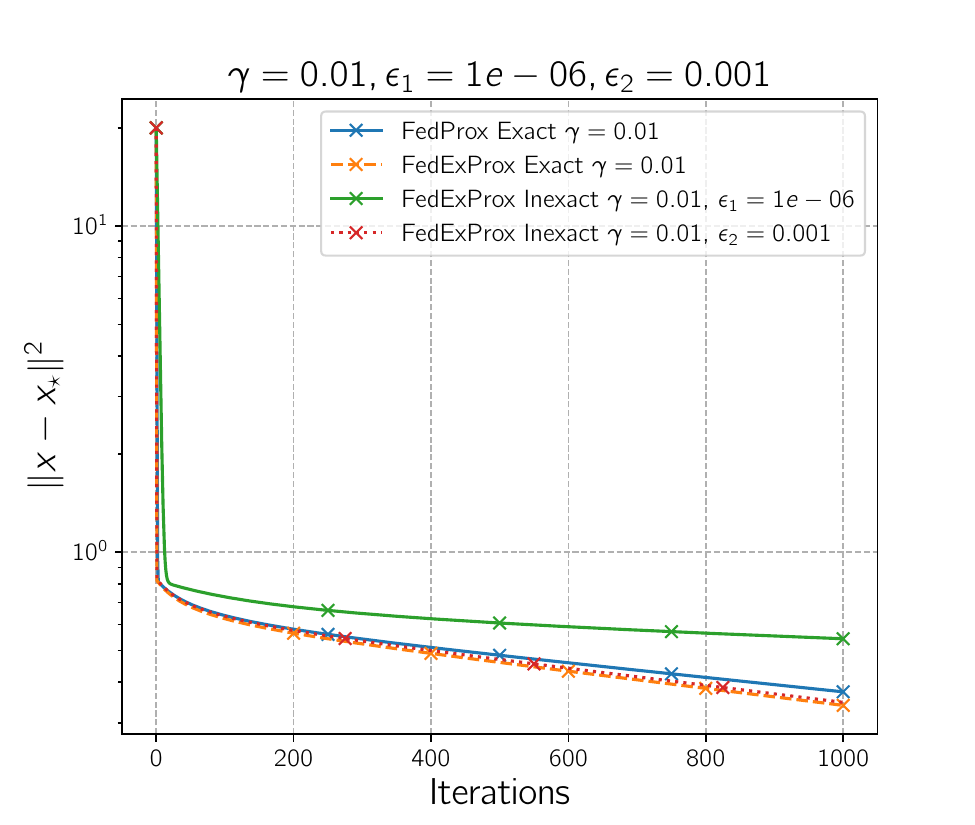} 
		   \includegraphics[width=0.33\textwidth]{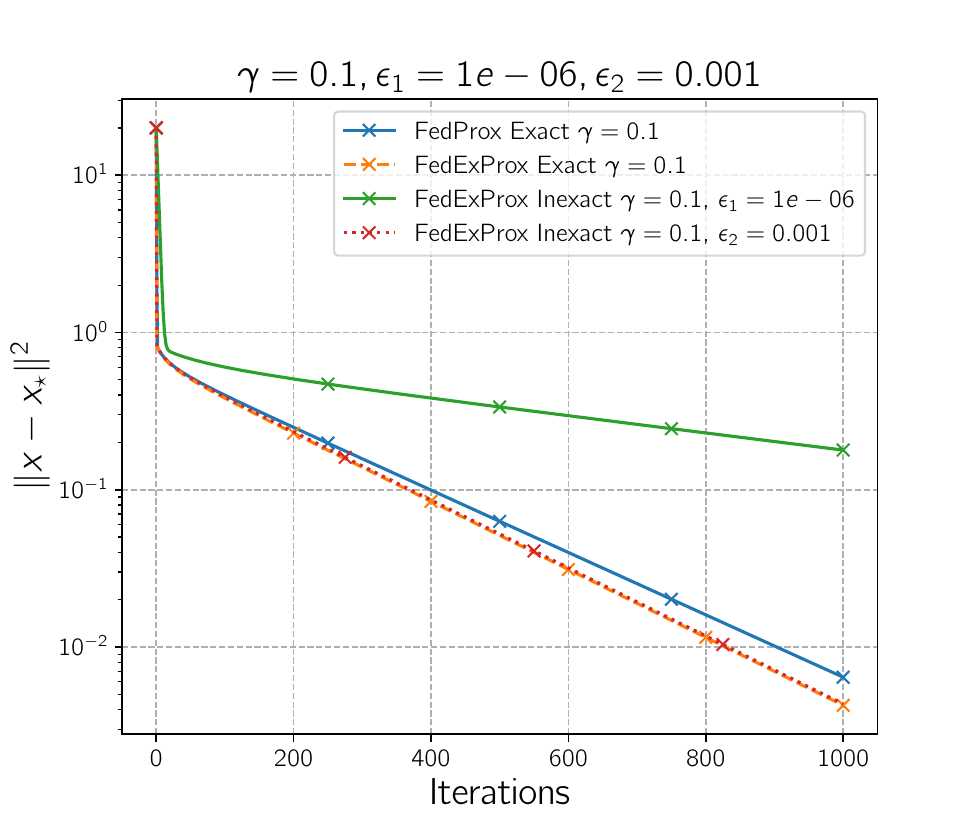}
         \includegraphics[width=0.33\textwidth]{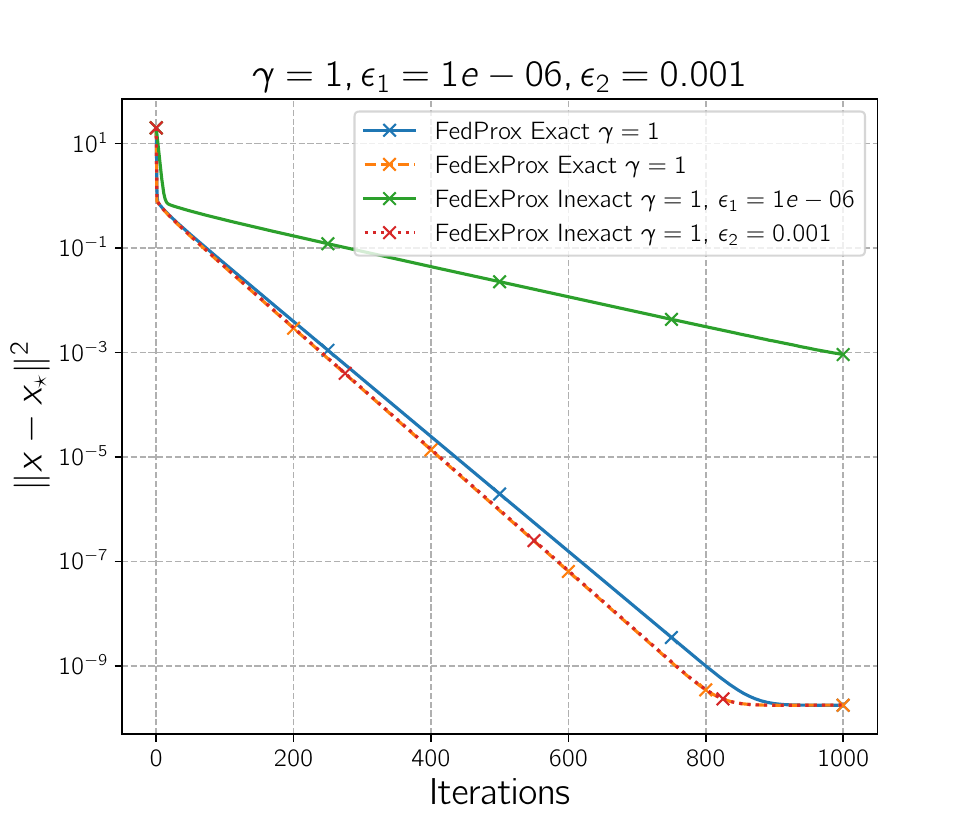}
	   \end{minipage}
   }

   \subfigure{
	   \begin{minipage}[t]{0.98\textwidth}
		   \includegraphics[width=0.33\textwidth]{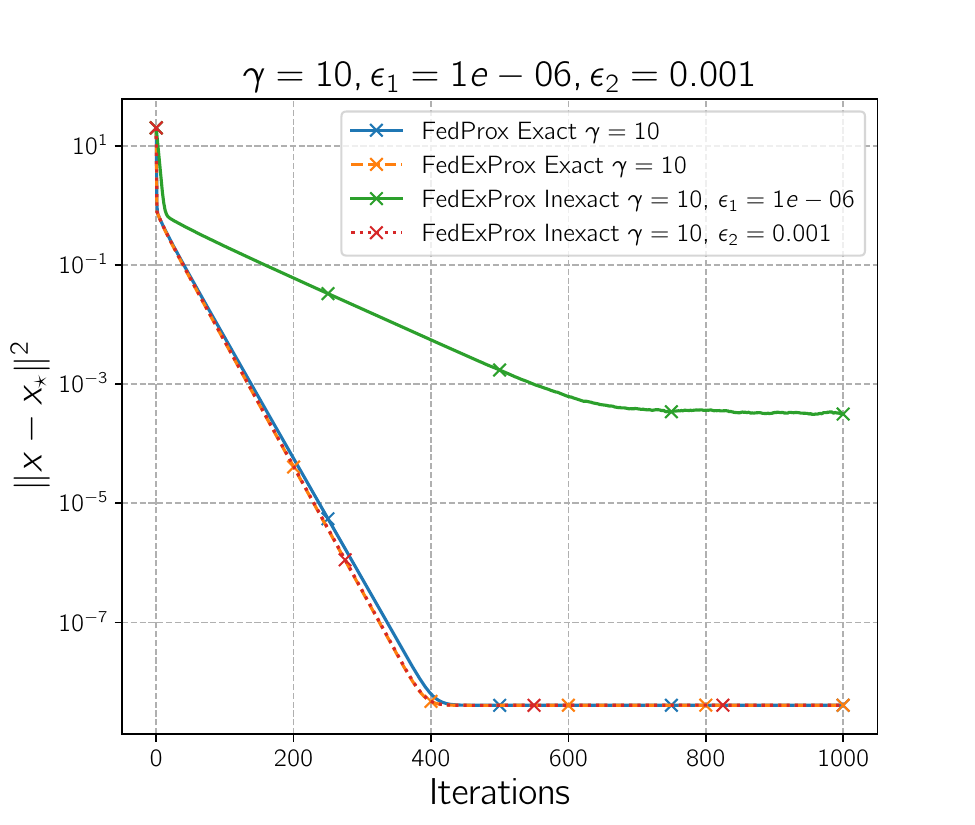} 
		   \includegraphics[width=0.33\textwidth]{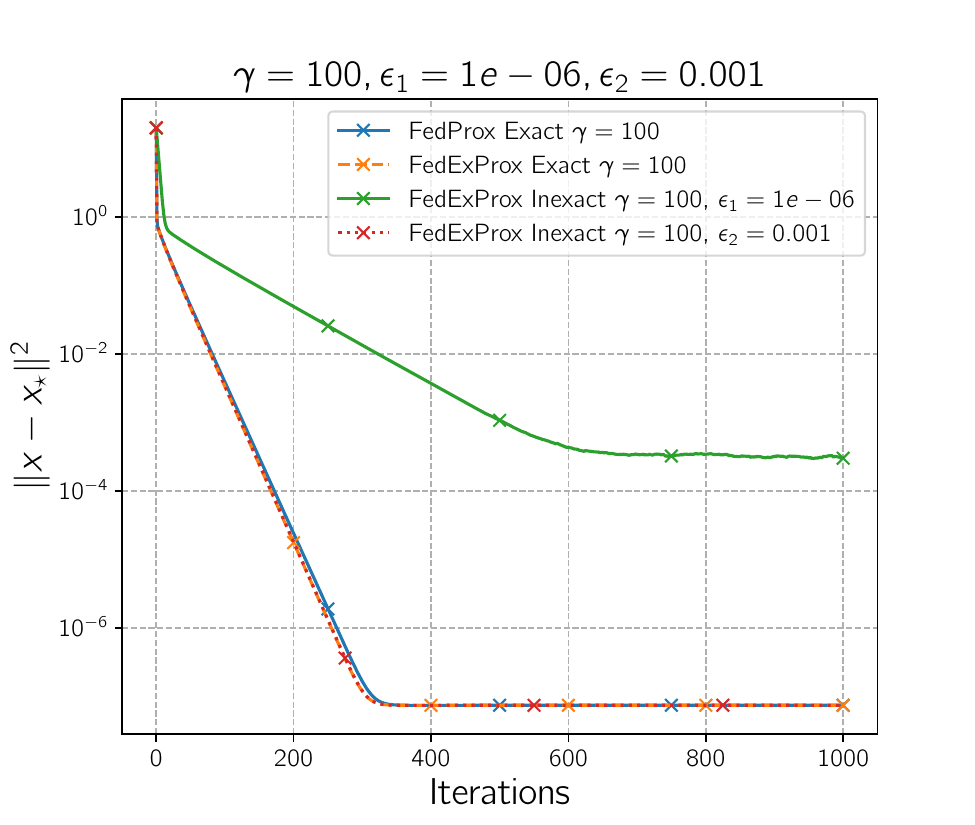}
         \includegraphics[width=0.33\textwidth]{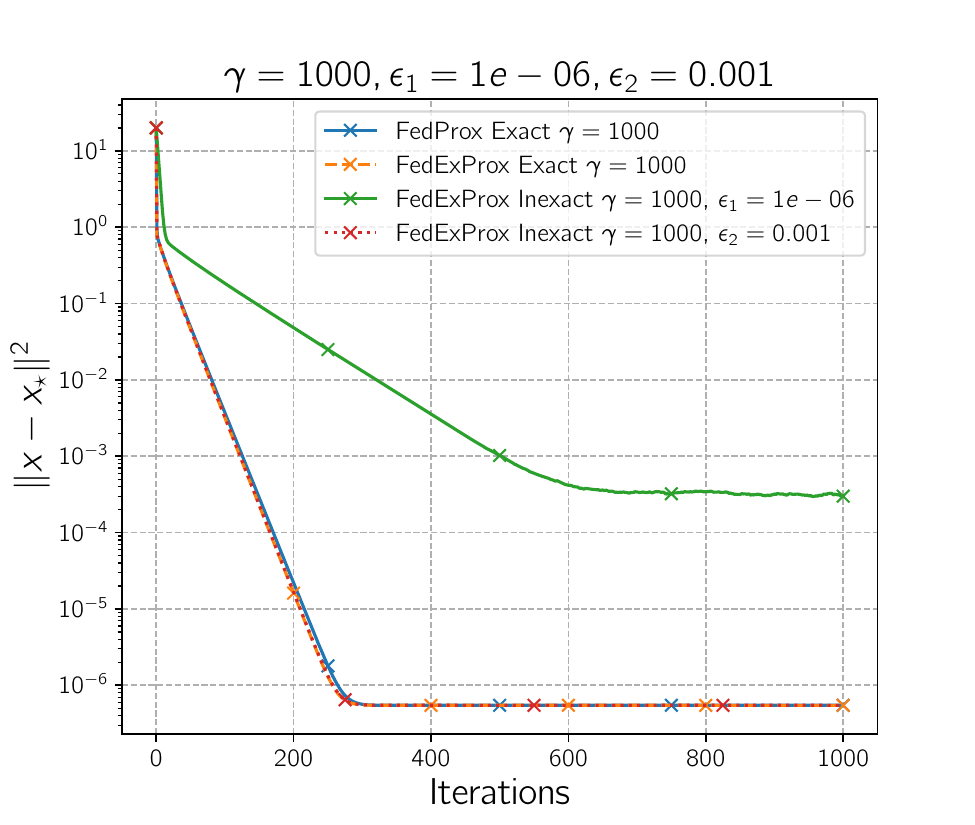}
	   \end{minipage}
   }
    
   \caption{Comparison of {\FEDPROX}, {\FEDEXPROX} with exact proximal evaluations, {\FEDEXPROX} with $\varepsilon_1$-absolute approximation and {\FEDEXPROX} with $\varepsilon_2$-relative approximation. 
   In this case, we fix $\varepsilon_1 = 1e-6$, $\varepsilon_2 = 0.001$ and pick the local step size $\gamma \in \cbrac{1000, 100, 10, 1, 0.1. 0.01}$.
   The $y$-axis is the squared distance to the minimizer of $f$, and the $x$-axis denotes the iterations. 
   }
   
   \label{fig:1-2}
\end{figure}

Notably in \Cref{fig:1-1} and \Cref{fig:1-2}, in all cases, {\FEDEXPROX} with absolute approximation exhibits the poorest performance and converges only to a neighborhood of the solution.
This is expected, since the bias in this case does not go to zero as the algorithm progresses.
It is worth mentioning that as the local step size $\gamma$ increases, the size of the neighborhood decreases, which supports our claim in \Cref{thm:1:conv-full-batch-stncvx}.
As anticipated, in all cases, {\FEDEXPROX} outperforms {\FEDPROX} due to server extrapolation. 
However, as $\gamma$ increases, the performance gap between them diminishes.
The performance of {\FEDEXPROX} with relative approximation is surprisingly good, outperforming {\FEDPROX} in several cases.
This suggests the effectiveness of server extrapolation even when the proximal evaluations are inexact.

\subsection{\texorpdfstring{Comparison of {\FEDEXPROX} with absolute approximation under different inaccuracies}{Comparison of FedExProx with absolute approximation under different inaccuracies}}

In this section, we compare {\FEDEXPROX} with absolute approximations under different level of inaccuracies. 
We fix the local step size $\gamma$ to be a reasonable value, and we vary the level of inexactness for the algorithm.
Specifically, we select $\gamma$ from the set $\cbrac{0.1, 1, 10}$ and for each choice of $\gamma$, we select $\varepsilon_1$ from the set $\cbrac{0.001, 0.005, 0.01, 0.05, 0.1}$.

\begin{figure}
	\centering
   \subfigure{
	   \begin{minipage}[t]{0.98\textwidth}
		   \includegraphics[width=0.33\textwidth]{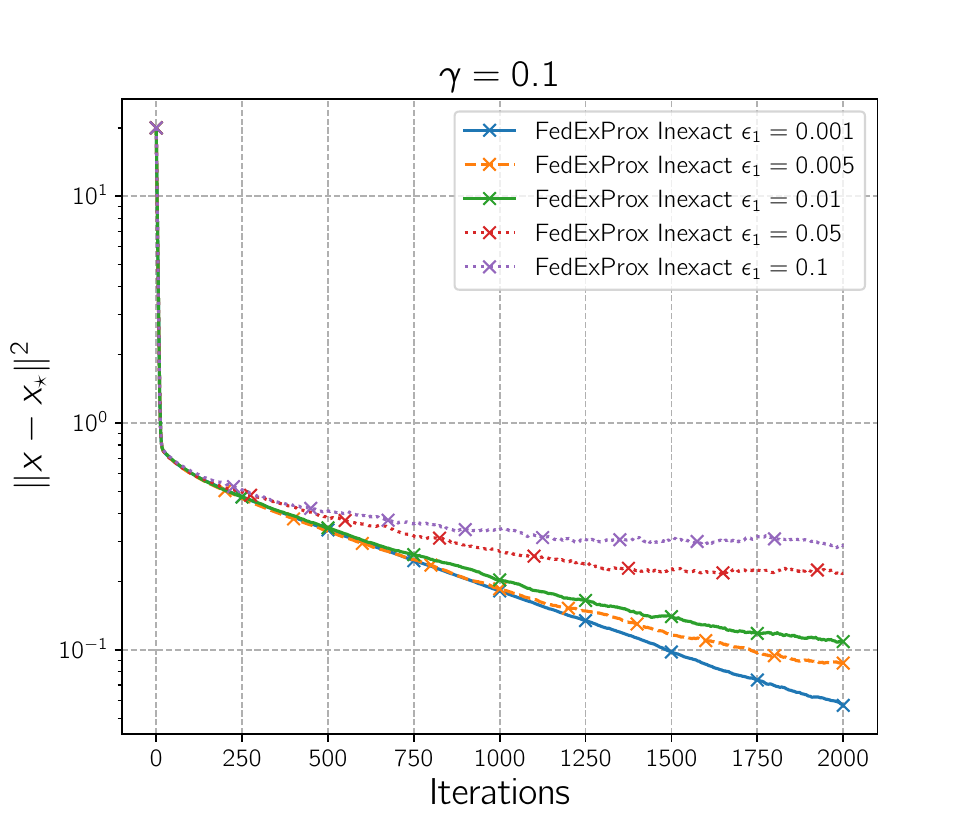} 
		   \includegraphics[width=0.33\textwidth]{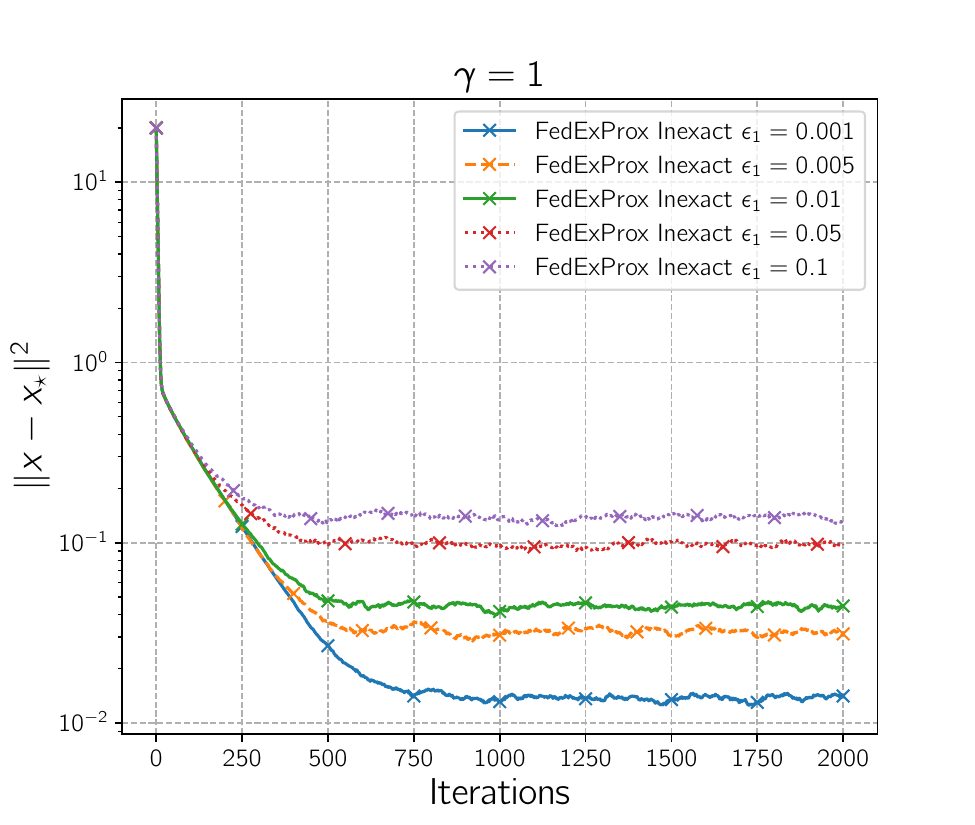}
         \includegraphics[width=0.33\textwidth]{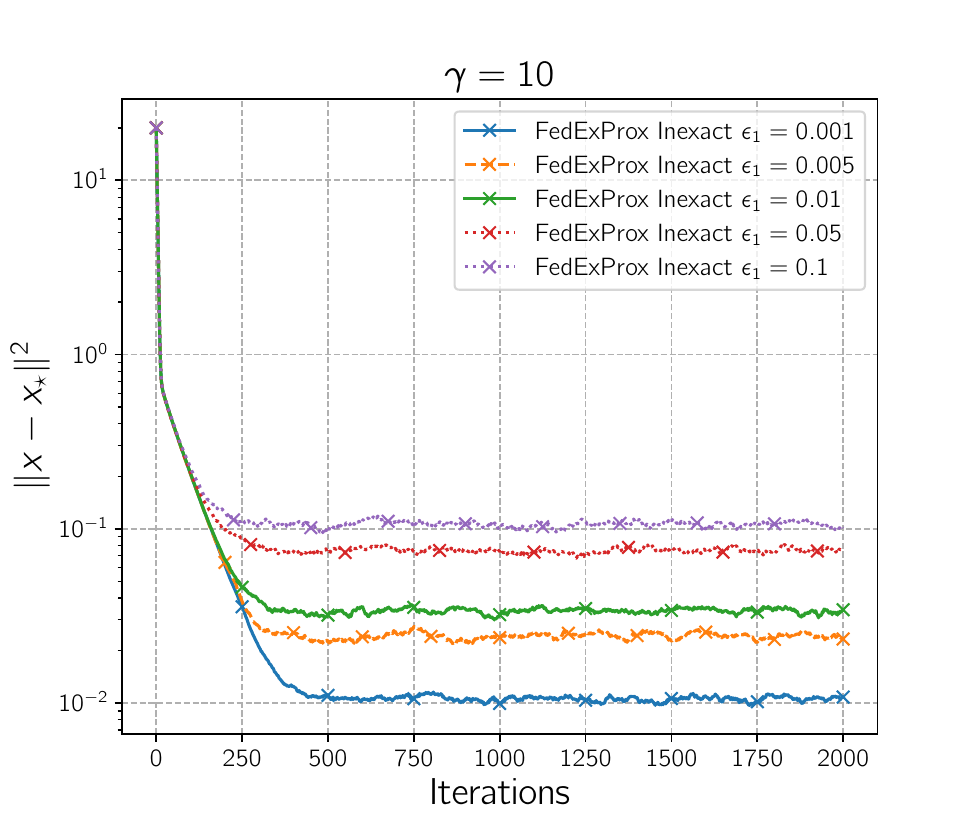}
	   \end{minipage}
   }
    
   \caption{Comparison of {\FEDEXPROX} with $\varepsilon_1$-absolute approximation under different level of inexactness.
   We select $\gamma$ from the set $\cbrac{0.1, 1, 10}$ and for each choice of $\gamma$, we select $\varepsilon_1$ from the set $\cbrac{0.001, 0.005, 0.01, 0.05, 0.1}$.
   The $y$-axis denotes the squared distance to the minimizer and the $x$-axis is the number of iterations.
   }
   
   \label{fig:2-1}
\end{figure}

As observed in \Cref{fig:2-1}, the size of the neighborhood increases with $\varepsilon_1$, further corroborating our theoretical findings in \Cref{thm:1:conv-full-batch-stncvx}. 
Before reaching the neighborhood, the convergence rates of {\FEDEXPROX} with different level of inexactness are similar, which is expected.

\subsection{\texorpdfstring{Comparison of {\FEDEXPROX} with relative approximation under different inaccuracies}{Comparison of FedExProx with relative approximation under different inaccuracies}}

In this section, we compare {\FEDEXPROX} with relative approximations under different level of relative inaccuracies. 
We fix the local step size $\gamma$ to be a reasonable value, and we vary the level of inexactness for the algorithm.
Specifically, we select $\gamma$ from the set $\cbrac{0.1, 0.05, 0.01}$ and for each choice of $\gamma$, we select $\varepsilon_2$ from the set $\cbrac{0.001, 0.005, 0.01, 0.05, 0.1}$.

\begin{figure}
	\centering
   \subfigure{
	   \begin{minipage}[t]{0.98\textwidth}
		   \includegraphics[width=0.33\textwidth]{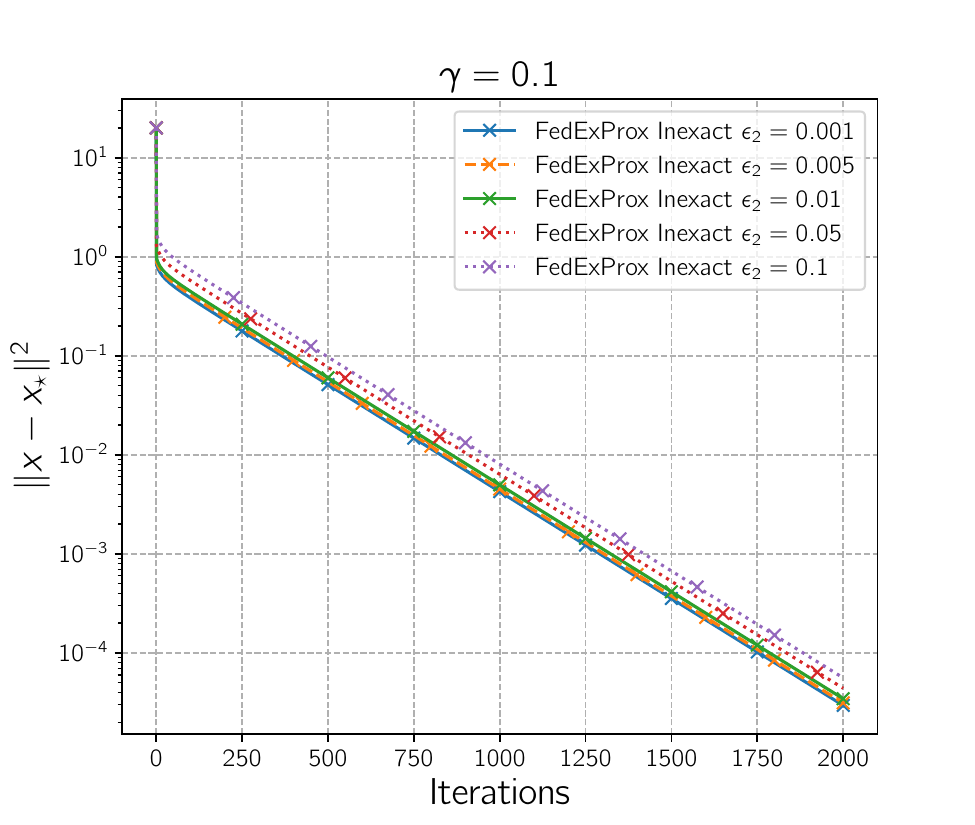} 
		   \includegraphics[width=0.33\textwidth]{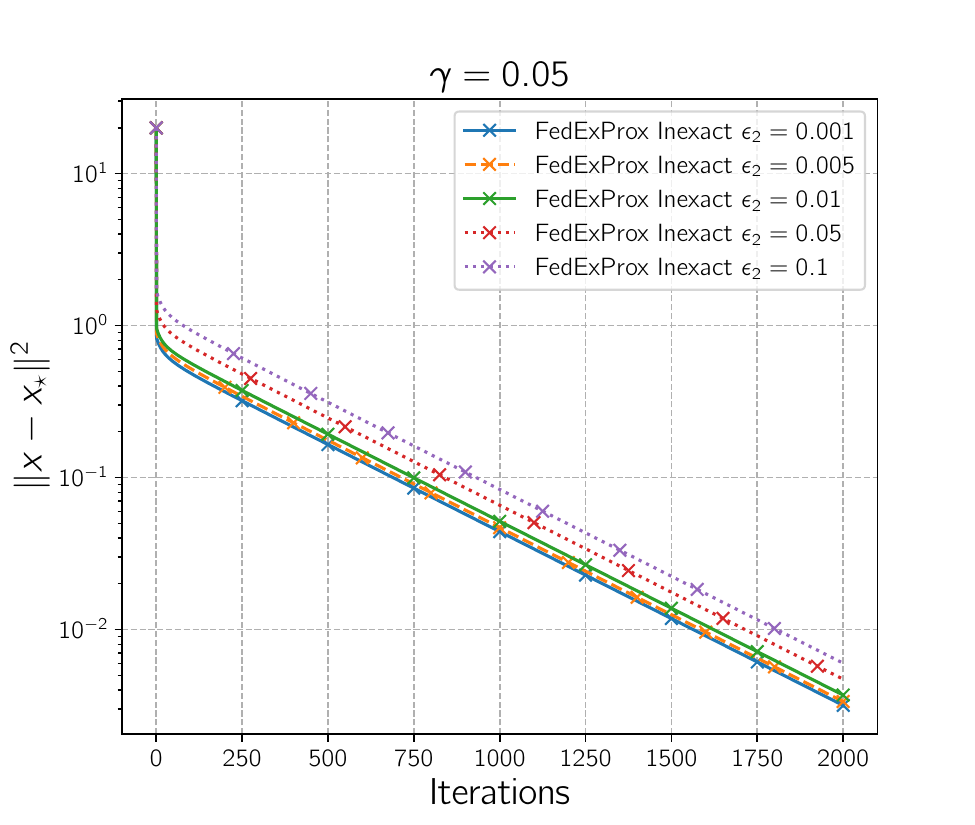}
         \includegraphics[width=0.33\textwidth]{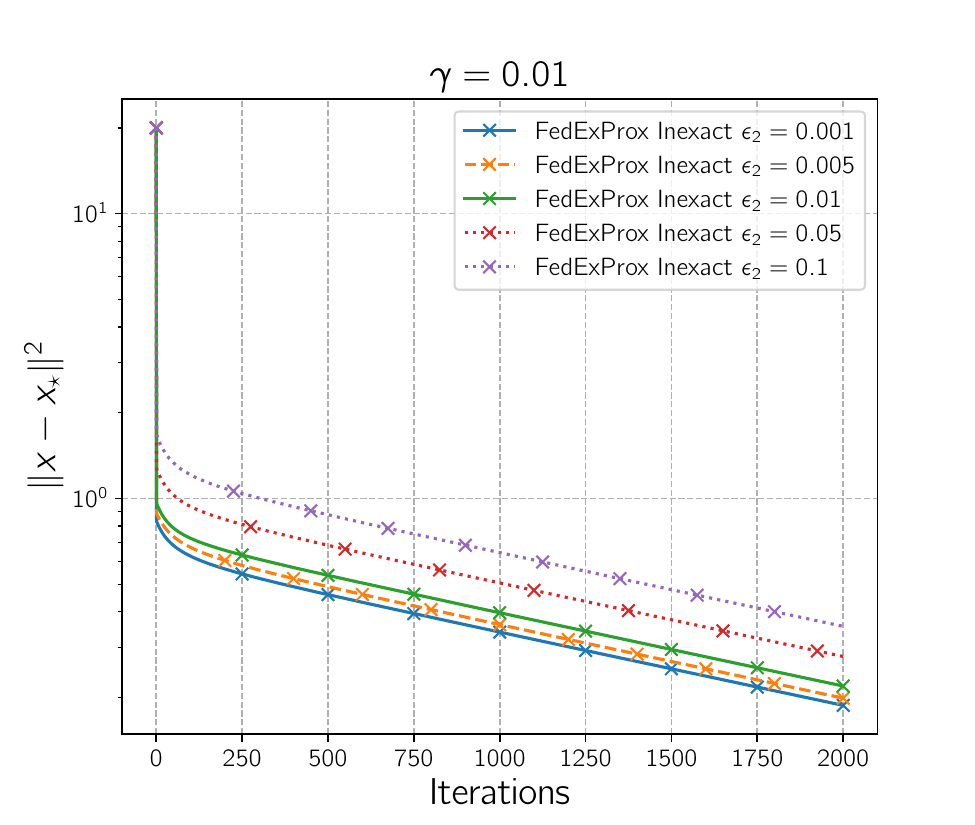}
	   \end{minipage}
   }
    
   \caption{Comparison of {\FEDEXPROX} with $\varepsilon_2$-relative approximation under different level of inexactness.
   We select $\gamma$ from the set $\cbrac{0.01, 0.05, 0.1}$ and for each choice of $\gamma$, we select $\varepsilon_2$ from the set $\cbrac{0.001, 0.005, 0.01, 0.05, 0.1}$.
   The $y$-axis denotes the squared distance to the minimizer and the $x$-axis is the number of iterations.
   }
   
   \label{fig:3-1}
\end{figure}

As observed in \Cref{fig:3-1}, in all cases, a smaller $\varepsilon_2$ corresponds to faster convergence of the algorithm. 
This supports the claim of \Cref{thm:010-stncvx-rela}.
All the tested algorithm converges to the exact solution linearly, which validates the effectiveness of the proposed technique of relative approximation to reduce the bias term.

\subsection{Adaptive extrapolation for inexact proximal evaluations}

In this section, we study the possibility of applying adaptive extrapolation to {\FEDEXPROX} with relative approximations.
We do not consider the case of absolute approximation since it converges only to a neighborhood, which causes problems when combined with adaptive step sizes such as gradient diversity and Polyak step size.

We are using the following definition of gradient diversity based extrapolation,
\begin{align*}
   \alpha_k = \alpha_{k, G} \eqdef \frac{1 + \gamma L_{\max}}{\gamma L_{\max}}\cdot\frac{\frac{1}{n}\sum_{i=1}^{n}\norm{x_k - \ProxSub{\gamma f_i}{x_k}}^2}{\norm{\frac{1}{n}\sum_{i=1}^{n}\rbrac{x_k - \ProxSub{\gamma f_i}{x_k}}}^2}.
\end{align*}
for Polyak type extrapolation, we use 
\begin{align*}
   \alpha_k = \alpha_{k, S} \eqdef \frac{\frac{1}{n}\sum_{i=1}^n\rbrac{\MoreauSub{\gamma}{f_i}{x_k} - \inf \Moreau^{\gamma}_{f_i}}}{\gamma\norm{\frac{1}{n}\sum_{i=1}^n\nabla \MoreauSub{\gamma}{f_i}{x_k}}^2}.
\end{align*}

\begin{figure}
	\centering
   \subfigure{
	   \begin{minipage}[t]{0.98\textwidth}
		   \includegraphics[width=0.33\textwidth]{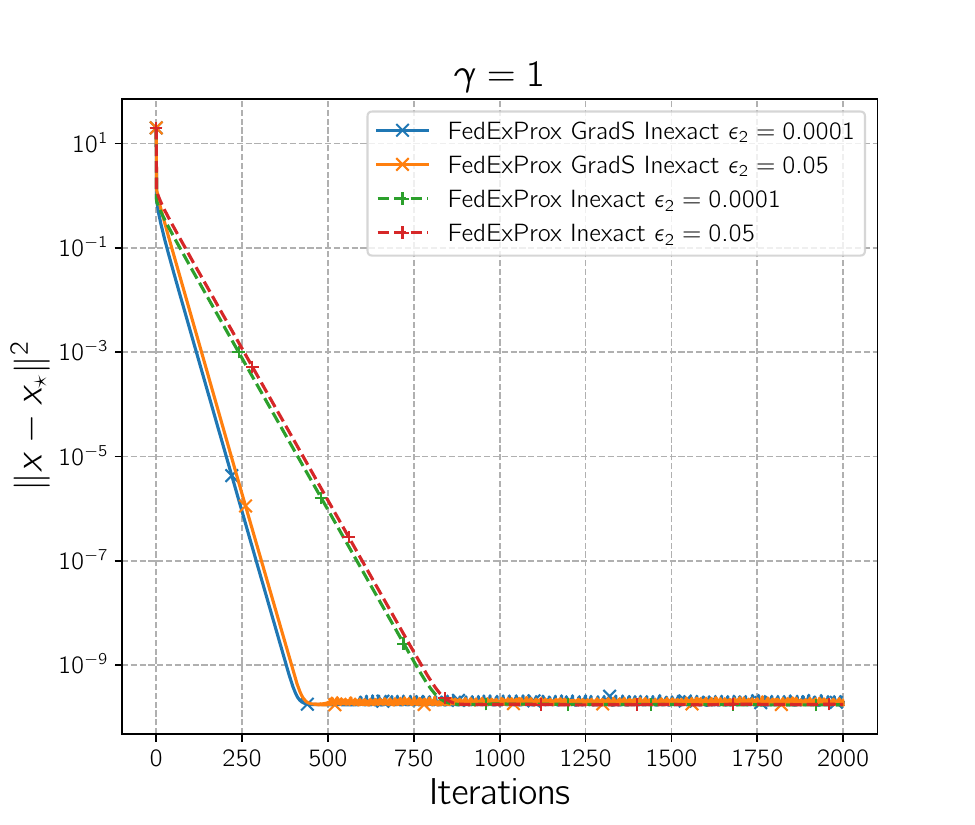} 
		   \includegraphics[width=0.33\textwidth]{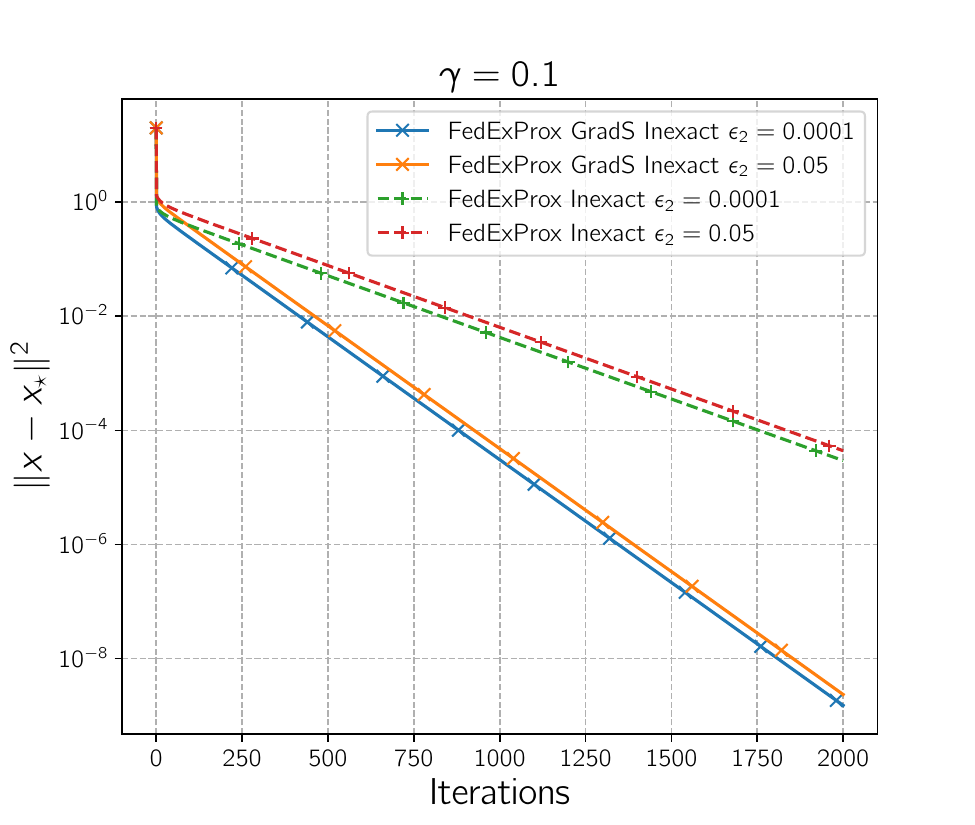}
         \includegraphics[width=0.33\textwidth]{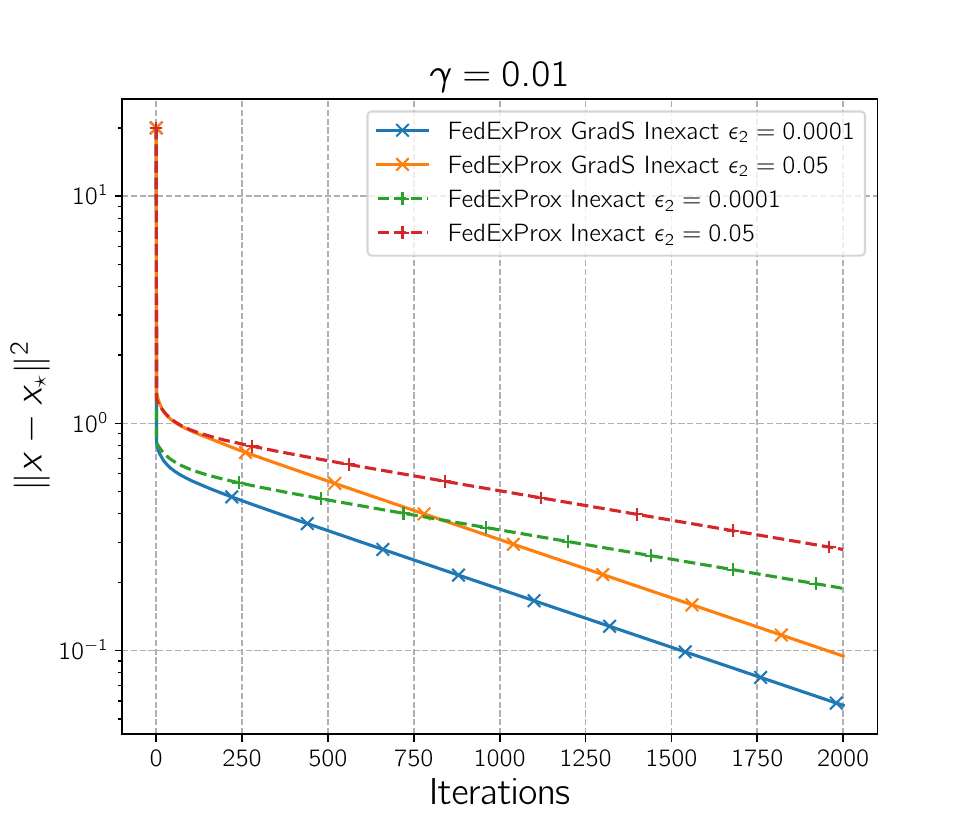}
	   \end{minipage}
   }
    
   \caption{Comparison of {\FEDEXPROX} with $\varepsilon_2$-relative approximation under different level of inexactness using gradient diversity based extrapolation.
   we select $\gamma$ from the set $\cbrac{1, 0.1, 0.01}$ and for each choice of $\gamma$, we select $\varepsilon_2$ from the set $\cbrac{0.0001, 0.05}$.
   The $y$-axis denotes the squared distance to the minimizer and the $x$-axis is the number of iterations.
   }
   
   \label{fig:4-1}
\end{figure}

\begin{figure}
	\centering
   \subfigure{
	   \begin{minipage}[t]{0.98\textwidth}
		   \includegraphics[width=0.33\textwidth]{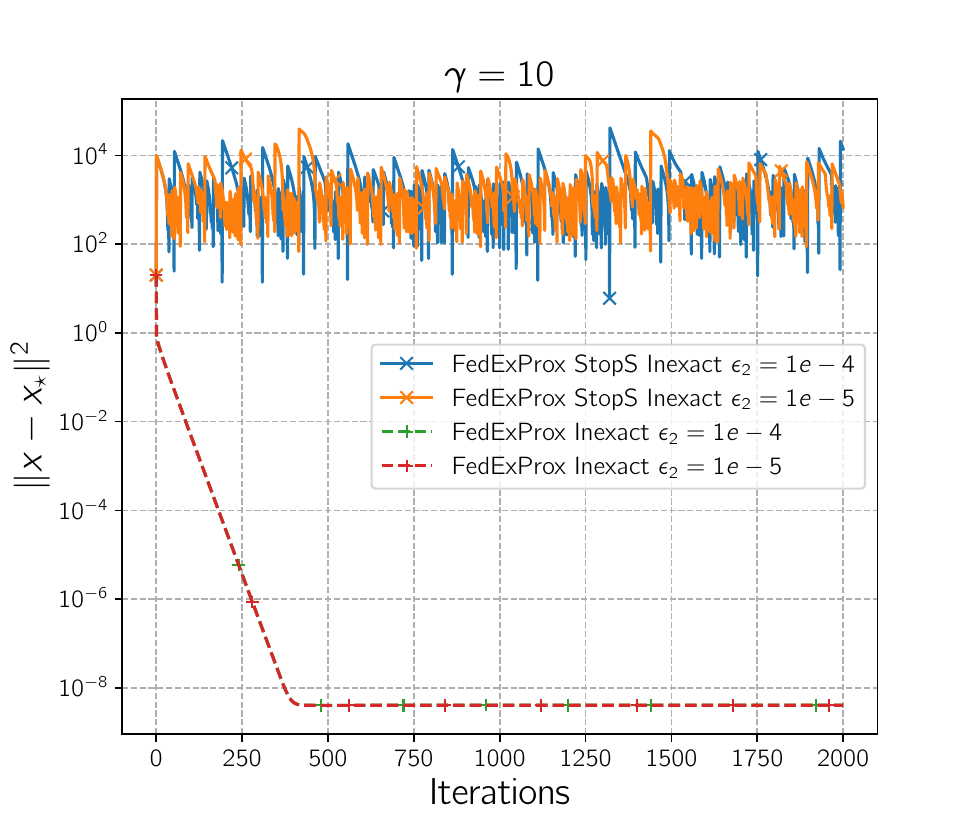} 
		   \includegraphics[width=0.33\textwidth]{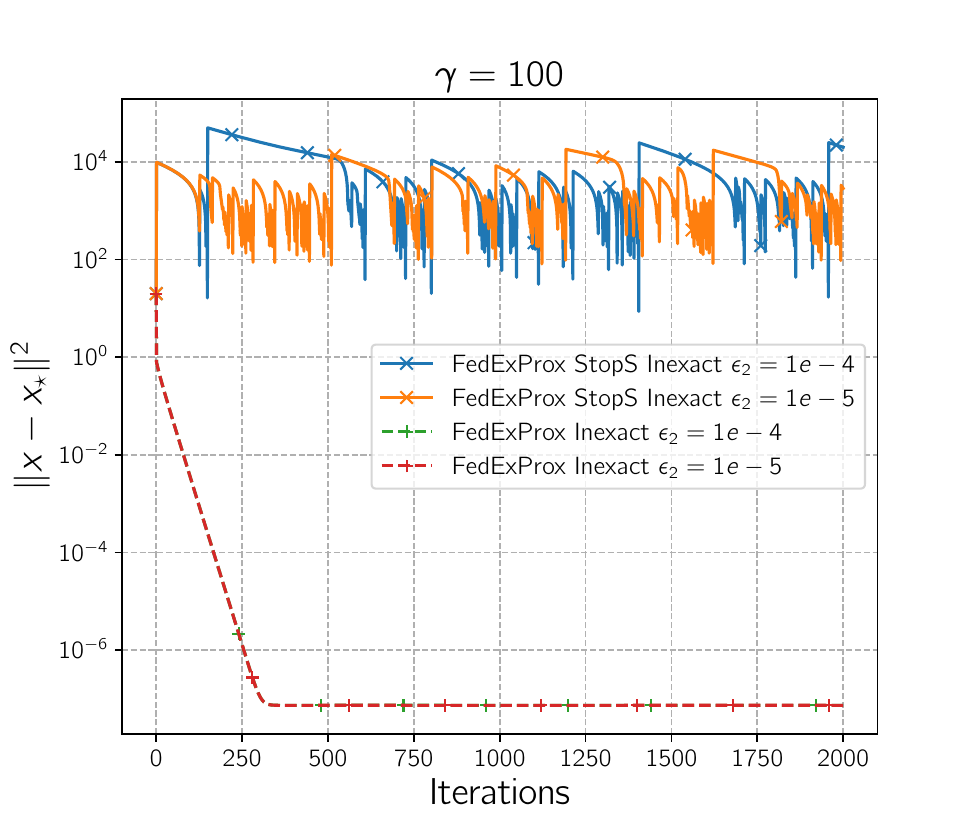}
         \includegraphics[width=0.33\textwidth]{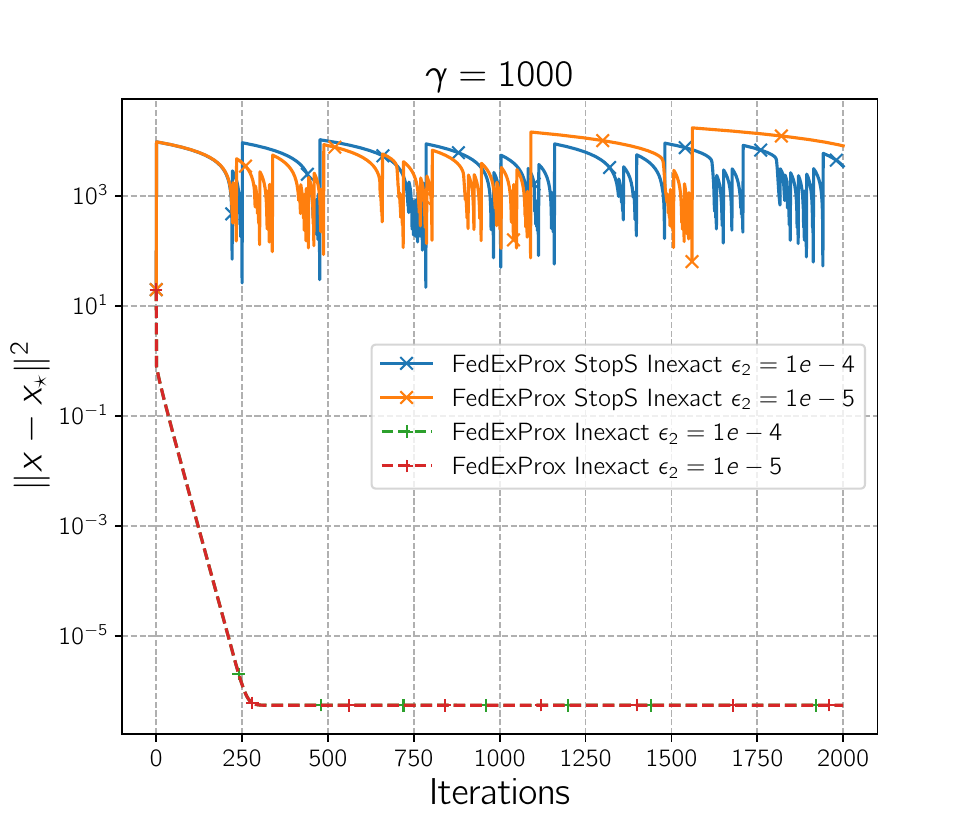}
	   \end{minipage}
   }
    
   \caption{Comparison of {\FEDEXPROX} with $\varepsilon_2$-relative approximation under different level of inexactness using Polyak step size based extrapolation.
   we select $\gamma$ from the set $\cbrac{10, 100, 1000}$ and for each choice of $\gamma$, we select $\varepsilon_2$ from the set $\cbrac{1e-4, 1e-5}$.
   The $y$-axis denotes the squared distance to the minimizer and the $x$-axis is the number of iterations.
   }
   
   \label{fig:4-2}
\end{figure}

As it can be observed from \Cref{fig:4-1}, in all cases, the use of a gradient diversity based adaptive extrapolation results in faster convergence of the algorithm.
This suggests the possibility of developing an adaptive extrapolation for our methods.
However, as we can see from \Cref{fig:4-2}, a direct implementation of Polyak step size type extrapolation results in divergence of the algorithm, indicating that the challenge may be more complex than anticipated.
In our case, this is equivalent to designing adaptive step sizes for {\SGD} with biased updates or {\CGD} with biased compression.
To the best of our knowledge, this field remains open and requires further investigation, as biased updates are quite common in practice.

\end{document}